\documentclass[11pt,a4paper]{article}
\usepackage[T1]{fontenc}
\usepackage[utf8]{inputenc}
\usepackage[margin=1in]{geometry}
\usepackage{tabularx, enumerate,enumitem, authblk,soul}
\usepackage[dvipsnames]{xcolor}
\usepackage{amsmath,amssymb,dsfont,subcaption}
\usepackage{mathtools, graphicx,color,dsfont,amsthm,bbm, mathrsfs}
\usepackage{lipsum,bm,comment}
\usepackage{natbib}
\usepackage{longtable}
\usepackage{array}
\usepackage{caption}
\usepackage[colorlinks=true,linkcolor=blue,anchorcolor=blue,citecolor=blue,filecolor=black,menucolor=black,runcolor=black,urlcolor=black]{hyperref}
\usepackage{cleveref}

\newcommand{\dH}{\dot{H}}
\newcommand{\divH}{\dot{H}_{\diamond}}

\theoremstyle{plain}
\newtheorem{Theor}{Theorem}[section]
\newtheorem{Lem}[Theor]{Lemma}

\newtheorem{Prop}[Theor]{Proposition}

\newtheorem{Condition}[Theor]{Condition}

\theoremstyle{definition}

\newtheorem{Example}[Theor]{Example}

\newtheorem{Remark}[Theor]{Remark}

\newcommand{\A}{\mathbb{A}}
\newcommand{\B}{\mathbb{B}}

\newcommand{\E}{\mathbb{E}}

\newcommand{\I}{\mathbb{I}}

\newcommand{\N}{\mathbb{N}}

\renewcommand{\P}{\mathbb{P}}

\newcommand{\R}{\mathbb{R}}

\newcommand{\T}{\mathbb{T}}

\newcommand{\X}{\mathbb{X}}

\newcommand{\Z}{\mathbb{Z}}

\newcommand{\BB}{\mathcal{B}}

\newcommand{\DD}{\mathcal{D}}

\newcommand{\HH}{\mathcal{H}}
\newcommand{\II}{\mathcal{I}}
\newcommand{\JJ}{\mathcal{J}}

\newcommand{\LL}{\mathcal{L}}

\newcommand{\NN}{\mathcal{N}}

\renewcommand{\SS}{\mathcal{S}}

\newcommand{\VV}{\mathcal{V}}

\newcommand{\XX}{\mathcal{X}}

\newcommand{\ZZ}{\mathcal{Z}}

\newcommand{\CCC}{\mathscr{C}}

\newcommand{\EEE}{\mathscr{E}}

\newcommand{\GGG}{\mathscr{G}}

\newcommand{\RRR}{\mathscr{R}}

\def\1{{\mathbbm{1}}}

\newcommand{\ie}{\textit{i.e.} }

\newcommand{\ps}[2]{\left\langle #1,#2 \right\rangle}
\newcommand{\sm}{\setminus}

\newcommand{\ve}{\varepsilon}
\newcommand{\Om}{\Omega}

\newcommand{\spa}{\quad\quad}

\let\emptyset\varnothing

\newcommand{\embed}{\hookrightarrow}

\usepackage{scalerel,stackengine}
\stackMath
\newcommand\wwidehat[1]{%
	\savestack{\tmpbox}{\stretchto{%
			\scaleto{%
				\scalerel*[\widthof{\ensuremath{#1}}]{\kern-.6pt\bigwedge\kern-.6pt}%
				{\rule[-\textheight/2]{1ex}{\textheight}}
			}{\textheight}%
		}{0.5ex}}%
	\stackon[1pt]{#1}{\tmpbox}%
}

\newcommand{\RN}[1]{%
	(\textup{\uppercase\expandafter{\romannumeral#1}})%
}
\newcolumntype{H}{>{\setbox0=\hbox\bgroup}c<{\egroup}@{}}

\newenvironment{Proof}[1]
{
	\begin{proof}[Proof of #1]
	}
	{
	\end{proof}
}

\allowdisplaybreaks

\title{Inverting the Fisher information operator in non-linear models}
\author{Dimitri Konen}
\affil{University of Cambridge}

\date{}

\begin{document}
	
\maketitle
	
\vspace{-10mm}
\begin{abstract}
We consider non-linear regression models corrupted by generic noise when the regression functions form a non-linear subspace of $L^2$, relevant in non-linear PDE inverse problems and data assimilation. We show that when the score of the model is injective, the Fisher information operator is automatically invertible between well-identified Hilbert spaces, and we provide an operational characterization of these spaces. This allows us to construct in broad generality the efficient Gaussian involved in the classical minimax and convolution theorems to establish information lower bounds, that are typically achieved by Bayesian algorithms as in \cite{Nickl2024} thus showing optimality of these methods. We illustrate our results on time-evolution PDE models for reaction-diffusion and Navier-Stokes equations.

\end{abstract}

\setcounter{tocdepth}{2}
{
    \hypersetup{linkcolor=black}
    \tableofcontents
}

\section{Introduction}\label{sec:Intro}

\subsection{Information geometry of non-linear models}\label{sec:IntroInfoGeo}

We consider i.i.d. data $(X_i, Y_i)_{i=1}^N$ obtained through the regression model 
\begin{equation}\label{model}
Y_i 
=
\GGG(\theta_0)(X_i) 
+
\ve_i
,\spa 
i=1,2,\ldots, N
,
\end{equation}
where $\GGG:\Theta\mapsto L^2_\lambda$ is a map between a Borel measurable subset  $\Theta\ni \theta_0$ of a separable normed vector space $(\VV, \|\cdot\|_\VV)$ and $$
L^2_\lambda 
\equiv 
L^2(\XX, \R^p, d\lambda)
,
$$
for a measurable set $\XX\subset \R^{d_\XX}$, an integer $p\ge 1$, and a probability measure $\lambda$ on $\XX$. We assume that the $X_i$'s have been independently sampled from $\lambda$ and that the errors $\ve_i$ are independent of $X_i$ and i.i.d. samples of a centred distribution on $\R^p$ with density $q_\ve$ such that $\sqrt{q_\ve}$ belongs to the Sobolev space $H^1(\R^p)$.
For $P_{\theta}$ the joint law on $\XX\times \R^p$ of a generic couple $(X, Y)$ with $X\sim \lambda$ and 
$
Y = \GGG(\theta)(X)+\ve
,
$
where $\ve\sim q_\ve$ is independent of $X$, then $P_\theta$ admits the density 
\begin{equation}\label{modelDensity}
p_{\theta}(x,y) 
\equiv
\frac{dP_\theta}{d\mu}(x,y)
=
q_\ve\big(y-\GGG(\theta)(x)\big)
,\spa 
x\in\XX,\ y\in\R^p 
,
\end{equation}
with respect to the measure $d\mu=d\lambda\otimes dy$, where~$dy$ stands for the Lebesgue measure on~$\R^p$. We denote by $P_{\theta}^N$ the $N$-fold product measure of $P_{\theta}$, corresponding to the law of the data $(X_i, Y_i)_{i=1}^N$, and we write the corresponding expectation as  $E^N_\theta[\cdot]$. This setting covers the classical theory of non-parametric estimation where $\GGG(\theta)=\theta$ for all $\theta\in\Theta$ and $\Theta$ is a linear space of regression functions in $L^2_\lambda$. It also encompasses scenarios where the collection of admissible regression functions does not necessarily form a linear space; this is particularly relevant in PDE models where we know \emph{a priori} that $\GGG(\theta)$ is the solution to a PDE that may have a non-linear dependence on $\theta$ so that the resulting collection $\{\GGG(\theta):\theta\in\Theta\}$ of regression maps is a non-linear `sub-manifold' of $L^2_\lambda$. 
Such models have been studied, for instance, to perform inference in non-linear inverse problems (\citet*{Nic20, VdvEtAl24}) as well as linear ones (\citet*{Trabs15, MonNicPat19, MonNicPat21, DonVan24, MagVdVZan25}), to infer the diffusion function in SPDEs models (\citet*{AltRei21, AltTieWah24, HofRay24, GioWan25, AltGau25, GioRay25}), or to infer the potential driving systems of interacting particles (\citet*{BelPodZho24, NicEtAl25, HeiPod25, PodEtAL25}). The field of non-parametric estimation in (S)PDE models has been very active and found many applications, for instance, in biology (\citet*{AltRei22}, \citet*{NicSei25}),  fluid mechanics and the geophysical sciences (\citet*{CDRS09}, \citet*{NicTit2024}, \cite{KonNic2025}, \citet*{SeizillesEtAl24}). We further refer to \citet*{Vdv2017} and \citet*{NicklEMS} for general overviews on non-parametric estimation methods in these contexts.\smallskip

The model (\ref{model}) is differentiable in quadratic mean  at $\theta_0$ when the `forward map' $\theta\mapsto \GGG(\theta)$ is differentiable in a suitable sense, and thus admits a continuous score operator $\A_{\theta_0}$ with corresponding adjoint $\A_{\theta_0}^*$; see Section \ref{sec:QMD} for details.
In this setting, for a smooth functional $F$ with continuous linearization $\dot F$ and corresponding adjoint $\dot F^*$, and assuming that the Fisher information operator $\A_{\theta_0}^* \A_{\theta_0}$ of the model is invertible, then \cite{Nickl2024} established general conditions ensuring the existence of $\sqrt{N}$-consistent estimators of $F(\theta_0)$ converging to a centred Gaussian with covariance  $\dot F (\A_{\theta_0}^* \A_{\theta_0})^{-1}\dot F^*$.
The Hájek-Le Cam convolution theorem (see, e.g., Theorem~8.8 in \cite{Vdv1998}) further establishes that the asymptotic mean square error (MSE) of any `regular' estimator of~$F(\theta_0)$ is lower bounded in terms of the same Gaussian, thus to be thought of as the distributional limit of asymptotically optimal (`efficient') estimators; the minimax theorem (\citet*{VdVWel1996}, Theorem~3.12.5) extends this result to arbitrary estimators in a local asymptotic minimax sense. For linear functionals of the form $F(\theta)=\ps{\psi}{\theta}_\VV$, Theorem 4.1 in \cite{VdV1991} yields that for the asymptotic MSE for estimating $\ps{\psi}{\theta_0}_\VV$ to be finite, it is necessary that $\psi$ belongs to the range $R(\A_{\theta_0}^*)$ of $\A_{\theta_0}^*$;
see also Theorem~3.1.4 in \cite{NicklEMS}.
Interpreting this range condition is not easy and understanding the exact set of $\psi$'s belonging to $R(\A_{\theta_0}^*)$ in specific models is a difficult question. It is showed in Lemma A.3 of \citet*{VdV1991} that $R(\A_{\theta_0}^*)=R( (\A_{\theta_0}^*\A_{\theta_0})^{1/2})$, but this result does not provide an operational way to determine the range in concrete examples.
In addition to understanding the range condition, one of the main challenges in practice to prove an information lower bound or convergence towards the efficient Gaussian is to show that the Fisher information operator $\A_{\theta_0}^* \A_{\theta_0}$ is invertible between explicit spaces. Although injectivity can usually be established on relevant subspaces in many cases, when $\Theta$ is infinite-dimensional this is no longer enough to ensure invertibility of $\A_{\theta_0}^* \A_{\theta_0}$ as injectivity and surjectivity become separate and challenging questions in their own right. This was investigated for concrete PDE models in \citet*{MonNicPat19, Nic20} and \citet*{NicPat2023} for inverse problems involving X-ray transforms, Schr\"odinger and elliptic equations, in \citet*{Nickl2024} and \citet*{KonNic2025} for prediction in  non-linear time-evolution equations, and studied in more generality in \citet*{MonNicPat21b, NicPat2023} and \citet*{NicklEMS}. In each case, establishing surjectivity of $\A_{\theta_0}^* \A_{\theta_0}$ leads to challenging PDE questions that require \emph{ad hoc} techniques on a case-by-case basis; see also Section 3.3.3 and Section 4.2 in \cite{NicklEMS}.

\subsection{Main contributions}\label{sec:Results}

We show that, when the forward regression map $\theta\mapsto \GGG(\theta)$ admits a linearization $\I_{\theta_0}$ at $\theta_0$ along a linear subspace $\VV_0$ of $\VV$, the injectivity and range properties of the score $\A_{\theta_0}$ and Fisher information $\A_{\theta_0}^* \A_{\theta_0}$ operators depend \emph{exclusively} on the linearization $\I_{\theta_0}$ and, in fact, only on its natural (in a certain sense) domain of definition. More precisely, we show that the \emph{qualitative} assumption that~$\I_{\theta_0}$ is injective on $\VV_0$ is enough to fully characterize the range of $\A_{\theta_0}^*$ and $\A_{\theta_0}^* \A_{\theta_0}$ and to establish that the latter is a homeomorphism between well-identified spaces. Our results further provide a more operational way to determine these spaces in practice. We do so in quite broad generality under weak assumptions, namely: the probability measure~$\lambda$ from which the design~$(X_i)_{i=1}^N$ is sampled can be arbitrary, the forward regression map $\theta\mapsto \GGG(\theta)$ is assumed to admit an $L^2_\lambda$-linearization $\I_{\theta_0}$ at $\theta_0$ only in each fixed direction belonging to some tangent space $\VV_0\subset \VV$ at $\theta_0$, and the density $q_\ve$ of the common law of the errors $\ve_i$ in~(\ref{model}) can be any density such that $\sqrt{q_\ve}$ belongs to the Sobolev space $H^1(\R^p)$.
To streamline the exposition and highlight the key ideas, we restrict to the case where $\VV$ is a separable Hilbert space, although this is not necessary as our results apply to arbitrary separable normed vector spaces $(\VV, \|\cdot\|_\VV)$ (see Appendix~\ref{sec:InvertibilityBanach}). 
We start by requiring that the forward regression map $\theta\mapsto \GGG(\theta)$ be directionally $L^2_\lambda$-differentiable, which is a standard assumption; see, e.g., Condition 3.1.1 in \citet*{NicklEMS} and Condition 1 in \citet*{NicPat2023} for analogous (although slightly stronger) assumptions.


\begin{Condition}\label{CondIntroLinModel}
    Fix $\theta_0\in\Theta$. Assume there exists a linear subspace $\VV_0\subset \VV$ such that for all $h\in \VV_0$ there exists $\delta(h)>0$ for which $\{\theta_0+ s h : |s| < \delta(h)\}\subset \Theta$. Further assume that there exists a linear and continuous operator 
    $$
    \I_{\theta_0}:(\VV_0,\ps{\cdot}{\cdot}_{\VV})\to L^2_\lambda
    $$
    such that for all $h\in \VV_0$ fixed we have
    $$
    \rho_{\theta_0,h}(s)
    \equiv
    \|\GGG(\theta_0+sh)-\GGG(\theta_0)-\I_{\theta_0}[sh]\|_{L^2_\lambda}
    =
    o(|s|)
    ,\spa 
    s\to 0 
    .
    $$
\end{Condition}

 Under Condition \ref{CondIntroLinModel}, $\I_{\theta_0}$ extends to a linear and continuous operator
	$$
	\I_{\theta_0} : (\overline \VV_0, \ps{\cdot}{\cdot}_{\VV}) \to L^2_\lambda
	,
	$$
	where $\overline \VV_0$ stands for the $\|\cdot\|_\VV$-closure of $\VV_0$. This yields  the continuous adjoint of $\I_{\theta_0}$
    $$
    \I_{\theta_0}^*
    :
    L^2_\lambda\to (\overline\VV_0,\ps{\cdot}{\cdot}_{\VV})
    .
    $$
We can then establish the fundamental information-geometric description of the statistical model~(\ref{model}) allowing us to define the classical notions of score and Fisher information operators.

\begin{Prop}\label{IntroQMD}
    Assume that Condition \ref{CondIntroLinModel} holds and that the density $q_\ve$ of the common law of the errors $\ve_i$ in (\ref{model}) belongs to the Sobolev space $H^1(\R^p)$. Then, the following holds.
\begin{itemize}
    \item[(i)] The model $(P_{\theta} : \theta\in\Theta)$ from~(\ref{modelDensity}) is differentiable in quadratic mean at $\theta_0$ with linear and continuous 
 score operator $\A_{\theta_0} : (\VV_0, \ps{\cdot}{\cdot}_\VV)\to L^2(\XX\times \R^p, dP_{\theta_0})$. 
    The information operator of the model is $\A_{\theta_0}^* \A_{\theta_0}: (\overline \VV_0, \ps{\cdot}{\cdot}_\VV)\to (\overline \VV_0, \ps{\cdot}{\cdot}_\VV)$ and we have
    \begin{equation*}
    \label{IntroInfo}
    \A_{\theta_0}^* \A_{\theta_0}
    =
    \I_{\theta_0}^* \II_\ve \I_{\theta_0}
    ,
    \end{equation*}    
    where $\II_\ve$ is the (symmetric non-negative definite) Fisher information matrix of $q_\ve$
    \begin{equation*}
    \label{IntroIepsilon}
    \II_\ve 
    \equiv 
    4 \int_{\R^p} \big(\nabla \sqrt{q_\ve}\big)(y)\big(\nabla \sqrt{q_\ve}\big)(y)^T\, dy
    =
    4\ \E_\ve\bigg[
    \frac{(\nabla \sqrt{q_\ve})}{\sqrt{q_\ve}} \frac{(\nabla \sqrt{q_\ve})^T}{\sqrt{q_\ve}}
    \bigg]
    .
    \end{equation*}

     \item[(ii)] The Fisher information matrix $\II_\ve$ is invertible.

    \item[(iii)] The sampling model $\{P^N_\theta : \theta\in\Theta\}$ from (\ref{model}) is locally asymptotically normal (LAN) at $\theta_0$ with respect to the tangent space $\VV_0$, and the LAN-norm is given by 
    \begin{equation}\label{IntroLANnorm}
    \|h\|_{\rm LAN}
    \equiv 
    \|\II_\ve^{1/2} \I_{\theta_0}[h]\|_{L^2_\lambda}
    ,\spa 
    \forall\ h\in\VV_0 
    .
    \end{equation}
\end{itemize} 
\end{Prop}

To address the question of inverting the Fisher information operator between well-identified spaces and to determine the range of the score operator, the only non-trivial assumption we will be making is the following.

\begin{Condition}\label{CondIntroInj}
    The operator $\I_{\theta_0} : (\VV_0, \ps{\cdot}{\cdot}_\VV)\to L^2_\lambda$ is injective.
\end{Condition}

Under Condition \ref{CondIntroInj}, the LAN-norm $\|\cdot\|_{\rm LAN}$ from (\ref{IntroLANnorm}) is indeed a \emph{norm} on $\VV_0$ since $\II_\ve$ is non-degenerate by virtue of Proposition \ref{IntroQMD}(iii). Therefore, let $(\HH, \|\cdot\|_\HH)$ denote the $\|\cdot\|_{\rm LAN}$-completion of $\VV_0$, which is a Banach space by construction. One can then extend the continuous operator $\I_{\theta_0} : (\VV_0, \|\cdot\|_\HH)\to L^2_\lambda$ to a continuous operator $\I_{\theta_0} : (\HH, \|\cdot\|_\HH)\to L^2_\lambda$. Since 
\begin{equation}\label{IntroLANH}
\|h\|_\HH
=
\|\II_\ve^{1/2} \I_{\theta_0}[h]\|_{L^2_\lambda}
,\spa 
\forall\ h\in \VV_0 
,
\end{equation}
then this equality also extends to all $h\in\HH$ by approximation and continuity. One can further define the inner-product on $\HH$
$$
\ps{u}{v}_\HH 
\equiv 
\ps{\I_{\theta_0}[u]}{\II_\ve \I_{\theta_0}[v]}_{L^2_\lambda}
,\spa 
u,v\in \HH 
,
$$
which provides $\ps{u}{u}_\HH=\|u\|^2_\HH$, and makes $(\HH, \ps{\cdot}{\cdot}_\HH)$ a separable Hilbert space; see Section~\ref{sec:Extend}. The space $(\HH, \ps{\cdot}{\cdot}_\HH)$ is an extension of the domain of $\I_{\theta_0}$ on which $\I_{\theta_0}$ remains injective by~(\ref{IntroLANH}).

\begin{Theor}\label{IntroHomeo}
    Assume that Condition \ref{CondIntroLinModel} and Condition \ref{CondIntroInj} hold. Then, the Fisher information operator $\I_{\theta_0}^* \II_\ve \I_{\theta_0}$ is a linear homeomorphism between the Hilbert space $(\HH, \ps{\cdot}{\cdot}_\HH)$ and the Banach space $(\SS, \|\cdot\|_\SS)$ defined as the collection of those $h\in\overline\VV_0$ such that $\|h\|_\SS<\infty$, where 
    \begin{equation}\label{IntroS}
    \|h\|_\SS 
    \equiv 
    \sup_{\substack{v\in \VV_0\\ \|v\|_{\HH} \le 1}} \lvert\ps{h}{v}_{\VV}\rvert < \infty
    .
    \end{equation}
    In addition, we have the following range equalities:
    $
    R(\A_{\theta_0}^*) 
    =
    R(\I_{\theta_0}^*) 
    =
    \SS
    .
    $
\end{Theor}

The normed space $(\SS, \|\cdot\|_\SS)$ has the natural characterization as the dual of $\HH$ through the $\SS$-$\HH$ pairing given by the $\overline\VV_0$-inner product in view of (\ref{IntroS}) (Proposition \ref{PropDual}), and we  have the continuous embeddings
\begin{equation*}\label{IntroEmbeddingSVVH}
(\SS, \|\cdot\|_\SS)
\embed 
(\overline \VV_0, \|\cdot\|_\VV) 
\embed 
(\HH, \|\cdot\|_\HH) 
.
\end{equation*}
The main appeal of Theorem \ref{IntroHomeo} is that it establishes that the information operator is \emph{always} invertible between adequate spaces when $\I_{\theta_0}$ is injective on $\VV_0$, and further provides an operational way to actually exhibit the set $\SS$ in many practical situations through (\ref{IntroS}). For instance, one can determine the sets $\HH$ and $\SS$ when, for some norm $\|\cdot\|_X$ on $\VV_0$, we have the two-sided estimates 
$$
\|h\|_X 
\lesssim 
\|\I_{\theta_0}[h]\|_{L^2_\lambda}\lesssim 
\|h\|_X 
,\spa 
\forall\ h\in \VV_0 
.
$$
In this case, Condition \ref{CondIntroInj} is satisfied, $\HH$ coincides with the $\|\cdot\|_X$-completion of $\VV_0$, and $\SS$ coincides with the collection of those $h\in\overline\VV_0$ such that
$$
\sup_{\substack{v\in \VV_0\\ \|v\|_X\le 1}} \lvert\ps{h}{v}_\VV\rvert
<
\infty
.
$$


These results allow us to construct in Section \ref{sec:LimitProcess} the efficient limiting Gaussian process with covariance given by the inverse $(\A_{\theta_0}^* \A_{\theta_0})^{-1}$ of the Fisher information as in classical multivariate statistics, under the sole assumption that~$\I_{\theta_0}$ be injective over the tangent space $\VV_0\subset \VV$. From there, we are able to provide in Section \ref{sec:GeneralMinimax} another (somewhat easier and more intuitive) proof of the necessity of the range condition $\psi\in R(\A_{\theta_0}^*)$ from \citet*{VdV1991} (Theorem 4.1) for finiteness of the local asymptotic minimax MSE for estimating semi-parametric functionals $\ps{\psi}{\theta_0}_{\VV}$. The corresponding asymptotic lower bound from \citet*{VdVWel1996} (Theorem 3.12.5) for general functionals~$F(\theta)$  requires that the efficient limiting Gaussian above exists as a tight random variable in some Banach space, which is another challenge in practice. We use our results to provide concrete conditions on the ambient space and $\I_{\theta_0}$ to guarantee that this random field exists and, hence, that the local asymptotic minimax risk admits a finite lower bound. We further provide the first converse parts to the minimax theorem by showing that existence of this tight process and some continuity for the linearization of $F$ are necessary for the minimax risk to be finite. Finally, in Section~\ref{sec:Appli} we illustrate these results on non-linear data assimilation models in the case of reaction-diffusion and Navier-Stokes equations.

\subsection{Data assimilation in non-linear dynamical systems}

Consider the states $\{u_{\theta}(t):\Om\subset \R^d\to \R^p,\ t\in [0,T]\}$ of a dynamical system parametrized by its initial state $\theta\equiv u_{\theta}(0)$, whose evolution is governed by the
parabolic non-linear equation
\begin{align}\label{eq:Para}
\frac{\partial u}{\partial t}(t,x) + \LL[u](t,x) + N[u](t,x) &= f(t,x)
~~~ \textrm{on } (0,T]\times\Om, \\ \nonumber
u(0,x)&=\theta(x)\spa 
 \textrm {on } \Om
,
\end{align}
where $\LL$ is a linear self-adjoint positive-definite operator with compact inverse (typically, one considers $\LL=-\Delta$, the usual negative Laplacian, or some fractional power or anisotropic version thereof), and $N$ is some non-linear operator specific to the physical phenomenon at hand. Common non-linear systems include reaction-diffusion equations when $N[u](t,x)=F(u(t,x))$ for some reaction vector-field $F:\R^p \to \R^p$ or Navier-Stokes equations obtained when $N[u]=(u\cdot \nabla)u$, on which this paper focuses in  Section \ref{sec:Appli}.
In most practical scenarios, only noisy and discrete samples of the system are available. Thus, we consider data $(t_i,\omega_i, Y_i)_{i=1}^N$ arising from 
$$
Y_i 
=
u_{\theta}(t_i, \omega_i)
+
\ve_i 
,\spa 
i=1,\ldots, N
$$
where $X_i\equiv (t_i, \omega_i)$ are i.i.d. samples drawn from a distribution $\lambda$ on the time-space cylinder $\XX\equiv [0,T]\times \Om$, and $(\ve_i)$ are i.i.d. samples drawn from $q_\ve$ independently from the $X_i$'s. Further assume that $\lambda$ has a density such that $0<\lambda_{\min}\le \lambda(t,x)\le \lambda_{\max}<\infty$ for all $(t,x)\in\XX$.
In this context, an important area of data assimilation---see, e.g., \citet*{E09, LSZ15, RC15} and \citet*{EVvL22}---is concerned with designing statistical procedures that consistently recover the ground truth trajectory from samples corrupted by noise; see also \citet*{Nickl2024} and \citet*{KonNic2025}. 
We illustrate our results in the case of reaction-diffusion equations, covered in \cite{Nickl2024} and detailed in Section~\ref{sec:AppliReac} where we also consider the  incompressible Navier-Stokes equations. For simplicity, we take $\Om \equiv \T^d$ the $d$-dimensional torus with unit length and $d\in\{1,2,3\}$, and consider the unique strong solution $\GGG(\theta)\equiv u_\theta :[0,T]\times \Om\to\R$ of the reaction-diffusion equation
\begin{align*}\label{eq:IntroRD}
\frac{\partial u_\theta}{\partial t}(t,x)
-
\Delta u_\theta(t,x) &= f(u_\theta(t,x))
~~~ \textrm{on } (0,T]\times\Om, \\ \nonumber
u_\theta(0,x)&=\theta(x)\spa 
~~~~~~ \textrm {on } \Om
,
\end{align*}
where $\theta\in H^1(\Om)$ is the initial condition and the reaction term $f$ is taken in $C^\infty_c(\R)$. Then, take $\Theta=H^1(\Om)$ as a subset of the separable Hilbert space $(\VV, \ps{\cdot}{\cdot}_\VV)=(L^2(\Om), \ps{\cdot}{\cdot}_{L^2(\Om)})$. This defines a map $\GGG:\Theta\to L^2_\lambda \equiv L^2([0,T]\times \Om, \R, d\lambda)$. For smooth `ground truth' $\theta_0\in C^\infty(\Om)$ and tangent space $\VV_0 = H^2(\Om)$, then $\GGG:\Theta\to L^2_\lambda$ satisfies Condition~\ref{CondIntroLinModel} at $\theta_0$ with linearization $\I_{\theta_0} : H^2(\Om)\to L^2_\lambda$ given, for all $h\in H^2(\Om)$, by the unique solution $U\equiv \I_{\theta_0}[h]$ to
\begin{align*}
\frac{\partial U}{\partial t}(t,x)
-
\Delta U(t,x) - f'(u_{\theta_0}(t,x)) U(t,x) &= 0\spa\spa\textrm{on } (0,T]\times\Om, \\ 
U(0,x)&=h(x)\spa 
~~ \textrm {on } \Om
.
\end{align*}
The linear operator $\I_{\theta_0}$ is continuous from $(H^2(\Om), \ps{\cdot}{\cdot}_{L^2(\Om)})$ to $L^2_\lambda$, and further satisfies
\begin{equation*}\label{eq:IntroRDEstimate}
\|h\|_{H^1(\Om)^*}
\lesssim 
\|\I_{\theta_0}[h]\|_{L^2([0,T]\times \Om)} 
\lesssim 
\|h\|_{H^1(\Om)^*} 
,\spa 
\forall\ h\in H^2(\Om) 
,
\end{equation*}
where, under the usual representation of $H^1(\Om)^*$ as a space of rough `functions' (distributions),
$$
\|h\|_{H^1(\Om)^*}
=
\sup_{\substack{f\in H^1(\Om)\\ \|f\|_{H^1(\Om)}\le 1}} \Big| \int_{\Om} h(x)f(x)\, dx \Big|
.
$$
In the notation of Theorem \ref{IntroHomeo} above, this provides $\HH=H^1(\Om)^*$ and $\SS=H^1(\Om)$, so that the Fisher information $\A_{\theta_0}^* \A_{\theta_0}$ of the statistical model is a linear homeomorphism from $H^1(\Om)^*$ to~$H^1(\Om)$. The regularity of the model established in Proposition \ref{IntroQMD} above and the results alluded to at the end of Section \ref{sec:Results} permit to establish a lower bound on the local asymptotic minimax risk for estimating a linear functional $\ps{\psi}{\theta_0}_{L^2(\Om)}$, for $\psi\in L^2(\Om)$, defined as
$$
\RRR_{\psi,\theta_0} 
\equiv 
\inf_{(\hat \psi_N)} 
\sup_{\substack{J\subset C^\infty_c(\Om)\\ |J|<\infty}}
\liminf_{N\to \infty} 
\max_{h\in J} 
E^N_{\theta_0+\frac{h}{\sqrt{N}}}
\Big| \sqrt{N}\Big(\hat \psi_N-\big\langle\psi, \theta_0+h/\sqrt{N}\big\rangle_{L^2(\Om)}\Big) \Big|^2
,
$$
where the infimum ranges over all sequences $(\hat \psi_N)$ of estimators of the data $(t_i, \omega_i, Y_i)_{i=1}^N$, and the supremum ranges over finite collections of smooth and compactly supported functions. Then, 
$$
\RRR_{\psi,\theta_0} 
\ge 
\ps{\psi}{(\A_{\theta_0}^*\A_{\theta_0})^{-1}\psi}_{L^2(\Om)}
=
\|\psi\|^2_\SS 
,
$$
if $\psi\in H^1(\Om)=\SS$, whereas $\RRR_{\psi,\theta_0}=\infty$ if $\psi \notin H^1(\Om)$. An analogous lower bound holds for general smooth functionals $F(\theta_0)$ taking values in Banach spaces $\B$ (see Section \ref{sec:MinimaxBound}) thus showing that the Bayesian estimators from \cite{Nickl2024} asymptotically achieve optimal performance in an objective information-theoretic sense.\smallskip

Since all our results only rely on injectivity of $\I_{\theta_0}$, it is clear that one can use the proof technique from \cite{KonNic2025} (more specifically, Proposition~B.2 there) to establish the desired injectivity and, hence, obtain analogous results for general non-linear parabolic equations (whenever they are well-posed) of the form~(\ref{eq:Para}) for non-linearities~$N$ that are compact enough compared to the diffusion operator~$\LL$.

\section{LAN expansion}\label{sec:LAN}


When the errors $\ve_i=Y_i - \GGG(\theta_0)(X_i)$ from (\ref{model}) are normally distributed, one can provide a direct proof of the LAN expansion in Proposition \ref{PropLAN} below as in Theorem 3.1.3 of \citet*{NicklEMS}. In this section, we establish a LAN expansion in the non-Gaussian case by establishing that the model is differentiable in quadratic mean under Condition \ref{CondIntroLinModel} only.

\subsection{Quadratic mean differentiability}\label{sec:QMD}

When the errors are non-Gaussian, a LAN expansion can be established provided the base model $(P_{\theta} : \theta\in\Theta)$ from~(\ref{modelDensity}) is differentiable in quadratic mean, \ie if there exists a `score operator' $\A_{\theta_0} : \VV_0 \to L^2(\XX\times \R^p, P_{\theta_0})$ such that, for all $h\in\VV_0$ fixed we have
\begin{equation}\label{eq:QMD}
\lim_{s\to 0}
\frac{1}{s^2}
\int_{\XX\times \R^p} 
\Big( 
\sqrt{p_{\theta_0+sh}(x,y)} - \sqrt{p_{\theta_0}(x,y)}- \frac{s}{2} \A_{\theta_0}[h](x,y) \sqrt{p_{\theta_0}(x,y)}
\Big)^2\, d\lambda(x)\, dy
=
0
,
\end{equation}
where we recall that $p_{\theta}(x,y)=q_\ve(y-\GGG(\theta)(x))$ for all $\theta\in\Theta$, $x\in\XX$, and $y\in\R^p$. Since quadratic mean differentiability (QMD) is a local statement in nature, after linearization of the parametrization $\theta\mapsto \GGG(\theta)$ of regression functions near $\theta_0$, we expect the model to look like the location model $\{y\mapsto q_\ve(y-z):z\in\R^p\}$ up to a `Jacobian' transformation $\I_{\theta_0}$ arising from the linearization of $\theta\mapsto \GGG(\theta)$. A common difficulty to establish the QMD property is when the densities in the model have different support. In the location model alluded to above with base density $q_\ve$, this occurs if the support of $q_\ve$ is not all of $\R^p$; when $p=1$ and $q_\ve(y)=\1_{[0,1]}(y)$, for instance, the corresponding location model is not QMD. This phenomenon happens because the score function is undetermined on the set $\{q_{\ve}=0\}$ so that the difference $\sqrt{q_{\ve}(y-sh)}-\sqrt{q_{\ve}(y)}$ should be approximated by $0$ when $q_{\ve}(y)=0$. In other words, we should have $(\nabla \sqrt{q_{\ve}})(y)=0$ when $q_{\ve}(y)=0$. This will be the case if $\sqrt{q_{\ve}}$ is smooth, for if $q_{\ve}(y)=0$ then $\sqrt{q_{\ve}}$ attains a minimum at $y$ so that that necessarily $(\nabla \sqrt{q_\ve})(y)=0$. When dealing with a non-smooth $q_\ve$, one should thus argue that any reasonable notion of gradient of $\sqrt{q_\ve}$ has a version that vanishes on the set $\{q_\ve=0\}$. We show in the next Proposition that QMD holds when $\sqrt{q_\ve}$ belongs to the Sobolev space $H^1(\R^p)$, which is a non trivial fact and relies on an auxiliary lemma deferred to Appendix \ref{sec:AppendixSupport} (see Lemma \ref{LemSupport}) establishing that, in this case, the weak gradient $\nabla \sqrt{q_\ve}$ of $\sqrt{q_\ve}$ admits a version that vanishes on $\{q_\ve=0\}$.
\
\begin{Prop}\label{PropQMD}
    Let $\theta_0\in\Theta$ and $\VV_0$ satisfy Condition \ref{CondIntroLinModel}. 
    Assume that the density $q_\ve$ of the common law of the errors $\ve_i$ in (\ref{model}) belongs to the Sobolev space $H^1(\R^p)$.
    Then, the model $(P_{\theta} : \theta\in\Theta)$ from~(\ref{modelDensity}) is differentiable in quadratic mean at $\theta_0$ with 
 score operator
    $$
    \A_{\theta_0}[h](x,y)
    \equiv 
    -2\frac{\ps{(\nabla\sqrt q_\ve)\big(y - \GGG(\theta_0)(x)\big)}{\I_{\theta_0}[h](x)}_{\R^p}}{ 
    \sqrt{p_{\theta_0}(x,y)} 
    }\
    \1\big[p_{\theta_0}(x,y) > 0\big]
    .
    $$
\end{Prop}


\begin{Proof}{Proposition \ref{PropQMD}}
By Lemma \ref{LemSupport} in Appendix \ref{sec:AppendixSupport}, $\nabla \sqrt{q_\ve}$ admits a version such that $\nabla \sqrt{q_\ve}=(\nabla \sqrt{q_\ve}) \1[q_\ve>0]$. Thus, we have
\begin{align*}
\A_{\theta_0}[h](x,y) \sqrt{p_{\theta_0}(x,y)} 
&=
-2\ps{(\nabla\sqrt q_\ve)\big(y - \GGG(\theta_0)(x)\big)}{\I_{\theta_0}[h](x)}_{\R^p} \1\Big[q_{\ve}\big(y-\GGG(\theta_0)(x)\big)>0\Big]
\\[2mm]
&=
-2\ps{(\nabla\sqrt q_\ve)\big(y - \GGG(\theta_0)(x)\big)}{\I_{\theta_0}[h](x)}_{\R^p}
\end{align*}
Writing $\varphi(y)\equiv \sqrt{q_\ve(y)}$ and denoting by $I_\varphi(s)$ the integral in (\ref{eq:QMD}), a change of variables yields
$$
I_\varphi(s) 
=
\int_{\XX\times \R^p} 
\bigg\{
\varphi\Big(y-\big(\GGG(\theta_0+sh)(x) - \GGG(\theta_0)(x)\big)\Big) - \varphi(y)
+ 
s\ps{\nabla \varphi(y)}{\I_{\theta_0}[h](x)}_{\R^p}
\bigg\}^2\, d\lambda(x)\, dy
.
$$
Using that $(a+b)^2\le 2(a^2+b^2)$, we have
$$
I_\varphi(s)
\le 
2(I_1(s) + I_2(s)) 
,
$$
where 
$$
I_1(s) 
\equiv 
\int_{\XX\times \R^p} 
\bigg\{
\varphi \Big(y-\big(\GGG(\theta_0+sh)(x) - \GGG(\theta_0)(x)\big)\Big)
- 
\varphi\big(y-s\I_{\theta_0}[h](x)\big)
\bigg)^2\, d\lambda(x)\, dy
,
$$
and 
$$
I_2(s) 
\equiv 
\int_{\XX\times \R^p} 
\Big\{
\varphi\big(y-s\I_{\theta_0}[h](x)\big) - \varphi(y)
+ 
s\ps{\nabla \varphi(y)}{\I_{\theta_0}[h](x)}_{\R^p}
\Big\}^2\, d\lambda(x)\, dy 
.
$$
Assume for now that $\varphi\in H^1(\R^p)\cap C^1(\R^p)$. The fundamental theorem of calculus combined with Cauchy-Schwarz and Jensen's inequalities then yields 
$$
\int_{\R^p} |\varphi(y-z) - \varphi(y)|^2\, dy 
\le 
\|z\|^2 \int_{\R^p} |\nabla \varphi(y)|^2\, dy
,\spa 
\forall\ z\in\R^p
.
$$
Thus, letting $R_s(x)\equiv \GGG(\theta_0+sh)(x) - \GGG(\theta_0)(x) - s\I_{\theta_0}[h](x)$, we find
$$
I_1(s) 
= 
\int_{\XX\times \R^p} 
\big(
\varphi\big(y-R_s(x)\big)
- 
\varphi(y)
\big)^2\, d\lambda(x)\, dy
\le 
\rho_{\theta_0, h}(s)^2 \|\nabla \varphi\|^2_{L^2(\R^p)}
,
$$
where $\rho_{\theta_0,h}(s)=o(s)$ is the remainder term from Condition~\ref{CondIntroLinModel}. The fundamental theorem and Jensen's inequality provide 
\begin{eqnarray*}
\lefteqn{
I_2(s) 
= 
\int_{\XX\times \R^p} \Big(
\int_0^1 \ps{\nabla \varphi(y) - \nabla \varphi\big(y-ts\I_{\theta_0}[h](x)\big)}{s\I_{\theta_0}[h](x)}_{\R^p}\, dt
\Big)^2\, d\lambda(x)\, dy
}
\\[2mm]
&&
\le 
s^2 \int_{\XX\times \R^p}
\int_0^1 \Big\lvert \ps{\nabla \varphi(y) - \nabla \varphi\big(y-ts\I_{\theta_0}[h](x)\big)}{\I_{\theta_0}[h](x)}_{\R^p}\Big\rvert^2\, dt\, d\lambda(x)\, dy
\equiv 
s^2 J_\varphi(s)
\end{eqnarray*}
We thus obtain 
\begin{equation}\label{eq:BoundI}
I_\varphi(s) 
\le 
2\big( \rho_{\theta_0,h}(s)^2 \|\nabla \varphi\|^2_{L^2(\R^p)} 
+ 
s^2 J_\varphi(s) \big)
,
\end{equation}
for any $\varphi\in C^1(\R^p)\cap H^1(\R^p)$. This bound extends to any $\varphi\in H^1(\R^p)$ by approximation, hence in particular to $\varphi=\sqrt{q_\ve}$. In the general case where $\sqrt{q_\ve}$ is only assumed to belong to $H^1(\R^p)$, we have for any $R>0$
\begin{eqnarray*}
\lefteqn{
J_{\varphi}(s)
\le
\int_0^1 \int_{\XX} |\I_{\theta_0}[h](x)|^2 \bigg( \int_{\R^p} \big|\nabla \varphi(y) - \nabla \varphi\big(y-ts\I_{\theta_0}[h](x)\big)\big|^2\, dy\bigg)\, d\lambda(x)\, dt
}
\\[2mm]
&&\le 
4 \|\nabla \varphi\|^2_{L^2(\R^p)} \int_{|\I_{\theta_0}[h](x)|>R} |\I_{\theta_0}[h](x)|^2\, d\lambda(x) 
\\[2mm]
&&\hspace{10mm}
+
\bigg(
\int_{|\I_{\theta_0}[h](x)|\le R} |\I_{\theta_0}[h](x)|^2\, d\lambda(x) \bigg) 
\sup_{|v|\le sR} \int_{\R^p} \big|\nabla \varphi(y) - \nabla \varphi\big(y+v\big)\big|^2\, dy
.
\end{eqnarray*}
Since the map $v\mapsto \int_{\R^p} |\psi(y+v)-\psi(y)|^2\, dy$ is continuous over $\R^p$ for any $\psi\in L^2(\R^p)$, and recalling that $\rho_{\theta_0,h}(s)=o(s)$ as $s\to 0$, we find for any $R>0$
$$
\limsup_{s\to 0} \frac{I_\varphi(s)}{s^2} 
\le
8 \|\nabla \varphi\|^2_{L^2(\R^p)} \int_{|\I_{\theta_0}[h](x)|>R} |\I_{\theta_0}[h](x)|^2\, d\lambda(x) 
.
$$
Since $\I_{\theta_0}[h]\in L^2_\lambda$, taking $R\to \infty$ yields the conclusion.
\end{Proof}
\vspace{3mm}

Under the assumptions of Proposition \ref{PropQMD} the Fisher information matrix of the law of $\ve$, given by
\begin{equation}\label{eq:defIepsilon}
\II_\ve 
\equiv 
4 \int_{\R^p} \big(\nabla \sqrt{q_\ve}\big)(y)\big(\nabla \sqrt{q_\ve}\big)(y)^T\, dy
=
4\ \E_\ve\bigg[
\frac{(\nabla \sqrt{q_\ve})}{\sqrt{q_\ve}} \frac{(\nabla \sqrt{q_\ve})^T}{\sqrt{q_\ve}}
\bigg]
,
\end{equation}
is well-defined. When $q_\ve\in C^1(\R^p)$ and $\sqrt{q_\ve}\in H^1(\R^p)$, then $\II_\ve$ admits the usual representation 
$$
\II_\ve 
=
\E_\ve\Big[ \big(\nabla \log q_\ve\big) \big(\nabla \log q_\ve\big)^T \Big]
,
$$
by a straightforward application of the chain rule. 

\subsection{LAN expansion of the model}\label{sec:LANlan}

A direct application of Fubini's theorem and Cauchy-Schwarz's inequality entails that the linear score operator $\A_{\theta_0}$ from Proposition \ref{PropQMD} is continuous from $(\VV_0, \ps{\cdot}{\cdot}_\VV)$ to $L^2(\XX\times\R^p, dP_{\theta_0})$. In particular, it admits an adjoint $\A_{\theta_0}^* : L^2(\XX\times\R^p, dP_{\theta_0}) \to (\overline \VV_0, \ps{\cdot}{\cdot}_\VV)$ and the information operator of the model is the linear and continuous operator $\A_{\theta_0}^* \A_{\theta_0} : (\overline \VV_0, \ps{\cdot}{\cdot}_\VV) \to (\overline \VV_0, \ps{\cdot}{\cdot}_\VV)$; see  Section 25.5 in \citet*{Vdv1998} and Section 2 in \citet*{NicPat2023}. Arguing as in the proof of Proposition 1 of \citet*{NicPat2023}, we have $\A_{\theta_0}^* = \I_{\theta_0}^* \circ \EEE_{\theta_0}$, where $\EEE_{\theta_0}:L^2(\XX\times\R^p, dP_{\theta_0})\to L^2(\XX, d\lambda)^p$ is the continuous operator
$$
\EEE_{\theta_0}[w](x)
\equiv 
-2\ \E_{\theta_0}\Big[ w(X,Y)\frac{\nabla \sqrt{q_\ve}}{\sqrt{q_\ve}} \big(Y-\GGG(\theta_0)(X)\big)\ |\ X=x \Big]
,
$$
and satisfies 
$$
\EEE_{\theta_0}[\A_{\theta_0}[h]]
=
\II_\ve \I_{\theta_0}[h]
,\spa 
\forall\ h\in\VV_0 
,
$$
where the matrix $\II_\ve$ from (\ref{eq:defIepsilon}) acts on $\I_{\theta_0}[h](x)\in \R^p$ pointiwse for each $x\in\XX$ through usual matrix multiplication. Consequently, we have
\begin{equation}\label{eq:Info}
\A_{\theta_0}^* \A_{\theta_0}
=
\I_{\theta_0}^* \II_\ve \I_{\theta_0}
,
\end{equation}
The matrix $\II_\ve$  is non-negative definite so that, in particular, it admits the usual square root matrix $\II_\ve^{1/2}$. We introduce the LAN-norm $\|\cdot\|_{\rm LAN}$ that features in Proposition \ref{PropLAN} below as
\begin{equation}\label{eq:LANnorm}
\|h\|_{\rm LAN}
\equiv 
\|\II_\ve^{1/2} \I_{\theta_0}[h]\|_{L^2_\lambda}
,\spa 
\forall\ h\in\VV_0 
.
\end{equation}
Relying on the QMD property established in Proposition \ref{PropQMD}, then arguing as in the proof of Theorem 7.2 in \citet*{Vdv1998} leads to the following LAN expansion.
	
\begin{Prop}\label{PropLAN}
Under the assumptions of Proposition \ref{PropQMD}, the model $\{P^N_\theta : \theta\in\Theta\}$ from (\ref{model}) is locally asymptotically normal (LAN) at $\theta_0$ with respect to the tangent space $\VV_0$, and for all $h\in \VV_0$ fixed we have as $N\to\infty$
\begin{eqnarray*}
\lefteqn{
    \hspace{-32mm}
    \log \frac{dP^N_{\theta_0+h/\sqrt{N}}}{dP^N_{\theta_0}}((X_i, Y_i)_{i=1}^N)
    =
    \frac{2}{\sqrt{N}}\sum_{i=1}^N \ps{\I_{\theta_0}[h](X_i)}{\frac{\nabla \sqrt{q_\ve}}{\sqrt{q_\ve}}(\ve_i)}_{\R^p}
    - \frac12 \|h\|^2_{\rm LAN}
    +
    o_{P^N_{\theta_0}}(1)
}
\\[2mm]
&&
\hspace{20mm}
\overset{d}{\longrightarrow} 
\NN\Big(-\frac12 \|h\|^2_{\rm LAN}, \|h\|^2_{\rm LAN}\Big)
.
\end{eqnarray*}
\end{Prop}

Note that the quotient $(\nabla \sqrt{q_\ve} / \sqrt{q_\ve})(\ve_i)$ in the statement of Proposition \ref{PropLAN} is almost surely well-defined since $q_\ve(\ve_i)>0$ with $q_\ve$-probability one.\smallskip

We will see in Section \ref{sec:FisherInfo} and Section \ref{sec:LimitProcess} below that the information geometry of the model $\{P_{\theta}^N : \theta\in\Theta\}$ around $\theta_0$ crucially relies on the assumption that the {\rm LAN}-norm is an actual norm on $\VV_0$, \ie the map $h\mapsto \II_\ve^{1/2} \I_{\theta_0}[h]$ is injective over $\VV_0$. While we will assume that $\I_{\theta_0}$ defines an injective operator on $\VV_0$, injectivity of $h\mapsto \II_\ve^{1/2} \I_{\theta_0}[h]$ may fail if $\II_\ve$ is degenerate in some directions. The following result establishes that this does not happen, hence $\II_\ve$ is an invertible matrix. 

\begin{Prop}\label{PropInvertIepsilon}
    Under the assumptions of Proposition~\ref{PropLAN}, the Fisher information matrix $\II_\ve$ from~(\ref{eq:defIepsilon}) is positive definite, hence invertible.
\end{Prop}

\begin{Proof}{Proposition \ref{PropInvertIepsilon}}
    Assume, \textit{ad absurdum}, that there exists $v\in\R^p\sm\{0\}$ such that $v^T\II_\ve v=0$ or, equivalently, $\ps{(\nabla \sqrt{q_\ve})(y)}{v}_{\R^p}=0$ for Lebesgue-almost every $y\in\R^p$. Up to rotating the axes, let us assume that $v=(1,0,\ldots,0)\in\R^p$, \ie $(\partial \sqrt{q_\ve}/\partial y_1)(y)=0$ for almost every $y\in\R^p$. Let $\rho\in C^\infty_c(\R^p)$ be non-negative with $\int_{\R^p} \rho(y)\, dy = 1$. Then, $\sqrt{q_\ve}*\rho$ is a smooth function over $\R^p$ (see Theorem 2.1.7 in \citet*{Bog2010}) and the equality
    $$
    \frac{\partial (\sqrt{q_\ve} * \rho)}{\partial y_1}
    =
    \sqrt{q_\ve} * \frac{\partial\rho}{\partial y_1}
    =
    \Big(\frac{\partial \sqrt{q_\ve}}{\partial y_1}\Big) * \rho
    =
    0
    $$
    holds pointwise on $\R^p$. It follows that the map $y\mapsto (\sqrt{q_\ve}*\rho)(y)$ is constant on each line $\{(y_1, z) : y_1\in \R\}$ for $z\in \R^{p-1}$. Consequently, 
    $$
    \int_\R |(\sqrt{q_\ve}*\rho)(y_1, z)|^2\, dy_1
    =
    \infty
    ,\spa 
    \forall\ z\in\R^{p-1} 
    .
    $$
    It follows that 
    $$
    \int_{\R^p} |(\sqrt{q_\ve}*\rho)(y)|^2\, dy 
    =
    \int_{\R^{p-1}} \bigg(\int_\R |(\sqrt{q_\ve}*\rho)(y_1, z)|^2\, dy_1\bigg)\, dz
    =
    \infty 
    .
    $$
    Now, Young's convolutional inequality entails that $\|\sqrt{q_\ve}*\rho\|_{L^2}
    \le 
    \|\sqrt{q_\ve}\|_{L^2} \|\rho\|_{L^1}
    =1
    ,
    $
    a contradiction. We deduce that $v^T \II_\ve v>0$ for all $v\in\R^p\sm\{0\}$, which concludes the proof.
\end{Proof}

\section{Inverting the Fisher information operator}\label{sec:FisherInfo}

Throughout this section, we assume that Condition \ref{CondIntroLinModel} holds for $\theta_0\in\Theta$ fixed and tangent space at $\theta_0$ given by a linear subspace $\VV_0$ of the separable normed vector space $(\VV,\|\cdot\|_\VV)$. To streamline the exposition and highlight the key ideas, in this section we focus on the case where $\VV$ is a separable Hilbert space $(\VV, \ps{\cdot}{\cdot}_\VV)$. The proofs in the case of a general separable normed vector space $(\VV, \|\cdot\|_\VV)$ are deferred to Appendix \ref{sec:InvertibilityBanach}. We thus have a pre-Hilbertian space $(\VV_0, \ps{\cdot}{\cdot}_\VV)$ along which the map $\GGG:\Theta\mapsto L^2_\lambda$ is Gâteaux-differentiable at $\theta_0$ in every direction of $\VV_0$ with continuous (Gâteaux) linearization operator $\I_{\theta_0} : (\VV_0, \ps{\cdot}{\cdot}_\VV)\to L^2_\lambda$. By continuity, $\I_{\theta_0}$ admits a unique continuous extension 
$$
\I_{\theta_0} : (\overline \VV_0, \ps{\cdot}{\cdot}_\VV)\to L^2_\lambda
$$ 
to the $\|\cdot\|_\VV$-closure $\overline \VV_0$ of $\VV_0$, with corresponding adjoint operator
$$
\I_{\theta_0}^* : L^2_\lambda \to (\overline \VV_0, \ps{\cdot}{\cdot}_\VV)
,
$$
characterized by 
\begin{equation}\label{eq:AdjointProp}
\ps{\I_{\theta_0}[h]}{w}_{L^2_\lambda}
=
\ps{h}{\I_{\theta_0}^*[w]}_{\VV}
,\spa 
h\in \overline \VV_0,\ w\in L^2_\lambda 
.
\end{equation}
We further assume that Condition \ref{CondIntroInj} holds, \ie $\I_{\theta_0}$ is injective over $\VV_0$. Under these assumptions, we establish in Theorem \ref{TheorI*IHomeo} below that the resulting information operator initially given, by virtue of (\ref{eq:Info}), by 
$$
\I_{\theta_0}^* \II_\ve \I_{\theta} : (\overline \VV_0, \ps{\cdot}{\cdot}_\VV) \to (\overline \VV_0, \ps{\cdot}{\cdot}_\VV)
,
$$
can be extended to a `natural domain' from which it will define a homeomorphism onto its $\overline\VV_0$-dual. 
This shows that the question of inverting the Fisher information operator between appropriate spaces can \emph{always} be solved as long as the starting linearization $\I_{\theta_0} : (\VV_0, \ps{\cdot}{\cdot}_\VV)\to L^2_\lambda$ is injective and, hence, that the question of surjectivity can always be bypassed for information operators as soon as injectivity holds. When Condition \ref{CondIntroInj} is too difficult to establish for a given $\VV_0$, one can always restrict to a linear subspace $\VV_0'$ of $\VV_0$ and try to establish injectivity of $\I_{\theta_0}$ on $\VV_0'$ instead. The results below would thus apply to $\VV_0'$. Although our results still reliy on the non-benign assumption of injectivity, this provides a principled way to establish invertibility of the information operator between explicitly constructed spaces. We will illustrate this on concrete PDE examples in Section \ref{sec:Appli}.



\subsection{Extension to the natural domain}\label{sec:Extend}

To characterize the notion of `natural domain' precisely, let us introduce the natural norm induced by $\I_{\theta_0}$ on $\VV_0$: since $\I_{\theta_0}$ is injective, we can endow $\VV_0$ with the norm 
$$
\|v\|_{\HH}
\equiv 
\|\II_\ve^{\frac12} \I_{\theta_0}[v]\|_{L^2_\lambda}
,\spa 
v\in \VV_0
,
$$
where $\II_\ve$ is the symmetric matrix from (\ref{eq:defIepsilon}).
The resulting map $\|\cdot\|_{\HH}:\VV_0\to [0,\infty)$ is indeed a norm on $\VV_0$ since $\I_{\theta_0}$ is linear and injective by Condition \ref{CondIntroInj} and $\II_\ve$ is an invertible matrix by Proposition \ref{PropInvertIepsilon}. Note that we have 
\begin{equation}\label{eq:VVdomH}
\|v\|_{\HH} 
= 
\|\II_\ve^{\frac12} \I_{\theta_0}[v]\|_{L^2_\lambda}
\lesssim 
\|\I_{\theta_0}[v]\|_{L^2_\lambda}
\lesssim 
\|v\|_{\VV}
,\spa 
v\in \VV_0
,
\end{equation}
by continuity  of $\I_{\theta_0}:(\VV_0,\|\cdot\|_{\VV})\to L^2_\lambda$. 
Denoting by $(\HH, \|\cdot\|_{\HH})$ the completion of $\VV_0$ with respect to the $\|\cdot\|_{\HH}$-norm, we then have the continuous embedding of Banach spaces
\begin{equation}\label{VVEmbedH}
(\overline \VV_0,\|\cdot\|_{\VV})
\embed 
(\HH, \|\cdot\|_{\HH})
.
\end{equation}
Note that $(\VV_0, \|\cdot\|_\VV)$ is separable since it is a linear subspace of the separable normed space $(\VV, \|\cdot\|_\VV)$. Consequently, $\overline\VV_0$ and $\HH$ are also separable Banach spaces, 
since $\VV_0$ is dense in both $\overline \VV_0$ and $\HH$ by construction. The (densely defined) linear map $\I_{\theta_0}:\VV_0\to L^2_\lambda$ satisfies 
$$
\|\II_\ve^{\frac12} \I_{\theta_0}[v]\|_{L^2_\lambda}
=
\|v\|_{\HH} 
,\spa 
v\in \VV_0
.
$$
Since $L^2_\lambda$ is a Banach space and $\VV_0$ is dense in $\HH$, then $\I_{\theta_0}$ admits a unique extension 
\begin{equation}\label{eq:AExtended}
\I_{\theta_0}:(\HH, \|\cdot\|_{\HH}) \to L^2_\lambda
\end{equation}
as a linear and continuous operator satisfying
$$
\|\II_\ve^{\frac12}\I_{\theta_0}[v]\|_{L^2_\lambda}
=
\|v\|_{\HH}
,\spa 
v\in \HH
,
$$
where, for $v\in \HH\sm \VV_0$, $\I_{\theta_0}[v]$ is defined as the (unique) limit of an arbitrary sequence $(v_k)\subset \VV_0$ such that $\|v_k-v\|_{\HH}\to 0$. In particular, the operator $\I_{\theta_0} : \HH \to L^2_\lambda$ is an \emph{injective} linear and continuous operator between Banach spaces.\smallskip

Our choice to denote the natural domain of $\I_{\theta_0}$ by $\HH$ is not meaningless since the corresponding $\|\cdot\|_\HH$-norm is the `LAN norm' $\|\II_\ve^{1/2}\I_{\theta_0}[\cdot]\|_{L^2_\lambda}$ featuring in the LAN expansion of the log-likelihood ratio in Proposition \ref{PropLAN}; as such, $\HH$ will be seen to be the RKHS (or Cameron-Martin space) of the optimal limiting Gaussian process in Section \ref{sec:LimitProcess}.

The information operator $\I_{\theta_0}^* \II_\ve \I_{\theta_0}$ also extends to a continuous operator 
\begin{equation}\label{eq:A*Aextended}
\I_{\theta_0}^* \II_\ve \I_{\theta_0} : \HH \to \overline \VV_0 
.
\end{equation}
As a consequence, the Banach space $\HH$ is in fact a Hilbert space: since $\I_{\theta_0}:\HH\to L^2_\lambda$ is injective, then 
\begin{equation}\label{defHInnerProd}
\ps{u}{v}_{\HH}
\equiv 
\ps{\I_{\theta_0}[u]}{\II_\ve \I_{\theta_0}[v]}_{L^2_\lambda} 
=
\big\langle\II_\ve^{\frac12} \I_{\theta_0}[u] , \II_\ve^{\frac12} \I_{\theta_0}[v]\big\rangle_{L^2_\lambda} 
,\spa 
u,v\in \HH
,
\end{equation}
defines an inner-product on $\HH$ such that $\ps{u}{u}_{\HH}=\|u\|^2_{\HH}$. For arbitrary $u,v\in \HH$, the $\ps{\cdot}{\cdot}_\HH$ inner-product can further be written as 
\begin{equation}\label{defHInnerProd2}
\ps{u}{v}_{\HH}
= 
\ps{u}{\I_{\theta_0}^* \II_\ve \I_{\theta_0}[v]}_{\VV} 
,\spa 
u,v\in \HH
,
\end{equation}
where the previous $\ps{\cdot}{\cdot}_\VV$ inner-product is understood as the limit 
$$
\ps{u}{\I_{\theta_0}^* \II_\ve \I_{\theta_0}[v]}_{\VV} 
\equiv 
\lim_{k\to \infty}
\ps{u_k}{\I_{\theta_0}^* \II_\ve \I_{\theta_0}[v]}_{\VV}
$$
for a sequence $(u_k)\subset \VV_0$ such that $\|u_k-u\|_\HH\to 0$.

\begin{Example}\label{ex:Stability}
    Although $\HH$ is defined intrinsically from $\I_{\theta_0}$, in concrete models it is usually possible to identify $\HH$ through two-sided estimates 
    \begin{equation}\label{twoSidedEstimates}
    \|h\|_X
    \lesssim
    \|\I_{\theta_0}[h]\|_{L^2_\lambda}
    \lesssim 
    \|h\|_X
    ,
    \end{equation}
    for some norm $\|\cdot\|_X$ on $\VV_0$.
    Indeed, in this case the set $\HH$ coincides with the completion of~$\VV_0$ for the $\|\cdot\|_X$-norm. For instance, consider $\VV=L^2(\ZZ, \R^m)$, $m\ge 1$, where $\ZZ\subset \R^n$ is a bounded domain with smooth boundary and $n\ge 1$, and $\VV_0\subset \VV$ a linear subspace containing~$C^\infty_c(\ZZ, \R^m)$ (or any dense subset of $L^2(\ZZ,\R^m)$).
    If (\ref{twoSidedEstimates}) holds for the dual norm $\|\cdot\|_X=\|\cdot\|_{(H^\kappa_0(\ZZ,\R^m))^*}$, where $\kappa\in\N$, then $\HH = (H_0^{\kappa}(\ZZ,\R^m))^*$. Here, the space $H^{\kappa}_0(\ZZ,\R^m)$ is the closure of $C^\infty_c(\ZZ,\R^m)$ with respect to the usual $H^\kappa(\ZZ,\R^m)$-norm, and is a separable Hilbert space for the inner-product 
    $$
    \ps{u}{v}_{H^\kappa_0(\ZZ, \R^m)} 
    \equiv  
    \sum_{\substack{\alpha\in\N^n\\ |\alpha|=\kappa}} \ps{\partial^\alpha u}{\partial^\alpha v}_{L^2(\ZZ,\R^m)}
    ,
    $$
    with induced norm equivalent to the usual $H^\kappa(\ZZ^, \R^m)$-norm on $H^\kappa_0(\ZZ,\R^m)$ as a consequence of Poincaré's inequality; see, e.g., Corollary 6.31 in \citet*{Adams2003}. It should also be emphasized that the equality $\HH=(H^\kappa_0(\ZZ,\R^m))^*$ holds as sets but their norms might be different, although they are equivalent. For basic definitions and properties of (positive and negative order) Sobolev spaces, we refer the reader to Chapter 5 in \citet*{Evans1998} or Chapter 2 in \citet*{Adams2003}, and more specifically to §3.13 in the latter for details on $(H^\kappa_0(\ZZ,\R^m))^*$.
    
\end{Example}

\begin{Remark}
In the setting of Example {\ref{ex:Stability}} above, when $\VV_0$ consists of maps satisfying further constrained, the resulting space $\HH$ typically reflects these constraints when they are `closed' properties. For instance, when $\GGG(\theta)$ is the solution to the Navier-Stokes equations considered in Section~\ref{sec:AppliNS} below, $\I_\theta[h]$ is a divergence-free vector field for any $h$ in the relevant tangent space. In this case, the two-sided estimates (\ref{twoSidedEstimates}) hold with $\|\cdot\|_X=(H^1_0)^*$, and thus $\HH$ consists of those elements of $(H^1_0)^*$ which, as generalized functions, have zero divergence since convergence in the distributional sense (hence also in $\|\cdot\|_{(H^1_0)^*})$ preserves divergence-freeness.
\end{Remark}

\subsection{Characterization of the range and surjectivity}\label{sec:Range}

To characterize the range of $\I_{\theta_0}^* \II_\ve \I_{\theta_0} : \HH\to \overline \VV_0$, our goal is to define a norm on $\overline \VV_0$ that should take finite values only for elements of $\overline \VV_0$ that lie in $\I_{\theta_0}^* \II_\ve \I_{\theta_0}(\HH)$. At a heuristic level, we aim to define the norm of such an element $s=\I_{\theta_0}^* \II_\ve \I_{\theta_0}[\bar s]$, for some $\bar s\in\HH$, simply as $\|\bar s\|_\HH=\|(\I_{\theta_0}^* \II_\ve \I_{\theta_0})^{-1}[s]\|_\HH$. Of course, this definition only makes sense if we already know that $s\in \I_{\theta_0}^* \II_\ve\I_{\theta_0}(\HH)$, so we must find another definition that does not explicitly make use of the inverse of $\I_{\theta_0}^* \II_\ve\I_{\theta_0}$. To do so, observe that since $\VV_0$ is dense in $\HH$, we formally have by (\ref{defHInnerProd2})
\begin{align*}
\|(\I_{\theta_0}^* \II_\ve\I_{\theta_0})^{-1}[s]\|_\HH
&=
\sup_{\substack{v\in \VV_0\\ \|v\|_\HH\le 1}} \lvert\ps{(\I_{\theta_0}^* \II_\ve\I_{\theta_0})^{-1}[s]}{v}_\HH\rvert
\\[2mm]
&=
\sup_{\substack{v\in \VV_0\\ \|v\|_\HH\le 1}} \lvert\ps{\I_{\theta_0}^* \II_\ve \I_{\theta_0}(\I_{\theta_0}^* \II_\ve\I_{\theta_0})^{-1}[s]}{v}_{\VV}\rvert
\\[2mm]
&=
\sup_{\substack{v\in \VV_0\\ \|v\|_\HH\le 1}} \lvert\ps{s}{v}_{\VV}\rvert
.
\end{align*}
Observe now that the last expression can be defined for any $s\in \overline \VV_0$ without assuming any mapping property of $\I_{\theta_0}^* \II_\ve\I_{\theta_0}$. This motivates the following definition: let 
\begin{equation}\label{DefS}
\SS
\equiv 
\Big\{ 
h\in \overline \VV_0 : 
\|h\|_\SS 
\equiv 
\sup_{\substack{v\in \VV_0\\ \|v\|_{\HH} \le 1}} \lvert\ps{h}{v}_{\VV}\rvert < \infty 
\Big\} 
,
\end{equation}
which we refer to as the $\overline\VV_0$-dual of $\HH$ in view of Proposition \ref{PropDual} in Section \ref{sec:duality} below. By definition, $\SS$ is a subset of $\overline \VV_0$. Recalling the embedding $\overline \VV_0\embed \HH$ from (\ref{eq:VVdomH}) and (\ref{VVEmbedH}), we further have for all $h\in \SS$
$$
\|h\|_\SS 
=
\sup_{\substack{v\in \VV_0\\ \|v\|_\HH\le 1}} \lvert \ps{h}{v}_\VV \rvert 
\gtrsim
\sup_{\substack{v\in \VV_0\\ \|v\|_\VV\le 1}} \lvert \ps{h}{v}_\VV \rvert 
=
\|h\|_\VV
,
$$
where we used that $\VV_0$ is dense in $(\overline \VV_0, \ps{\cdot}{\cdot}_\VV)$ in the last equality. This yields the continuous embeddings 
\begin{equation}\label{EmbeddingSVVH}
(\SS, \|\cdot\|_\SS)
\embed 
(\overline \VV_0, \|\cdot\|_\VV) 
\embed 
(\HH, \|\cdot\|_\HH) 
.
\end{equation}
In fact, $\SS$ identifies isometrically with the dual of $\HH$ for the $\SS$-$\HH$ pairing given by the $\overline\VV_0$-inner product (see Proposition \ref{PropDual} in Section \ref{sec:duality} below). Then, $\overline \VV_0$ acts as a pivot space about which dual spaces are represented; this is similar to Sobolev spaces of opposite smoothness which are dual to each other when the duality pairing is represented by the $L^2$-inner product.

\begin{Example}\label{ex:HandS}
    In the setting of Example \ref{ex:Stability}, with $\VV=L^2(\ZZ, \R^m)$, $\ZZ\subset \R^n$ a bounded domain with smooth boundary, and $\VV_0\subset \VV$ a linear subspace containing~$C^\infty_c(\ZZ, \R^m)$, if 
    $$
    \|h\|_{(H^\kappa_0(\ZZ, \R^m))^*}
    \lesssim
    \|\I_{\theta_0}[h]\|_{L^2_\lambda}
    \lesssim 
    \|h\|_{(H^\kappa_0(\ZZ, \R^m))^*}
    ,\spa 
    h\in \VV_0
    ,
    $$
    then $\HH = (H^\kappa_0(\ZZ, \R^m))^*$ and $\SS=H^\kappa_0(\ZZ,\R^m)$, in view of the $L^2$-duality between $H^\kappa_0(\ZZ,\R^m)$ and $(H^\kappa_0(\ZZ,\R^m))^*$; see also Proposition \ref{PropDual} in Section \ref{sec:duality} below.
\end{Example}

\begin{Example}\label{ex:HandSBis}
Consider the case
$$
\|h\|_{\DD^{-\kappa}}
\lesssim
\|\I_{\theta_0}[h]\|_{L^2_\lambda}
\lesssim 
\|h\|_{\DD^{-\kappa}}
,\spa 
h\in \VV_0
,
$$
where, for an orthonormal basis $(e_j)$ of $\overline \VV_0$, and a positive sequence $(\tau_j)\to\infty$
$$
\DD^s
\equiv 
\Big\{ \sum_{j} u_j e_j : \|u\|^2_{\DD^s}\equiv \sum_j \tau_j^{s} u_j^2 < \infty \Big\} 
,\spa 
s\in\R 
.
$$
The series $\sum_j u_j e_j$ is to be understood in the sense of the convergence in the $\|\cdot\|_{\DD^s}$-norm. Each $\DD^s$ is a separable Hilbert space for the inner-product $\ps{u}{v}_{\DD^s}\equiv  \sum_j \tau_j^s u_j v_j$, and $(\DD^s)_{s\in\R}$ forms a nested family with $\DD^0 = \overline \VV_0$ and $\DD^\beta\subset \DD^\alpha$ for any real numbers $\alpha\le \beta$. It is straightforward to see that $\DD^{-s}=(\DD^s)^*$ for any $s\in\R$, so that the two-sided estimates above yield $\HH=\DD^{-\kappa}$ and $\SS=\DD^\kappa$ (see Remark \ref{RemEquivNorms} below).
\end{Example}

An important case of Example \ref{ex:HandSBis} above is when $\VV=L^2(\ZZ, \R^m)$, where $\ZZ\subset \R^n$ is a bounded domain with smooth boundary, and when $(e_j, \tau_j)$ are the eigenpairs of the Dirichlet Laplacian over $\ZZ$; then, $(e_j)\subset H^1_0(\ZZ,\R^m)\cap C^\infty(\overline \ZZ, \R^m)$ is an orthonormal system of $L^2(\ZZ,\R^m)$, and orthogonal in $H^1_0(\ZZ, \R^m)$, with $-\Delta e_j=\tau_j e_j$ over $\ZZ$ and $e_j=0$ over $\partial \ZZ$; see, e.g., Theorem~6.5.1 in \citet*{Evans1998}. Weyl's asymptotics yields $\tau_j\sim j^{2/n}$ as $j\to\infty$ (see Corollary~8.3.5 in \citet*{TaylorPDE2}), and we further have $\DD^1=H^1_0(\ZZ, \R^m)$ and $\DD^2=H^2_0(\ZZ,\R^m)$, as can be seen from integration by parts and Parseval's formula in $L^2(\ZZ,\R^m)$, hence also $\DD^{-1}=(H^1_0(\ZZ,\R^m))^*$ and $\DD^{-2}=(H^2_0(\ZZ, \R^m))^*$. Unless $\overline{\Om}$ is a smooth manifold with no boundary, we only have $\DD^j \subset H^j_0(\ZZ, \R^m)$ in general for $j\ge 2$; see, e.g., Section 5.A in \citet*{Taylor2011}.

\begin{Remark}\label{RemEquivNorms}
Assume that, for some norm $\|\cdot\|_X$ on $\VV_0$, we have 
$$
\|h\|_X 
\lesssim 
\|\I_{\theta_0}[h]\|_{L^2_\lambda}\lesssim 
\|h\|_X 
,\spa 
\forall\ h\in \VV_0 
.
$$
Since $\SS$ is the $\overline \VV_0$-dual of $\HH$ (see Proposition \ref{PropDual} in Section \ref{sec:duality} below), we expect $\|\cdot\|_\SS$ to be equivalent to the $\overline \VV_0$-dual norm of $\|\cdot\|_X$. Indeed, we have 
$$
\|h\|_\SS 
=
\sup_{\substack{v\in \VV_0\\ \|v\|_\HH\le 1}} \lvert\ps{h}{v}_\VV\rvert 
\asymp
\sup_{\substack{v\in \VV_0\\ \|v\|_X\le 1}} \lvert\ps{h}{v}_\VV\rvert
.
$$
In particular, in Example \ref{ex:HandSBis} above, we have $\|\cdot\|_\SS\asymp \|\cdot\|_{\DD^\kappa}$.
\end{Remark}

Let us now turn to studying the range $R(\A_{\theta_0}^*)$ of $\A_{\theta_0}^* : L^2(\XX\times\R^p, dP_{\theta_0}) \to (\overline\VV_0, \ps{\cdot}{\cdot}_\VV)$. The equality $\A_{\theta_0}^* = \I_{\theta_0}^*\circ \EEE_{\theta_0}$ from Section \ref{sec:LANlan} yields $R(\A_{\theta_0}^*)\subset R(\I_{\theta_0}^*)$. Arguing as in the proof of Proposition 1 in \citet*{NicPat2023}, then for any $h\in L^2_\lambda$, taking 
$$
w(x,y)
=
-2\Big\langle \II_\ve^{-1} h(x) , 
\frac{\nabla \sqrt{q_\ve}}{\sqrt{q_\ve}}\big(y-\GGG(\theta_0)(x)\big)\Big\rangle_{\R^p}
\in L^2(\XX\times\R^p, dP_{\theta_0})
,
$$
yields $\EEE_{\theta_0}[w]=h$ to the effect that $R(\I_{\theta_0}^*)\subset R(\A_{\theta_0}^*)$, hence $R(\I_{\theta_0}^*)= R(\A_{\theta_0}^*)$.
For all $w\in L^2_\lambda$ and $v\in \VV_0$, observe now that
\begin{eqnarray}\label{ContA*toS}
\lefteqn{
\nonumber
\|\I_{\theta_0}^*[w]\|_\SS 
=
\sup_{\substack{v\in \VV_0\\ \|v\|_\HH\le 1}} 
\lvert\ps{\I_{\theta_0}^*[w]}{v}_{\VV}\rvert
=
\sup_{\substack{v\in \VV_0\\ \|v\|_\HH\le 1}} 
\lvert\ps{w}{\I_{\theta_0}[v]}_{L^2_\lambda}\rvert 
}
\\[2mm]
&&\le
\sup_{\substack{v\in \VV_0\\ \|v\|_\HH\le 1}} 
\|w\|_{L^2_\lambda} \|\I_{\theta_0}[v]\|_{L^2_\lambda} 
=
\sup_{\substack{v\in \VV_0\\ \|v\|_\HH\le 1}} 
\|w\|_{L^2_\lambda} \|v\|_{\HH}
\le
\|w\|_{L^2_\lambda}
.
\end{eqnarray}
In particular, $\I_{\theta_0}^*[w]\in \SS$ for all $w\in L^2_\lambda$. Consequently, we have the following inclusions of sets 
\begin{equation}\label{InclusionSets}
\I_{\theta_0}^* \II_\ve \I_{\theta_0}(\HH) 
\subset 
R(\I_{\theta_0}^*) 
=
R(\A_{\theta_0}^*)
\subset 
\SS
.
\end{equation}
The following result establishes that $\SS$ exactly captures all elements of $\overline \VV_0$ that can be attained by $\I_{\theta_0}^*$ and $\I_{\theta_0}^* \II_\ve \I_{\theta_0}$ to the effect that all inclusions in the previous display are in fact equalities of sets.

\begin{Prop}\label{PropEqualityRanges}
We have 
$
R(\A_{\theta_0}^*)
=
R(\I_{\theta_0}^*)
= 
\I_{\theta_0}^* \II_\ve \I_{\theta_0}(\HH) 
=
\SS
.
$
\end{Prop}

\begin{Proof}{Proposition \ref{PropEqualityRanges}}
In view of (\ref{InclusionSets}), it is enough to establish that $\SS\subset \I_{\theta_0}^* \II_\ve \I_{\theta_0}(\HH)$. Then, fix $h\in \SS$, so that 
$$
\lvert\ps{h}{v}_\VV\rvert
\le 
\|h\|_\SS \|v\|_{\HH}
,\spa 
v\in \VV_0
.
$$
Since $\VV_0$ is dense in $\HH$, then the linear and continuous map $\ell : (\VV_0, \ps{\cdot}{\cdot}_{\HH})\to \R,\ v\mapsto \ps{h}{v}_{\VV}$ extends to a continuous map defined on all of $\HH$. Since $(\HH,\ps{\cdot}{\cdot}_\HH)$ is a Hilbert space for $\ps{\cdot}{\cdot}_\HH$ defined in (\ref{defHInnerProd}), the Riesz representation theorem entails that there exists $\bar h\in \HH$ such that 
$$
\ps{h}{v}_\VV
=
\ell(v)
=
\ps{\bar h}{v}_{\HH}
=
\ps{\I_{\theta_0}^* \II_\ve \I_{\theta_0}[\bar h]}{v}_{\VV}
,\spa 
v\in \VV_0
,
$$
by virtue of (\ref{defHInnerProd2}).
In particular, $h = \I_{\theta_0}^* \II_\ve \I_{\theta_0}[\bar h]$ in $\overline \VV_0$. We deduce that $h \in \I_{\theta_0}^* \II_\ve \I_{\theta_0}(\HH)$. Together with (\ref{InclusionSets}), this yields the result. 
\end{Proof}
\vspace{3mm}

In the previous examples, we have seen that $\SS$ is always a `nice' space in the sense that it is a well-identified subspace of $L^2$ equipped with a natural norm that makes it a Banach (even Hilbert) space. In the following result, we show that $\SS$ endowed with the $\|\cdot\|_\SS$-norm is indeed always a Banach space.

\begin{Prop}\label{PropBanach}
    The set $\SS$ endowed with the norm $\|\cdot\|_\SS$ from (\ref{DefS}) is a Banach space. 
    
\end{Prop}

\begin{Proof}{Proposition \ref{PropBanach}}
    The fact that $\SS$ is a linear space follows immediately from the definition (\ref{DefS}) and the triangle inequality. The homogeneity and triangle inequality for $\|\cdot\|_\SS$ are obvious, and if $w\in \SS$ satisfies $\|w\|_\SS=0$, then $\ps{w}{v}_{\VV}=0$ for all $v\in \VV_0$ with $\|\I_{\theta_0}[v]\|_{L^2_\lambda}\le 1$; by linearity, this implies that $\ps{w}{v}_{\VV}=0$ for all $v\in \VV_0$ hence also for all $v\in \overline \VV_0$, to the effect that $w=0$ in $\overline \VV_0$. Consequently, $(\SS, \|\cdot\|_\SS)$ is a normed vector, and it remains to see that it is complete. For this purpose, observe that, for all $s\in \SS$ with corresponding $\bar s\in \HH$ such that $s=\I_{\theta_0}^* \II_\ve \I_{\theta_0}[\bar s]$ (Proposition \ref{PropEqualityRanges}), we have by virtue of (\ref{defHInnerProd2})
    \begin{equation}\label{eq:SIsometricH}
    \|s\|_\SS 
    =
    \sup_{\substack{v\in \VV_0\\ \|v\|_\HH\le 1}} \lvert \ps{\I_{\theta_0}^* \II_\ve \I_{\theta_0}[\bar s]}{v}_\VV\rvert
    =
    \sup_{\substack{v\in \VV_0\\ \|v\|_\HH\le 1}} \lvert \ps{\bar s}{v}_{\HH}\rvert
    =
    \|\bar s\|_\HH
    ,
    \end{equation}
    where the last equality follows since $\VV_0$ is dense in $\HH$. Thus, let $(w_n)\subset \SS$ be a Cauchy sequence in $(\SS,\|\cdot\|_\SS)$. For each $n$, let $\bar w_n\in \HH$ be such that $w_n=\I_{\theta_0}^* \II_\ve \I_{\theta_0}[\bar w_n]$, by virtue of Proposition~\ref{PropEqualityRanges}. Then (\ref{eq:SIsometricH}) applied to $s=w_n-w_m\in\SS$ entails that $(\bar w_n)$ is a Cauchy sequence in the Banach space $(\HH, \|\cdot\|_\HH)$. Thus, there exists $\bar w_\infty\in\HH$ such that $\|\bar w_n - \bar w_\infty\|_\HH \to 0$. Defining $w_\infty \equiv \I_{\theta_0}^* \II_\ve \I_{\theta_0}[\bar w_\infty]$, then Proposition \ref{PropEqualityRanges} entails that $w_\infty \in \SS$ so that, by virtue of (\ref{eq:SIsometricH}), we have 
    $$
    \|w_n - w_\infty\|_\SS 
    =
    \|\bar w_n - \bar w_\infty\|_\HH 
    \to 
    0
    .
    $$
    We deduce that $(\SS, \|\cdot\|_\SS)$ is complete, hence is a Banach space.
\end{Proof}

\begin{Remark}
    Although $(\SS,\|\cdot\|_\SS)$ is a Banach space, and the ranges of $\I_{\theta_0}^*:L^2_\lambda\to \overline \VV_0$ and $\I_{\theta_0}^* \II_\ve \I_{\theta_0} : \HH\to \overline \VV_0$ are equal to $\SS$, this does not necessarily mean that $\I_{\theta_0}^*$ and $\I_{\theta_0}^* \II_\ve \I_{\theta_0}$ have closed range in $\overline \VV_0$ as $\SS$ need not be a closed subset of $(\overline \VV_0, \|\cdot\|_{\VV})$.
\end{Remark}



\subsection{Injectivity and invertibility}

We have seen that the range of the continuous linear operator $\I_{\theta_0}^* : L^2_\lambda\to \overline \VV_0$ is equal to a Banach space $(\SS,\|\cdot\|_\SS)$, given by (\ref{DefS}), and (\ref{ContA*toS}) establishes that 
$$
\I_{\theta_0}^* : (L^2_\lambda,\|\cdot\|_{L^2_\lambda})\to (\SS,\|\cdot\|_\SS)
$$ 
is a surjective and continuous mapping between Banach spaces. In particular, Proposition \ref{PropEqualityRanges} entails that the information operator
$$
\I_{\theta_0}^* \II_\ve \I_{\theta_0} : (\HH,\|\cdot\|_{\HH}) \to (\SS,\|\cdot\|_\SS)
$$
is also a surjective and continuous mapping between Banach spaces.

\begin{Theor}\label{TheorI*IHomeo}
    The linear operator $\I_{\theta_0}^* \II_\ve \I_{\theta_0} : (\HH,\|\cdot\|_{\HH}) \to (\SS,\|\cdot\|_\SS)$ is an isometric homeomorphism.
\end{Theor}

\begin{Proof}{Theorem \ref{TheorI*IHomeo}}
    We already established that $\I_{\theta_0}^* \II_\ve \I_{\theta_0} : (\HH,\|\cdot\|_{\HH}) \to (\SS,\|\cdot\|_\SS)$ defines a surjective linear and continuous operator. Since $\HH$ and $\SS$ are Banach spaces, it follows from the inverse mapping theorem that $\I_{\theta_0}^* \II_\ve \I_{\theta_0}$ is a homeomorphism if and only if it is a continuous bijection. Thus, it remains to establish injectivity. Let $h\in \HH$ be such that $\I_{\theta_0}^* \II_\ve \I_{\theta_0}[h]=0$ in $\SS$. Recalling (\ref{eq:SIsometricH}) from the proof of Proposition \ref{PropBanach}, we find
    $$
    0
    =
    \|\I_{\theta_0}^* \II_\ve \I_{\theta_0}[h]\|_\SS 
    =
    \|h\|_\HH
    ,
    $$
    which yields injectivity of $\I_{\theta_0}^* \II_\ve \I_{\theta_0} : (\HH, \|\cdot\|_\HH) \to (\SS, \|\cdot\|_\SS)$. This further establishes that $\I_{\theta_0}^* \II_\ve \I_{\theta_0}$ is an isometry between $(\HH, \|\cdot\|_\HH)$ and $(\SS, \|\cdot\|_\SS)$, which concludes the proof.
\end{Proof}

\subsection{Duality}\label{sec:duality}

\begin{Prop}\label{PropDual}
The Banach space $(\SS, \|\cdot\|_\SS)$ identifies isometrically with the topological dual of $(\HH, \ps{\cdot}{\cdot}_\HH)$, \ie for any continuous and linear map $\ell:(\HH, \ps{\cdot}{\cdot}_\HH)\to \R$ there exists a unique $\psi\in \SS$ such that $\ell(h)=\ps{\psi}{h}_\VV$ for all $h\in \overline\VV_0$, and the operator norm of $\ell$ is equal to $\|\psi\|_\SS$. Conversely, for any $\psi\in\SS$ the linear map $\ell_\psi(h)\equiv \ps{\psi}{h}_\VV$, $h\in\VV_0$, extends uniquely to a continuous and linear map $\ell:(\HH, \ps{\cdot}{\cdot}_\HH)\to \R$ with operator norm equal to $\|\psi\|_\SS$.
\end{Prop}

\begin{Proof}{Proposition \ref{PropDual}}
Let $\ell:(\HH, \ps{\cdot}{\cdot}_\HH)\to \R$ be a continuous and linear map. The Riesz representation theorem entails that there exists $\bar \psi\in\HH$ such that $\ell(h)=\ps{\bar \psi}{h}_\HH$ for all $h\in\HH$. In particular, we have 
$$
\ell(h)
=
\ps{\I_{\theta_0}^* \II_\ve \I_{\theta_0}[\bar \psi]}{h}_\VV 
, \spa 
\forall\ h\in \VV_0
,
$$
by definition of $\ps{\cdot}{\cdot}_\HH$ (see (\ref{defHInnerProd})) and the adjoint property. Theorem \ref{TheorI*IHomeo} entails that $\psi \equiv \I_{\theta_0}^* \II_\ve \I_{\theta_0}[\bar \psi]\in \SS$ and the equality in the previous display extends to all $h\in\overline\VV_0$ by approximation. In addition, the operator norm of $\ell$ is equal to $\|\bar\psi\|_\HH=\|\psi\|_\SS$ by virtue of (\ref{eq:SIsometricH}). Now let $\psi\in\SS$ and define the linear map $\ell_\psi(h)\equiv \ps{\psi}{h}_\VV$ for all $h\in \VV_0$. By definition of $\SS$ (see (\ref{DefS})), we have $|\ell_\psi(h)|\le \|\psi\|_\SS \|h\|_\HH$ for all $h\in \VV_0$. We deduce that $\ell_\psi$ extends uniquely to a continuous and linear map $\ell:(\HH, \ps{\cdot}{\cdot}_\HH)\to \R$.
\end{Proof}

In particular, $\SS$ can be endowed with a Hilbert space structure inherited from that of $\HH$ through $\I_{\theta_0}^* \II_\ve \I_{\theta_0}$ by letting 
\begin{equation}\label{SHilbert}
\ps{u}{v}_{\SS}
\equiv 
\ps{u}{(\I_{\theta_0}^* \II_\ve \I_{\theta_0})^{-1}[v]}_{\VV}
,\spa 
u,v\in \SS  
.
\end{equation}
As in (\ref{defHInnerProd2}), the previous $\ps{\cdot}{\cdot}_\VV$-inner product is understood as the limit 
$$
\ps{u}{(\I_{\theta_0}^* \II_\ve \I_{\theta_0})^{-1}[v]}_{\VV}
\equiv 
\lim_{k\to \infty} \ps{u}{y_k}_{\VV}
,
$$
for a sequence $(y_k)\subset \VV_0$ such that $\|y_k-(\I_{\theta_0}^* \II_\ve \I_{\theta_0})^{-1}[v]\|_\HH\to 0$. Note that such a sequence $(y_k)$ always exists since $(\I_{\theta_0}^* \II_\ve \I_{\theta_0})^{-1}[v]\in\HH$ by Theorem \ref{TheorI*IHomeo} and $\VV_0$ is dense in $\HH$. Furter observe that the $\ps{\cdot}{\cdot}_\SS$-inner product is compatible with the $\|\cdot\|_\SS$-norm. Indeed, we have 
\begin{equation}\label{eq:CompatS}
\ps{u}{(\I_{\theta_0}^* \II_\ve \I_{\theta_0})^{-1}[u]}_{\VV}
=
\ps{(\I_{\theta_0}^* \II_\ve \I_{\theta_0})^{-1}[u]}{(\I_{\theta_0}^* \II_\ve \I_{\theta_0})^{-1}[v]}_{\HH}
=
\|(\I_{\theta_0}^* \II_\ve \I_{\theta_0})^{-1}[u]\|^2_\HH 
=
\|u\|^2_\SS
,
\end{equation}
following the arguments at the beginning of Section \ref{sec:Range}. 

Notice that an orthonormal basis of $\HH$ can be realised as a subset of $\overline \VV_0$. Indeed, there exists a countable subset of $\VV_0$ that is dense in $\HH$ in view of the remarks after (\ref{VVEmbedH}), so that one can apply the Gram-Schmidt procedure to a maximal subset of independent vectors in the space $(\HH, \ps{\cdot}{\cdot}_\HH)$ to obtain an orthonormal basis $(h_j)$ of $\HH$ such that $h_j\in \VV_0$ for all $j$.

\begin{Prop}\label{PropRangeBasis}
    Let $(h_j)\subset \VV_0$ be an orthonormal basis of $(\HH, \ps{\cdot}{\cdot}_\HH)$. For all $\psi\in \overline \VV_0$, we have $\psi\in\SS$ if and only if 
    $
    \sum_{j} \lvert\ps{\psi}{h_j}_\VV\rvert^2 
    <
    \infty 
    .
    $
    In this case, we have $\|\psi\|^2_\SS=\sum_{j} \lvert\ps{\psi}{h_j}_\VV\rvert^2$.
\end{Prop}

\begin{Proof}{Proposition \ref{PropRangeBasis}}
Fix $k\in\N$ with $k\le \dim \HH$ smaller than the dimension of $\HH$. Then, letting $\HH_k\equiv {\rm span}\{h_1,\ldots, h_k\}$, we have
$$
\sum_{k=1}^k \lvert\ps{\psi}{h_j}_\VV\rvert^2 
=
\sup_{\substack{\alpha\in\R^k\\ |\alpha|\le 1}} \Big| \sum_{j=1}^k \alpha_j \ps{\psi}{h_j}_{\VV}\Big| 
=
\sup_{\substack{\alpha\in\R^k\\ |\alpha|\le 1}} \Big| \Big\langle \psi , \sum_{j=1}^k \alpha_j h_j \Big\rangle_{\VV}\Big| 
=
\sup_{\substack{h\in \HH_k\\ \|h\|_\HH\le 1}} \lvert \ps{\psi}{h}_\VV\rvert 
.
$$
The r.h.s. converges to $\|\psi\|_\SS$ defined in (\ref{DefS}) as $k\uparrow \dim\HH$ since $\HH_k\subset \VV_0$.
\end{Proof}

\section{The optimal limiting Gaussian process}\label{sec:LimitProcess}


\subsection{Covariance and series expansion}\label{sec:CovSeries}

Let $\{W(f) : f\in \SS\}$ be the centred Gaussian process with covariance 
\begin{equation}\label{CovW}
\E\big[W(f)W(g)\big]
\equiv 
\ps{f}{(\I_{\theta_0}^* \II_\ve\I_{\theta_0})^{-1}[g]}_\VV
,\spa 
f,g\in \SS 
.
\end{equation}
Note that this covariance is well-defined by (\ref{SHilbert}) and subsequent remarks. This process is an infinite-dimensional version of the `optimal' limiting Gaussian random vector with efficient covariance given by the inverse of the Fisher information matrix from parametric statistics. In fact, this process is nothing but the isonormal (or white noise) process on $\SS$ in view of the Hilbertian structure of $\SS$ from (\ref{SHilbert}). This process will be seen below to have RKHS equal to $(\HH, \ps{\cdot}{\cdot}_\HH)$, so to construct a Gaussian Borel-measurable random element with that covariance, consider the random series
\begin{equation}\label{LimitProc}
G 
\equiv 
\sum_{k} \gamma_k h_k 
,
\end{equation}
where $(h_k)$ is an orthonormal basis of $(\HH, \ps{\cdot}{\cdot}_\HH)$ since $\HH$ is separable (see Section \ref{sec:Extend}), $(\gamma_k)$ is an i.i.d. sequence of $\NN(0,1)$-distributed random variables, and the series runs for $k$ less than or equal to the dimension of $\HH$. Assume that $(\X, \|\cdot\|_\X)$ is a Banach space containing $\HH$ as a dense subset such that the canonical injection $(\HH, \|\cdot\|_\HH) \embed (\X, \|\cdot\|_\X)$ is continuous (in particular, since $(\HH, \|\cdot\|_\HH)$ is separable, then $(\X, \|\cdot\|_\X)$ is also separable), and assume that the series (\ref{LimitProc}) defining $G$ converges in $\X$ almost surely. Then, $G$ defines an $\X$-valued Gaussian Borel random variable and the RKHS of the law of $G$ in $\X$ is indeed $\HH$ (see, e.g., Theorem I.23 in \citet*{Vdv2017}).\smallskip

We now argue that the covariance of $G$ is given by (\ref{CovW}) through restriction to a subset of~$\SS$. Similarly to $\SS$ and $\HH$ which are dual to each other with respect to $\overline \VV_0$ (see Section \ref{sec:duality}), we can identify the dual space $\X^*$ with a subspace of $\overline \VV_0$. Indeed, fix $\ell\in \X^*$, \ie $\ell:(\X, \|\cdot\|_\X)\to \R$ is a linear and continuous map, with operator norm $\|\ell\|_{\X^*}$. Since the inclusion $\HH\embed \X$ is continuous, we have 
\begin{equation}\label{eq:ellBound}
|\ell(h)| 
\le 
\|\ell\|_{\X^*} \|h\|_\X 
\lesssim 
\|\ell\|_{\X^*} \|h\|_\HH
,\spa 
v\in  \HH 
.
\end{equation}
In other words, the restriction of $\ell$ to the Hilbert space $( \HH, \ps{\cdot}{\cdot}_\HH)$ is continuous. Consequently, the Riesz representation theorem entails that there exists a unique $\tilde\ell\in\HH$ such that 
\begin{equation}\label{eq:B*VV_0}
\ell(v) 
=
\langle\tilde \ell , v\rangle_\HH
=
\big\langle(\I_{\theta_0}^* \II_\ve \I_{\theta_0})[\tilde\ell] , v\big\rangle_\VV 
\equiv 
\ps{\bar \ell}{v}_\VV
,\spa 
v\in  \HH 
,
\end{equation}
where we used (\ref{defHInnerProd}) and (\ref{defHInnerProd2}), and $\bar\ell\equiv (\I_{\theta_0}^* \II_\ve \I_{\theta_0})[\tilde\ell]$ belongs to $\SS$ by virtue of Theorem \ref{TheorI*IHomeo}.
The fact that $\bar\ell$ is enough to reconstruct $\ell$ follows from the fact that $\HH$ is dense in $\X$, and we have 
\begin{equation}\label{eq:B*B*S}
\|\ell\|_{\X^*}
\equiv 
\sup_{\substack{x\in \X\\ \|x\|_\X\le 1}} |\ell(x)| 
=
\sup_{\substack{h\in \HH\\ \|h\|_\X\le 1}} |\ell(h)| 
=
\sup_{\substack{h\in \HH\\ \|h\|_\X \le 1}} \lvert\ps{\bar\ell}{h}_\VV\rvert 
\equiv 
\|\bar \ell\|_{\X^*_{\overline \VV_0}}
.
\end{equation}
Therefore, $\X^*$ identifies with the subset $\X^*_{\overline \VV_0} \equiv \{\bar\ell : \ell\in \X^*\}$ of $\SS$. Using that $\VV_0$ is dense in $\HH$ and $\HH\embed \X$ continuously yields 
$$
\|\bar\ell\|_{\X^*_{\overline\VV_0}} 
=
\sup_{\substack{v\in \VV_0\\ \|v\|_\X\le 1}} \lvert\ps{\bar\ell}{v}\rvert 
\gtrsim 
\sup_{\substack{v\in \VV_0\\ \|v\|_\HH\le 1}} \lvert\ps{\bar\ell}{v}\rvert
=
\|\bar\ell\|_\SS 
.
$$
Consequently, we have the continuous embeddings 
\begin{equation}\label{allEmbeddings}
(\X^* 
\cong )~
\X^*_{\overline \VV_0} 
\embed 
\SS~ 
(\cong 
\HH^*)
\embed 
\overline\VV_0 
\embed 
\HH 
\embed 
\X 
.
\end{equation}
We may now compute, for all $\ell,\ell'\in \X^*$
\begin{align*}
\E[\ell(G)\ell'(G)]
&=
\sum_{i\ge 1} \ell(h_i)\ell'(h_i)
\\[2mm]
&=
\sum_{i\ge 1} \ps{\bar \ell}{h_i}_\VV \ps{\bar \ell'}{h_i}_\VV 
\\[2mm]
&=
\sum_{i\ge 1} \ps{(\I_{\theta_0}^*\II_\ve \I_{\theta_0})^{-1}[\bar \ell]}{h_i}_\HH 
\ps{(\I_{\theta_0}^* \II_\ve \I_{\theta_0})^{-1}[\bar \ell']}{h_i}_\HH
\\[2mm]
&=
\ps{(\I_{\theta_0}^* \II_\ve \I_{\theta_0})^{-1}[\bar \ell]}{(\I_{\theta_0}^* \II_\ve \I_{\theta_0})^{-1}[\bar \ell']}_\HH
\\[2mm]
&=
\ps{\bar \ell}{(\I_{\theta_0}^* \II_\ve \I_{\theta_0})^{-1}[\bar \ell']}_\VV
,
\end{align*}
where we used Parseval's formula in the Hilbert space $(\HH, \ps{\cdot}{\cdot}_\HH)$ with orthonormal basis $(h_k)$. We deduce that the process $\{\ell(G) : \ell\in \X^*\}$ is a version of the restriction $\{W(\bar \ell) : \bar \ell\in \X^*_{\overline \VV_0}\}$ of $W$ to $\X^*\cong \X^*_{\overline \VV_0} \embed \SS$. From there, we can alternatively show that the RKHS of $W$, or equivalently $G$, equals $\HH$ by looking only at the covariance, instead of the series representation, by applying, for instance, Proposition 2.6.8 in \citet*{GinNic2016}, which essentially amounts to showing that $\HH$ is the $\overline \VV_0$-dual of $\SS$ (Proposition \ref{PropDual} in Section \ref{sec:duality}).\smallskip

The construction of a centred $\X$-valued Gaussian random variable $G$ with covariance (\ref{CovW}) as in (\ref{LimitProc}) requires an orthonormal basis of $(\HH, \ps{\cdot}{\cdot}_\HH)$. Although $\HH$ might be known to us as a set, its metric often remains inaccessible. Instead, when $\HH$ is endowed with another (known) inner-product, we wish to use an orthonormal basis in that metric to construct a variable $\tilde G$ having the same law as $G$. This will be relevant in Section \ref{sec:GeneralMinimax} below (more specifically Corollary~\ref{CoroMinimaxHilbert}), so we describe in the following example how this can be done.

\begin{Example}\label{ex:SeriesG}
Consider the setting of Example \ref{ex:HandSBis}, where $\|\cdot\|_\HH\asymp \|\cdot\|_{\DD^{-\kappa}}$, hence also $\|\cdot\|_\SS\asymp \|\cdot\|_{\DD^{\kappa}}$ (see Remark \ref{RemEquivNorms}), with
$$
\DD^s
\equiv 
\Big\{ \sum_{j} u_j e_j : \|u\|^2_{\DD^s}\equiv \sum_j \tau_j^{s} u_j^2 < \infty \Big\} 
,\spa 
s\in\R 
,
$$
for an orthonormal basis $(e_j)\subset \SS$ of $\overline \VV_0$ 
and $\tau_j>0$ with $\tau_j\to\infty$. Then, the law of $G$ can be realized as the law of 
\begin{equation}\label{defGbis}
\tilde G 
\equiv 
\sum_{k} \tilde \gamma_k e_k
,
\end{equation}
where $(\tilde \gamma_k : k\in\N)$ are centred normal random variables such that 
$$
\E\big[\tilde\gamma_i \tilde\gamma_j\big]
=
\ps{e_i}{(\I_{\theta_0}^* \II_\ve \I_{\theta_0})^{-1}e_j}_\VV 
.
$$
To see this, first define centred normal random variables $\{\tilde W(e_n) : n\in\N\}$ such that 
$$
\E\big[\tilde W(e_i) \tilde W(e_j)\big]
=
\ps{e_i}{(\I_{\theta_0}^* \II_\ve \I_{\theta_0})^{-1}e_j}_\VV 
.
$$
Then, for any $h=\sum_j c_j e_j\in \DD^{\kappa}$, define $\tilde W(h)$ by linearity as $\tilde W(h) \equiv \sum_j c_j \tilde W(e_j)$. Then, $\tilde W(h)$ is centred for all $h\in \DD^\kappa$, and for all $f=\sum_j c_j e_j\in \DD^\kappa$ and $g=\sum_j c_j' e_j\in \DD^\kappa$, we have
$$
\E\big[ \tilde W(f) \tilde W(g) \big]
=
\sum_{i,j} \lambda_i \lambda_j' \ps{e_i}{(\I_{\theta_0}^* \II_\ve \I_{\theta_0})^{-1} e_j}_\VV
=
\ps{f}{(\I_{\theta_0}^* \II_\ve \I_{\theta_0})^{-1} g}_\VV
.
$$
One thus obtain $\tilde G$ in (\ref{defGbis}) by taking $\tilde\gamma_k \equiv \tilde W(e_k)$. Indeed, we will see in Section \ref{sec:Support} below that $\tilde G$ is a centred $\DD^{-\beta}$-valued Gaussian random variable for $\beta>\kappa+\alpha$ provided $\lambda_k\gtrsim k^{1/\alpha}$ for some $\alpha>0$, so that  
$$
\langle \tilde G , h\rangle_\VV
=
\tilde W(h)
,\spa 
\forall\ h\in\DD^\beta
,
$$
hence $\tilde G$ has the right covariance.

\end{Example}


\subsection{Support of the process}\label{sec:Support}

Note that there always exists a Banach space $(\X, \|\cdot\|_\X)$ such that the formal series $G$ from~(\ref{LimitProc}) converges in $\X$ almost surely. When $\dim\HH<\infty$, then $G$ is obviously supported on $\HH$. When $\dim\HH=\infty$, define for any sequence $w=(w_k)\in \R^\infty$, with $w_k>0$, the separable Hilbert space 
\begin{equation}\label{Hweighted}
\HH_w 
\equiv 
\overline{
\Big\{ 
\sum_{k = 1}^n a_k h_k : a_1,\ldots, a_n\in\R,\ n\ge 1\Big\}
}
,
\end{equation}
where the completion is taken with respect to the norm 
$$
\Big\|\sum_{k=1}^n a_k h_k\Big\|^2_{\HH_w}
\equiv 
\sum_{k=1}^n w_k^2 a_k^2 
.
$$
Then, the series defining $G$ converges almost surely in $\HH_w$, hence defines a $\HH_w$-valued random variable, if and only if $w\in \ell^2(\N)$. Indeed, for $G_n=\sum_{j=1}^n \gamma_j h_j$ the partial sums of $G$, we have for all $n,p\in\N$
$$
\E\big[\|G_{n+p}-G_n\|^2_{\HH_w}\big]
=
\E\Big[\sum_{k=n+1}^{n+p} w_k^2 \gamma_k^2 \Big]
=
\sum_{k=n+1}^{n+p} w_k^2
.
$$
If $(w_k)\in \ell^2(\N)$, then Lévy's maximal inequality (see, e.g., Exercise 2.6.8 and Theorem 3.1.11 in \citet*{GinNic2016}) entails that, for all $k\in\N$ and $\ve>0$
$$
\P\big[\max_{p\le k} \|G_{n+p}-G_n\|^2_{\HH_w} > \ve\big]
\le 
2\P\big[\|G_{n+k}-G_n\|^2_{\HH_w} > \ve\big]
\le 
\frac{2}{\ve}\ \E\big[\|G_{n+k}-G_n\|^2_{\HH_w}\big]
$$
We deduce from the continuity of measures that 
$$
\limsup_{n\to\infty} \P\big[\sup_p \|G_{n+p}-G_n\|^2_{\HH_w} > \ve\big]
=
\limsup_{n\to\infty} \lim_{k\to\infty} \P\big[\max_{p\le k} \|G_{n+p}-G_n\|^2_{\HH_w} > \ve\big]
=
0
.
$$
Consequently, Lemma 9.2.4 in \citet*{Dudley2002} entails that the random sequence $(G_n)$ in $\HH_w$ converges almost surely in $\HH_w$. Necessity of $(w_k)\in \ell^2(\N)$ follows from Fernique's theorem.\smallskip 

Although this shows that a separable Hilbert space supporting the law of $G$ always exists, this does not provide an operational way to determine such a space in practice unless the metric $\ps{\cdot}{\cdot}_\HH$ is available to us. Instead, we would like to have an effective way to understand what spaces support the law of $G$, relying only on what $\HH$ is as a set. We illustrate this in the next example, in the setting of Example \ref{ex:HandSBis} and Example \ref{ex:SeriesG}.

\begin{Example}\label{ex:SuppHilbert}
Consider the setting of Example \ref{ex:SeriesG}, where $\|\cdot\|_\HH\asymp \|\cdot\|_{\DD^{-\kappa}}$ and $\|\cdot\|_\SS\asymp \|\cdot\|_{\DD^\kappa}$, and recall the formal series 
$
\tilde G 
=
\sum_k \tilde \gamma_k e_k
$,
where $\E[\tilde \gamma_k^2]=\|e_k\|^2_\SS$; see (\ref{eq:CompatS}). Reasoning as above yields that $\tilde G$ converges almost surely in the metric 
$
\|\sum_j u_j e_j\|^2_w 
\equiv 
\sum_j w_j^2 u_j^2
$ if and only if $\sum_j w_j^2 \|e_j\|^2_\SS<\infty$. Since $\|\cdot\|_\SS\asymp \|\cdot\|_{\DD^\kappa}$ and $\|e_j\|^2_{\DD^\kappa}=\tau_j^\kappa$, the latter condition is equivalent to 
$
\sum_j w_j^2 \tau_j^\kappa
<
\infty 
$. Taking $w_j = \tau_j^{-\beta/2}$ thus entails that $\tilde G$ is a $\DD^{-\beta}$-valued random variable if and only if 
$\sum_j \tau_j^{\kappa-\beta}<\infty.$
If we further assume that $\tau_j \gtrsim j^{1/\alpha}$ for some $\alpha>0$, then the previous series converges as soon as $\beta>\kappa+\alpha$, and this also becomes necessary if the reverse inequality $\tau_j\lesssim j^{1/\alpha}$ holds. In the setting where $(e_j, \tau_j)$ are the eigenpairs of the Dirichlet Laplacian on some regular domain $\ZZ\subset \R^n$ (see after Example \ref{ex:HandSBis}), yielding $\tau_j\asymp j^{2/n}$, we deduce that $\tilde G$ is a $\DD^{-\beta}$-valued random variable if and only if $\beta>\kappa+n/2$. 
\end{Example}

For general (separable) Banach spaces $\X$, determining when $\X$ supports the law of $G$ is more delicate as the geometry of such spaces is more flexible, hence can also be more complicated, than for Hilbert spaces. 
If $(\X, \|\cdot\|_\X)$ is a Banach space containing $\HH$ as a dense subset such that the canonical injection $(\HH, \|\cdot\|_\HH)\embed (\X, \|\cdot\|_\X)$ is continuous, and denoting by $B^*_1$ the unit ball of $(\X^*_{\overline\VV_0}, \|\cdot\|_{\X^*_{\overline\VV_0}})$ from (\ref{eq:B*B*S}), then the law of $G$ is supported on $\X$ provided 
$$
\int_0^1 
\sqrt{\log N(B^*_1, \|\cdot\|_\SS, \ve)}\, d\ve 
< 
\infty,
$$
where $N(B^*_1, \|\cdot\|_\SS, \ve)$ stands for the $\ve$-covering number of $B^*_1$ in $\|\cdot\|_\SS$-metric. This follows from observing that the process $W$ from (\ref{CovW}) is the white-noise process on $\SS$ by virtue of (\ref{SHilbert}) and that the diameter of $B_1^*$ in the covariance metric is finite by continuity of the inclusion 
$\X^*_{\overline \VV_0}\embed \SS$
(see (\ref{allEmbeddings})). The result thus follows from Dudley's theorem applied to the metric space $(B^*_1, \|\cdot\|_\SS)$ and sub-Gaussian process $\{W(b) : b\in B^*_1\}$; see, e.g., Theorem 2.3.7 in \citet*{GinNic2016}.

	\section{Asymptotic minimax efficiency lower bound}\label{sec:GeneralMinimax}

    For a appropriate loss functions $\rho:\B\to [0,\infty)$, and under regularity assumptions, the minimax theorem (see Theorem 3.12.5 in \cite{VdVWel1996}) establishes a lower bound on the local asymptotic minimax risk for estimating a functional $F(\theta)$ taking values in some Banach space $(\B, \|\cdot\|_\B)$, given by
    \begin{equation}\label{MinimaxRisk}
    \RRR_{\theta_0}(\rho, F, \VV_0) 
    \equiv 
    \inf_{(\hat F_N)} 
    \sup_{\substack{J\subset \VV_0\\ |J|<\infty}}
    \liminf_{N\to \infty} \max_{h\in I} E^N_{\theta_0+\frac{h}{\sqrt{N}}}
    \rho\bigg( \sqrt{N}\Big(\hat F_N-F\big[\theta_0+h/\sqrt{N}\big]\Big) \bigg)
    ,
    \end{equation}
    where the infimum ranges over all sequences $(\hat F_N)_N$ of Borel measurable maps $\hat F_N : (\XX\times \R^p)^N\to \B$, and the supremum ranges over all subsets $J\subset \VV_0$ with finite cardinality $|J|$. The lower bound is usually established for the quantity 
    $$
    \tilde\RRR_{\theta_0}(\rho, F, \VV_0) 
    \equiv
    \inf_{\hat\theta_N} 
    \limsup_{c\to\infty}
    \liminf_{N\to \infty} 
    \sup_{\substack{h\in\VV_0\\ \|h\|\le c/\sqrt{N}}}
    E^N_{\theta_0+h}\Big[\rho\Big( \sqrt{N}\big(\hat\theta_N-F[\theta_0+h]\big) \Big)\Big]
    ,
    $$
    for convenience but it is straightforward to see that $\tilde\RRR_{\theta_0}(\rho, F, \VV_0)\ge \RRR_{\theta_0}(\rho, F,\VV_0)$.


\subsection{Non-parametric lower bound for differentiable functionals}\label{sec:MinimaxBound}

Let us now introduce the functionals of the parameter space for which we wish to establish an asymptotic efficiency lower bound. Formally, consider (not necessarily linear) functionals $F(\theta)$ of parameters $\theta\in\Theta$, where $F:\Theta\mapsto \B$ for some Banach space $(\B, \|\cdot\|_\B)$. We require that these functionals be differentiable at $\theta_0$ along any direction from the tangent space $\VV_0$ at $\theta_0$, with continuous linearization $\dot{F}: (\VV_0,\ps{\cdot}{\cdot}_\VV)\mapsto \B$.

\begin{Condition}\label{CondRegularParameters}
	Fix $\theta_0\in\Theta$ with tangent space $ \VV_0$ and linearization $\I_{\theta_0}$ satisfying Condition \ref{CondIntroLinModel} and Condition \ref{CondIntroInj}. Let $(\HH, \ps{\cdot}{\cdot}_\HH)$ be the Hilbert space defined in (\ref{VVEmbedH}) and containing $\VV_0$ as a subspace. Let $(\B,\|\cdot\|_\B)$ be a Banach space and $F:\Theta\mapsto \B$ a measurable map. Assume that there exists a linear and continuous map $\dot{F}:(\VV_0,\ps{\cdot}{\cdot}_{\VV})\mapsto (\B,\|\cdot\|_\B)$ such that, as $s\to 0$,
	$$
	\frac{F(\theta_0+sh)-F(\theta_0)}{s}
	\to 
	\dot{F}[h]
	\quad
	\textrm{in }\ \B
	,\spa \forall\ h\in  \VV_0
	.
	$$
\end{Condition}

Note that under Condition \ref{CondIntroLinModel} the quantity $F(\theta_0+sh)$ is well-defined as soon as $|s|\le \delta(h)$. The lower bound on the asymptotic local minimax risk for estimating $F(\theta_0)$ in Theorem \ref{TheorGeneralMinimax} below involves the law of the centred Gaussian process $\dot F[G]$, where we recall that $G$ is defined in (\ref{LimitProc}) as the formal series $G=\sum_{j} \gamma_j h_j$ ranging over $j\le \dim \HH$. Similarly, we define $\dot F[G]$ as the formal series
\begin{equation}\label{DefFG}
    \dot F[G]
    \equiv 
    \sum_{j} \gamma_j \dot F[h_j]
    .
\end{equation}
Note that $\dot F[h_j]$ is well-defined since $h_j\in \VV_0$ for all $j$. If $\B$ is a Banach space and $F:\Theta\mapsto \B$ a measurable map satisfying Condition \ref{CondRegularParameters} are such that the series in (\ref{DefFG}) converges almost surely in $\B$, then the law of $\dot F[G]$ is described as follows. If we further assume that $\dot F : (\VV_0, \ps{\cdot}{\cdot}_\HH)\to \B$ is  continuous (which will turn out to be necessary, see Theorem \ref{TheorGeneralMinimax} below), then $\dot F$ admits a (continuous) adjoint 
$$
\dot{F}^*_{\rm LAN} : \B^*\to (\HH, \ps{\cdot}{\cdot}_{\HH})
,
$$
since $\VV_0$ is dense in $\HH$, by construction. The adjoint is characterized by 
\begin{equation}\label{adjFLan}
\big\langle \dot F^*_{\rm LAN}[\ell] , h\big\rangle_{\HH}
\equiv 
\ell(\dot F[h]) 
,\spa 
\forall\ \ell\in\B^*, h\in\HH 
.
\end{equation}
Then, $\dot F[G]$ is a centred $\B$-valued random variable with covariance given, for all $\ell,\ell'\in\B^*$, by
\begin{align*}
\E\big[\ell(\dot F[G])\ell'(\dot F[G])\big]
&=
\sum_{j} \ell(\dot F[h_j])\ell'(\dot F[h_j]) 
\\[2mm]
&=
\sum_{j} \big\langle\dot F^*_{\rm LAN}[\ell] , h_j\big\rangle_{\HH}
\big\langle\dot F^*_{\rm LAN}[\ell'], h_j\big\rangle_{\HH}
\\[2mm]
&=
\big\langle\dot F^*_{\rm LAN}[\ell] , \dot F^*_{\rm LAN}[\ell']\big\rangle_\HH
,
\end{align*}
by Parseval's formula in the Hilbert space $(\HH, \ps{\cdot}{\cdot}_\HH)$ with orthonormal basis $(h_j)$. It is instructive to express the covariance above when the adjoint of $\dot F$ is computed in the space $(\VV_0, \ps{\cdot}{\cdot}_\VV)$ instead. Indeed, for continuous and linear $\dot F : (\VV_0, \ps{\cdot}{\cdot}_\VV)\to \B$, there is an adjoint 
$$
\dot F^* : \B^* \to (\overline \VV_0,\ps{\cdot}{\cdot}_\VV)  
,
$$
characterized by 
\begin{equation}\label{AdjF}
\big\langle \dot F^*[\ell], h \big\rangle_\VV 
\equiv 
\ell(\dot F[h])
,\spa 
\forall\ \ell\in \B^*, h\in \overline \VV_0 
.
\end{equation}
Equating the l.h.s. of (\ref{adjFLan}) and (\ref{AdjF}), and recalling the definition of $\ps{\cdot}{\cdot}_\HH$ from (\ref{defHInnerProd}) provides
$$
\big\langle (\I_{\theta_0}^* \II_\ve \I_{\theta_0})[\dot F^*_{\rm LAN}[\ell]] , h \big\rangle_\VV 
=
\big\langle \dot F^*[\ell], h \big\rangle_\VV
,\spa 
\forall\ \ell\in\B^*, h\in \VV_0 
.
$$
It follows that in order to define the process $\{Z(\ell) : \ell\in \B^*\}$ with covariance 
\begin{align}\label{eq:covF}
\E\big[Z(\ell)Z(\ell')\big]
&\equiv 
\langle \dot F^*_{\rm LAN}[\ell], \dot F^*_{\rm LAN}[\ell']\rangle_\HH
,\spa 
\ell, \ell'\in \B^* 
\\[2mm]
\nonumber
&=
\big\langle \dot F^*[\ell], (\I_{\theta_0}^* \II_\ve \I_{\theta_0})^{-1}\big[\dot F^*[\ell']\big] \big\rangle_\VV
,
\end{align}
as above, it is necessary and sufficient that there exists a solution $\bar\ell\in\HH$, for any $\ell\in\B^*$, to the information equation 
\begin{equation}\label{InfoEq}
(\I_{\theta_0}^* \II_\ve \I_{\theta_0})[\bar \ell]
=
\dot F^*[\ell]
.
\end{equation}
It follows from Theorem \ref{TheorI*IHomeo} characterizing the range of the information operator featuring in~(\ref{InfoEq}), that such $\bar\ell\in\HH$ exists precisely when $\dot F^*[\ell]\in \SS$, where we recall that $\SS$ is the $\overline \VV_0$-dual of $\HH$ defined in (\ref{DefS}). It is straightforward to see that this is equivalent to continuity of $\dot F:(\VV_0, \ps{\cdot}{\cdot}_\HH)\to \B$, which is easier to establish since it does not require to compute adjoints.\smallskip

Recall that a Banach space $\B$ is reflexive if its second dual is isormophic to $\B$, \ie any continuous linear map $f:(\B^*,\|\cdot\|_{\B^*})\to \R$ is of the form $f(\ell)=\ell(x)$ for some $x\in\B$ and all $\ell\in \B^*$, where $\|\cdot\|_{\B^*}$ is the operator norm on $\B^*$. Examples of such spaces include $L^p$-spaces, $1<p<\infty$, any Hilbert space, or Sobolev spaces $W^{k,p}(U)$ on an open set $U\subset\R^n$, with $k\in\N$ and $1<p<\infty$. We say that a map $\rho : \B\to \R$ is sub-convex if the set $\{y\in \B : \rho(y)\le c\}$ is $\|\cdot\|_\B$-closed, convex, and symmetric for all $c\in\R$. We say that $\rho$ is coercive if $\rho(x)\to \infty$ as $\|x\|_\B\to \infty$. We can now state the main result of this section, the first part of which is a rephrasing of Theorem 3.12.5 in \citet*{VdVWel1996} in our setting. The proof is deferred to Section~\ref{sec:ProofMinimax} below.

\begin{Theor}\label{TheorGeneralMinimax}
Let $\theta_0\in\Theta$ with tangent space $ \VV_0$ and linearization $\I_{\theta_0}$ satisfying Condition \ref{CondIntroLinModel} and Condition \ref{CondIntroInj}. Let $(\B,\|\cdot\|_\B)$ be a Banach space and $F : \Theta \to \B$ a measurable map, with linearization $\dot{F}: \VV_0\to \B$ at $\theta_0$ satisfying Condition \ref{CondRegularParameters}. Let $\rho:\B\to [0,\infty)$ be a sub-convex map. (i) If there exists a tight centred $\B$-valued Gaussian random variable $Z$ with covariance as in~(\ref{eq:covF}), and if $\dot F$ is continuous from $(\VV_0, \ps{\cdot}{\cdot}_\HH)$ to $(\B, \|\cdot\|_\B)$, then the asymptotic local minimax risk from (\ref{MinimaxRisk}) for estimating $F(\theta_0)$ with tangent space $\VV_0$ in $\rho$-loss is lower bounded as 
$$
\RRR_{\theta_0}(\rho, F, \VV_0)
\geq 
\E\big[\rho(Z)\big]
,
$$
and the risk is infinite if the r.h.s. is infinite. (ii) If $\dot F$ is not continuous from $(\VV_0, \ps{\cdot}{\cdot}_\HH)$ to $(\B, \|\cdot\|_\B)$, and $\rho$ is coercive, then the risk above is infinite. (iii) If no variable $Z$ as in Part (i) exists, and if $\B$ is separable and reflexive, and $\rho$ is coercive, then the risk above is infinite. 
\end{Theor}

In statistical models satisfying standard regularity conditions, Theorem \ref{TheorGeneralMinimax} provides a lower bound on the asymptotic local minimax risk for estimating a differentiable functional $F(\theta)$ of the natural parameter $\theta\in\Theta$ as soon as the centred process (\ref{eq:covF}) exists in some Banach space. 

\begin{Example}\label{ExFId}
For $F\equiv {\rm Id}_{\Theta}$, the identity map on $\Theta$, let $(\X,\|\cdot\|_\X)$ be a Banach space supporting the law of the process $G=\sum_k \gamma_k h_k$ from (\ref{LimitProc}), or equivalently $\{W(f) : f\in \SS\}$ in (\ref{CovW}). We necessarily have a continuous embedding $\HH\embed \X$ since $\HH$ is the RKHS of $G$. Consequently, $G$ has covariance given by (\ref{eq:covF}) (see Section \ref{sec:CovSeries}) and $F\equiv {\rm Id}_{\Theta}$ satisfies Condition \ref{CondRegularParameters} with Banach space $\B\equiv \X$. Therefore, Theorem \ref{TheorGeneralMinimax} yields 
$$
\inf_{(\hat \theta_N)} 
\sup_{\substack{J\subset \VV_0\\ |J|<\infty}}
\liminf_{N\to \infty} \max_{h\in J} E^N_{\theta_0+\frac{h}{\sqrt{N}}}\rho\bigg( \sqrt{N}\Big(\hat \theta_N-\theta_0-\frac{h}{\sqrt{N}} \Big)\bigg)
\ge 
\E\big[\rho(G)\big]
,
$$
where the infimum ranges over all Borel measurable maps $\hat \theta_N : (\XX\times \R^p)^N\to \X$, and the supremum ranges over all subsets $J\subset \VV_0$ with finite cardinality $|J|$. 
\end{Example} 

\begin{Example} \label{ex:ContBanach}
For a Banach space $(\X,\|\cdot\|_\X)$ supporting the law of $G$ as in Example~\ref{ExFId}, if a map $F:\Theta\mapsto \B$ satisfies Condition \ref{CondRegularParameters} with continuous linearization $\dot F : (\HH, \ps{\cdot}{\cdot}_\HH)\to \B$ that further extends to a continuous operator $\dot F : (\X,\|\cdot\|_\X)\to (\B,\|\cdot\|_\B)$, then the composition $Z\equiv \dot F[G]$ defines a centred $\B$-valued Gaussian random variable with covariance (\ref{eq:covF}). Indeed, it is clear that $\dot F[G]$ is centred, $\B$-valued since $G$ is $\X$-valued, and Gaussian since $G$ is Gaussian and $\dot F$ is linear. Given the covariance of $G$ from (\ref{CovW}) and previous calculations, one checks directly that $\dot F[G]$ has covariance (\ref{eq:covF}), and $\dot F[G]$ is Borel measurable as the composition of the continuous map $\dot F$ with the Borel measurable random variable $G$.
\end{Example}


Note that the lower bound $\E[\rho(Z)]$ in Theorem \ref{TheorGeneralMinimax} may be infinite for some sub-convex maps $\rho$ even if $Z$ is a centred $\B$-valued Gaussian random variable. When $\rho(x)=\exp(\|x\|^2/2\alpha^2)$ for all $x$ in a separable Banach space $\B$, then the lower bound is infinite if and only if $\alpha\le \|\dot F\|_{\HH\to\B}$, the operator norm of $\dot F : \HH\to \B$, by virtue of Corollary 2.2.9 in \citet*{GinNic2016} since
$$
\sigma^2
\equiv 
\sup_{\substack{\ell\in \B^*\\ \|\ell\|_{\B^*}\le 1}} \E[\ell(Z)^2]
=
\sup_{\substack{\ell\in \B^*\\ \|\ell\|_{\B^*}\le 1}} \|\dot F^*_{\rm LAN}[\ell]\|^2_\HH 
=
\sup_{\substack{h\in\VV_0\\ \|h\|_{\HH}\le 1}}
\sup_{\substack{\ell\in \B^*\\ \|\ell\|_{\B^*}\le 1}} |\ell(\dot F[h])|^2
=
\|\dot F\|^2_{\HH\to\B}
.
$$


\begin{Example}\label{CoroMinimaxNorm}
Any convex map $\rho:(\B,\|\cdot\|_\B)\to\R$ is sub-convex, so in particular the power norm $\|\cdot\|^\alpha_\B$ is a sub-convex map for any $\alpha>0$. Let $\theta_0\in\Theta$ with tangent space $ \VV_0$ and linearization $\I_{\theta_0}$ satisfying Condition \ref{CondIntroLinModel} and Condition \ref{CondIntroInj}. Let $(\B,\|\cdot\|_\B)$ be a separable and reflexive Banach space, and $F : \Theta \to \B$ a measurable map with continuous linearization $\dot{F}: (\VV_0,\ps{\cdot}{\cdot}_\HH)\to \B$ at $\theta_0$ satisfying Condition \ref{CondRegularParameters}. 
    The local asymptotic minimax risk in $\|\cdot\|^\alpha_\B$-loss, $\alpha>0$, is lower bounded as
    $$
	\RRR_{\theta_0}(\|\cdot\|^\alpha_{\B}, F, \VV_0)
	\geq 
	\E\big[\|Z\|^\alpha_{\B}\big]
    ,
	$$
    for a centred $\B$-valued Gaussian random variable $Z$ with covariance as in (\ref{eq:covF}) as soon as such a variable exists. Otherwise, the risk is infinite.
\end{Example}

\begin{Remark}
    Once the regularity of the model has been established in some norm $\|\cdot\|_\B$, then the minimax theorem applies to the $\|\cdot\|_\B$-loss as well as to any $\|\cdot\|_{\B'}$-loss corresponding to a norm weaker than $\|\cdot\|_\B$, \ie such that $\|\cdot\|_{\B'}\lesssim \|\cdot\|_\B$ on $\B$.
\end{Remark}

\subsection{Semi-parametric lower bound}\label{sec:MinimaxSemiPar}

The results of Section \ref{sec:MinimaxBound} provide an alternative way to establish a lower bound on the asymptotic risk from (\ref{MinimaxRisk}) in the classical semi-parametric setting where the functional $F(\theta)$ to be estimated is of the form $F(\theta)=\ps{\psi}{\theta}_\VV$ for some $\psi\in \overline \VV_0$, associated with (reflexive and separable) Banach space $\B\equiv \R$ and Euclidean loss $\rho(t)=|t|$, $t\in\R$. For such functionals, it has been established, first in \citet*{VdV1991} (see Theorem 4.1 there), that the local asymptotic minimax risk at $\theta=\theta_0$ in Euclidean loss admits a finite lower bound if and only if $\psi$ belongs to the range $R(\I_{\theta_0}^*)$ of $\I_{\theta_0}^*$, \ie $\psi=\I_{\theta_0}^*[\bar \psi]$ for some $\bar\psi\in \overline \VV_0$. In view of the equality ${\rm Ker}(\I_{\theta_0})=\overline{R(\I_{\theta_0}^*)}$, this means that unless $R(\I_{\theta_0}^*)$ is closed in $\overline \VV_0$, the asymptotic local risk will be infinite for some $\psi\in \overline \VV_0$, even when $\I_{\theta_0}$ is injective over $\overline \VV_0$. Determining which $\psi$'s lead to a finite lower bound on the asymptotic risk thus becomes a problem in its own right, independently of the injectivity of~$\I_{\theta_0}$, and needs to be addressed on a case-by-case basis. Although we do not `solve' this question, we can use the results of Section \ref{sec:FisherInfo} to characterize exactly when $\psi\in R(\I_{\theta_0}^*)$ in a way that proves much easier in some infinite-dimensional settings, as will be demonstrated in Section \ref{sec:Appli} below.

\begin{Prop}\label{PropSemiPar}
Let $\theta_0\in\Theta$ with tangent space $ \VV_0$ and linearization $\I_{\theta_0}$ satisfying Condition~\ref{CondIntroLinModel} and Condition~\ref{CondIntroInj}. For any $\psi\in\overline\VV_0$ and $\theta\in\Theta$, let $F_\psi(\theta)\equiv \ps{\psi}{\theta}_{\VV}$. Let $\SS$ be defined as in (\ref{DefS}). If $\psi\in \SS$, then the asymptotic local minimax risk from (\ref{MinimaxRisk}) for estimating $F_\psi(\theta_0)$ with tangent space $\VV_0$ in Euclidean $|\cdot|^2$-loss is lower bounded as
$$
\RRR_{\theta_0}(|\cdot|^2, F_\psi, \VV_0)
\ge 
\|\psi\|^2_\SS
=
\ps{\psi}{(\I_{\theta_0}^* \II_\ve \I_{\theta_0})^{-1}[\psi]}_\VV
,
$$
and the risk above is infinite if $\psi\in \overline \VV_0\sm \SS$. 
\end{Prop}

\begin{Proof}{Proposition \ref{PropSemiPar}}
    For any $\psi\in\SS$, the map $F_\psi(\theta)=\ps{\psi}{\theta}_\VV$ satisfies Condition \ref{CondRegularParameters} with continuous linearization $\dot F_\psi: (\VV_0, \ps{\cdot}{\cdot}_\VV)\to \B\equiv \R$ given by $\dot F_\psi[h] = \ps{\psi}{h}_\VV$. It is also continuous from $(\VV_0, \ps{\cdot}{\cdot}_\HH)$ to $\R$ since $\lvert\ps{\psi}{h}_\VV\rvert \le \|\psi\|_\SS \|h\|_\HH$ for all $h\in \VV_0$ by definition of $\SS$ (see (\ref{DefS})). It follows from (\ref{AdjF}) that 
    $$
    \langle\dot F_\psi^*[a] , h\rangle_\VV
    =
    a \dot F[h]
    =
    \ps{a\psi}{h}
    ,\spa 
    \forall\ a\in \R(=\B^*),\ h\in \overline\VV_0 
    ,
    $$
    \ie $\dot F_\psi^*[a] = a\psi$ for all $a\in\R$. It follows that the $\B$-valued random variable $Z$ from Theorem \ref{TheorGeneralMinimax}(i) is nothing but the centred real-valued random variable with variance $\ps{\psi}{(\I_{\theta_0}^* \II_\ve \I_{\theta_0})^{-1}[\psi]}_\VV$. 
    Theorem \ref{TheorGeneralMinimax}(i) thus yields 
    $$
    \RRR_{\theta_0}(|\cdot|^2, F_\psi, \VV_0)
    \ge
    \E\big[ |Z|^2 \big]
    =
    \ps{\psi}{(\I_{\theta_0}^* \II_\ve \I_{\theta_0})^{-1}[\psi]}_\VV
    ,
    $$
    and the lower bound is precisely $\|\psi\|^2_\SS$ by virtue of (\ref{SHilbert}). When $\psi\notin\SS$, then the map $\dot F : (\VV_0, \ps{\cdot}{\cdot}_\HH)\to \R$ is not continuous, by definition of $\SS$ (see (\ref{DefS})), so that Theorem \ref{TheorGeneralMinimax}(ii) entails that the local asymptotic minimax risk is infinite.
\end{Proof}


\subsection{Proof of Theorem \ref{TheorGeneralMinimax}}\label{sec:ProofMinimax}

\begin{Proof}{Theorem \ref{TheorGeneralMinimax}}
For the first part, we use Theorem 3.12.5 in \citet*{VdVWel1996} and the terminology employed in Section 3.12 there, by first checking the conditions to apply the theorem. Recall from (\ref{model}) that the probability measures $\{P_{\theta}^N : \theta\in\Theta\}$ are defined on the Borel $\sigma$-field $\BB_N$ of $\XX_N\equiv (\XX\times \R^p)^N$. Proposition \ref{PropLAN} entails that the sequence of experiments 
$$
\Big\{
\big(\XX_N, \BB_N, P^N_{\theta_0+h/\sqrt{N}}
\big) : h\in \VV_0\Big\}_{N\geq 1}
$$ 
satisfies an asymptotic locally normal expansion for the Hilbert space norm given by the LAN norm $\|\cdot\|_{\rm LAN}$ from (\ref{eq:LANnorm}). Condition \ref{CondIntroInj} and Proposition \ref{PropInvertIepsilon} entail that $\|\cdot\|_{\rm LAN}$ is indeed a norm on the subspace $\VV_0$ of the Hilbert space $(\HH, \ps{\cdot}{\cdot}_\HH)$ from (\ref{eq:VVdomH})-(\ref{VVEmbedH}), and recall that, by definition of $\HH$, the $\|\cdot\|_{\rm LAN}$ and $\|\cdot\|_\HH$ norms coincide on $\VV_0$.
Condition \ref{CondRegularParameters} entails that the sequence of parameters $(F(\theta_0+h/\sqrt{N}))_{N\geq 1}$ is regular for any $h\in \VV_0$, with norming operators $r_N:\B\to\B, x\mapsto \sqrt{N}x$ for all $N\geq 1$, and continuous linearization $\dot{F}:(\VV_0, \ps{\cdot}{\cdot}_\HH)\to \B$.

For $k\le \dim \HH$, define 
$$
\VV_0^{(k)}\equiv {\rm span}\{h_1,\ldots, h_k\}
,
$$ 
where $(h_k)$ is the orthonormal basis of $(\HH, \ps{\cdot}{\cdot}_\HH)$ used in the series representation (\ref{LimitProc}) of $G$, and let 
$$
G_k 
\equiv 
\sum_{j=1}^k \gamma_j h_j
,
$$
where the $\gamma_j$'s are i.i.d. $\NN(0,1)$-distributed. Then $G_k$ defines a random element in $\VV_0$ since each $h_j\in\VV_0$. The proof of \citet*{VdVWel1996} for sub-convex loss $\rho$ proceeds by constructing approximations $(\rho_r)$ of $\rho$ of the form 
$$
\rho_r(y)
=
\sum_{i=1}^r \1_{K_i^c}(\ell_{i,1}(y), \ldots, \ell_{i,p_i}(y)) 
,
$$
for some $\ell_{i,j}\in\B^*$, $p_i\in\N$, and $K_i\subset \R^{p_i}$ compact, convex and symmetric, in such a way that $0\le \rho_r\le \rho$ pointwise on $\B$, and by showing that, for any $k\ge 1$, the asymptotic local minimax risk $\RRR(\ell_r, F, \VV_0^{(k)})$ of the submodel $\VV_0^{(k)}\subset \VV_0$ is lower bounded as
$$
\RRR_{\theta_0}(\rho_r, F, \VV_0^{(k)}\big)
\ge 
\E\big[ \rho_r(\dot F[G_k])\big]
,
$$
for all $r\ge 1$ and $k\le\dim\HH$. Since $\VV_0\supset \VV_0^{(k)}$ and $\ell\ge \ell_r$, we have
$$
\RRR_{\theta_0}(\rho, F, \VV_0) 
\ge
\RRR_{\theta_0}(\rho, F, \VV_0^{(k)}\big)
\ge 
\RRR_{\theta_0}(\rho_r, F, \VV_0^{(k)}\big)
,
$$
for all $r\ge 1$ and $k\le \dim\HH$.  We deduce that
\begin{equation}\label{eq:Lrbound}
\RRR_{\theta_0}(\rho, F, \VV_0) 
\ge
\E\big[ \rho_r(\dot F[G_k])\big]
=
\sum_{i=1}^r \P\Big[\big(\ell_{i,1}(\dot F[G_k]), \ldots, \ell_{i,p_i}(\dot F[G_k])\big)\in K_i^c\Big] 
,
\end{equation}
for all $r\ge 1$ and $k\le\dim\HH$. We may now proceed with the proof of (i). \smallskip 

(i) Let us now assume that there exists a centred $\B$-valued Gaussian random variable $Z$ with covariance as in (\ref{eq:covF}). Let us first note that the sequence $(\rho_r)$ constructed in the proof of Theorem 3.12.5 in \citet*{VdVWel1996} is such that $\rho_r(Z)\uparrow \rho(Z)$ almost surely. Now, for any $\ell\in \B^*$, observe that 
$$
\sigma_{\ell,k}^2
\equiv 
\E\big[\ell(\dot F[G_k])^2\big]
=
\sum_{j=1}^k \ell(\dot F[h_j])^2 
\stackrel{k\uparrow\dim\HH}{\longrightarrow} 
\sum_{j\ge 1} \ell(\dot F[h_j])^2  
=
\E\big[\ell(Z)^2]
\equiv 
\sigma_\ell^2 
, 
$$
so that
$$
\E\big[ e^{it \ell(\dot F[G_k])}\big]
=
e^{-\frac{t^2}{2}\sigma_{\ell,k}^2} 
\stackrel{k\uparrow\dim\HH}{\longrightarrow}
e^{-\frac{t^2}{2}\sigma_\ell^2}
=
\E\big[ e^{it \ell(Z)}\big]
,\spa 
\forall\ t\in\R 
.
$$
We deduce that $\ell(\dot F[G_k])$ converges in distribution to $\ell(Z)$ as $k\to\dim\HH$, for any $\ell\in \B^*$. By taking arbitrary linear combinations, we deduce that
$$
\big(\ell_1(\dot F[G_k]),\ldots, \ell_n(\dot F[G_k])\big) 
\leadsto
\big(\ell_1(Z),\ldots, \ell_n(Z)\big) 
,
$$
as $k\to\dim\HH$, for any $\ell_1, \ldots, \ell_n\in\B^*$ and $n\in\N$. Then, the Portmanteau theorem (see, e.g., Theorem 1.3.4(ii) in \citet*{VdVWel1996}) provides
$$
\RRR_{\theta_0}(\rho, F,\VV_0) 
\ge 
\sum_{i=1}^r \P\Big[\big(\ell_{i,1}(Z), \ldots, \ell_{i,p_i}(Z)\big)\in K_i^c\Big]
=
\E\big[\rho_r(Z)\big]
,
$$
where the r.h.s. in the last display is finite (and bounded by $r$) for any $r\ge 1$. Since $\rho_r(Z)\uparrow \rho(Z)$ almost surely, the monotone convergence theorem yields 
$$
\RRR_{\theta_0}(\rho, F,\VV_0) 
\ge 
\E\big[\rho(Z)\big]
,
$$
even if the r.h.s. in the last display is infinite (in which case the l.h.s. is also infinite).\medskip

(ii) Assume that $\dot F:\VV_0\to \B$ is not continuous from $(\VV_0, \ps{\cdot}{\cdot}_\HH)$ to $(\B, \|\cdot\|_\B)$. It follows from the proof of Theorem 3.12.5 in \citet*{VdVWel1996} that $\rho_r(\dot F[G_k])\uparrow \rho(\dot F[G_k])$ almost surely as $r\to\infty$ for all $k$ (the complement of the uncountable intersection $C$ in the construction of $\rho_r$ there, resulting from applying the Hahn-Banach theorem, can be reduced to a countable one after intersecting it with $S$ taken, instead, to be the $\|\cdot\|_\B$-closure of ${\rm span}\{\dot F[h_j] : j\ge 1\}$ which is a separable closed subspace of $\B$ and supports all the $\dot F[G_k]$'s).
Taking the limit as $r\to\infty$ in (\ref{eq:Lrbound}) yields, by the monotone convergence theorem,
\begin{equation}\label{eq:Lrbound2}
\RRR_{\theta_0}(\rho, F, \VV_0) 
\ge
\E\big[ \rho(\dot F[G_k])\big]
,\spa 
\forall\ k \ge 1
.
\end{equation}
Assume, \emph{ad absurdum}, that 
\begin{equation}\label{asmpFiniteE}
\limsup_{n\to\infty} \E\big[ \rho(\dot F[G_n]) \big]
<
\infty 
.
\end{equation}
Define the map $\phi:[0,\infty)\to [0,\infty]$ as 
$$
\phi(r)
\equiv 
\inf_{\|y\|_\B\ge r} \rho(y) 
,\spa 
\forall\ r\ge 0 
.
$$
Then $\phi$ is non-negative, non-decreasing, and satisfies 
$$
\phi(\|y\|_\B) 
\le 
\rho(y)
,\spa 
\forall\ y\in \B 
.
$$
For all $t>0$, the continuity of measures and Lévy's maximal inequality---see, e.g., Theorem~3.1.11 and Exercise 2.6.8 in \citet*{GinNic2016})---entail that
\begin{eqnarray*}
\lefteqn{
\P\big[ \sup_{k} \|\dot F[G_k]\|_\B > t \big]
=
\lim_{n\to\infty} \P\big[ \max_{k\le n} \|\dot F[G_k]\|_\B > t\big]
\le 
2 \limsup_{n\to\infty} \P\big[ \|\dot F[G_n]\|_\B > t \big]
}
\\[2mm]
&&
\le
2 \limsup_{n\to\infty} \P\big[ \phi\big(\|\dot F[G_n]\|_\B\big) > \phi(t) \big]
\le 
2 \limsup_{n\to\infty} \P\big[ \rho\big(\dot F[G_n]\big) > \phi(t) \big]
.
\end{eqnarray*}
We deduce from Markov's inequality 
$$
\P\big[ \sup_{k} \|\dot F[G_k]\|_\B = \infty \big]
\le 
\inf_{t>0} \frac{2}{\varphi(t)} \times \limsup_{n\to\infty} \E\big[ \rho\big(\dot F[G_n]\big) \big]
=
0
,
$$
where we used (\ref{asmpFiniteE}) and the fact that $\varphi(t)\to \infty$ as $t\to\infty$ by assumption. We thus have
$$
Y\equiv 
\sup_{k}  
\sup_{\substack{\ell\in\B^*\\ \|\ell\|_{\B^*}\le 1}}
|\ell(\dot F[G_k])| 
= 
\sup_{k} \|\dot F[G_k]\|_\B 
<
\infty
,
$$
almost surely. Let $Y_{k,\ell}=\ell(\dot F[G_k])$ and $\sigma_{k,\ell}^2\equiv \E[Y_{k,\ell}^2]$. Let $M\in (0,\infty)$ 
and $z_{1/2}\in (0,\infty)$ be such that $\P[|\NN(0,1)|>z_{1/2}] = 1/2$. Then 
$$
\P\big[ \sigma_{k,\ell} |\NN(0,1)| > M \big]
=
\P\big[ |Y_{k,\ell}| > M \big]
\le 
\P\big[ Y > M \big]
<
\frac12 
=
\P\big[ |\NN(0,1)| > z_{1/2} \big]
.
$$
We deduce that $\sigma_{k,\ell}\le M/z_{1/2}$, uniformly in $k\ge 1$ and $\ell\in\B^*$ with $\|\ell\|_{\B^*}\le 1$.
Let us now compute
$$
\sigma_{k,\ell}^2
=
\E\big[\ell(\dot F[G_k])^2\big]
= 
\sum_{j=1}^k \ell(\dot F[h_k])^2 
=
\sum_{j=1}^k \big\langle\dot F^*[\ell] , h_k\rangle_\VV^2 
.
$$ 
Since $\sup_{k} \sigma_{k,\ell}<\infty$, Proposition \ref{PropRangeBasis} entails that $\dot F^*[\ell]\in\SS$ and $\sigma_{k,\ell}\to \|\dot F^*[\ell]\|_\SS$ as $k\to\infty$. It follows from the uniform bound $\sigma_{k,\ell}\le M/z_{1/2}$ that 
\begin{eqnarray*}
\lefteqn{
\hspace{-5mm}
\sup_{\substack{h\in \VV_0\\ \|h\|_\HH\le 1}} 
\|\dot F[h]\|_\B
=
\sup_{\substack{h\in \VV_0\\ \|h\|_\HH\le 1}} 
\sup_{\substack{\ell\in\B^*\\ \|\ell\|_{\B^*}\le 1}} 
\ell(\dot F[h]) 
=
\sup_{\substack{h\in \VV_0\\ \|h\|_\HH\le 1}} 
\sup_{\substack{\ell\in\B^*\\ \|\ell\|_{\B^*}\le 1}} 
\langle\dot F^*[\ell] , h\rangle_\VV 
}
\\[2mm]
&&\hspace{20mm}=
\sup_{\substack{\ell\in\B^*\\ \|\ell\|_{\B^*}\le 1}} 
\sup_{\substack{h\in \VV_0\\ \|h\|_\HH\le 1}} \langle\dot F^*[\ell] , h\rangle_\VV 
=
\sup_{\substack{\ell\in\B^*\\ \|\ell\|_{\B^*}\le 1}} \|\dot F^*[\ell]\|_\SS 
\le 
\frac{M}{z_{1/2}}
.
\end{eqnarray*}
Hence $\dot F$ is continuous from $(\VV_0, \ps{\cdot}{\cdot}_\HH)$ to $(\B, \|\cdot\|_\B)$, a contradiction. We deduce that 
$$
\limsup_{n\to\infty} \E\big[ \rho(\dot F[G_n]) \big]
=
\infty 
.
$$
Taking the limit in (\ref{eq:Lrbound2}) yields 
$
\RRR_{\theta_0}(\rho, F, \VV_0) 
=
\infty
,
$
which concludes this part of the proof.\medskip

(iii) Assume that there exists no centred $\B$-valued random variable with covariance as in~(\ref{eq:covF}). We may further assume that $\dot F : (\VV_0, \ps{\cdot}{\cdot}_\HH)\to (\B, \|\cdot\|_\B)$ is continuous, since the conclusion follows immediately from Part (ii) of the proof otherwise. The equalities preceeding the conclusion of Part (ii) entail that $\dot F^* : \B^* \to (\overline \VV_0, \ps{\cdot}{\cdot}_\VV)$ in fact defines a continuous operator $\dot F^* : \B^*\to (\SS, \|\cdot\|_\SS)$. Recall from (\ref{eq:Lrbound2}) that 
\begin{equation}\label{eq:Lrbound3}
\RRR_{\theta_0}(\rho, F, \VV_0) 
\ge
\E\big[ \rho(\dot F[G_k])\big]
,\spa 
\forall\ k\ge 1
.
\end{equation}
Since $\B$ is separable and reflexive, then $\B^*$ is separable and so is its unit ball. Let $D\subset \B^*$ be an at most countable $\|\cdot\|_{\B^*}$-dense subset in the unit ball of $\B^*$. For all $\ell\in D$, the sequence $(\ell(\dot F[G_k]))_k$ is Cauchy in $L^2(\P)$ since 
$$
\E\big[\ell(\dot F[G_k])^2\big]
\to 
\sum_{k} \ell(\dot F[h_k])^2 
=
\sum_{k} \big\langle\dot F^*[\ell] , h_k\rangle_\VV^2 
=
\|\dot F^*[\ell]\|^2_\SS 
<
\infty
,
$$
by Parseval's formula in the Hilbert space $(\HH, \ps{\cdot}{\cdot}_\HH)$ and since $\dot F^* : \B^* \to \HH$. Thus, the sequence of partial sums defining $\ell(\dot F[G_k])$ converges in $L^2(\P)$ to a centred Gaussian random variable $W(\ell)$ with 
$$
W(\ell) 
=
\sum_{k} \gamma_k\ \ell(\dot F[G_k])
.
$$
Lévy's equivalence theorem entails that this convergence also holds almost surely. Since $D$ is at most countable, the previous convergence holds almost surely for all $\ell\in D$ simultaneously, to the effect that the variables $\{W(\ell) : \ell\in D\}$ are well-defined on a same event of probability one. From the series representation in the last display, we see that $W$ can be extended to ${\rm span}(D)$ by linearity, hence defines a linear map $W:{\rm span}(D)\to \R$ with probability one, and the corresponding process is centred and has covariance as in (\ref{eq:covF}). Assume, \emph{ad absurdum}, that 
\begin{equation}\label{asmpFiniteE2}
\limsup_{n\to\infty} \E\big[ \rho(\dot F[G_n]) \big]
<
\infty 
.
\end{equation}
We will show that this implies that $W$ extends to an element of $\B^{**}\cong \B$, which will provide the contradiction. Arguing as in Part (ii) of the proof, we deduce that
$$
\P\big[ \sup_{k} \|\dot F[G_k]\|_\B < \infty \big]
=
1
,
$$
Now observe that the inequalities
$$
\sup_{k} \|\dot F[G_k]\|_\B 
=
\sup_{\ell\in D} \sup_{k} |\ell(\dot F[G_k])| 
\ge 
\sup_{\ell\in D} \limsup_{k\to \infty} |\ell(\dot F[G_k])| 
\ge 
\sup_{\ell \in D} |W(\ell)| 
,
$$
hold almost surely. We deduce that 
$$
\sup_{\ell\in D} |W(\ell)| 
<
\infty 
$$
almost surely. Consequently, with probability one, the linear map $W:{\rm span}(D)\to \R$ extends to a continuous and linear map on $\B^*$ since ${\rm span}(D)$ is dense in $\B^*$, hence $\P[W\in \B^{**}]=1$. Because $\B^{**}$ is reflexive, with probability one there exists an element $Z\in \B$ such that $W(\ell)=\ell(Z)$ for all $\ell\in\B^*$. As a consequence of the Hahn-Banach theorem, the constructed $\B$-valued map $Z$ defined on the underlying probability space is Borel measurable since $\B$ is separable and $\ell(Z)$ is Borel measurable for all $\ell\in\B^*$. Then, $Z$ is a centred $\B$-valued Gaussian random variable with covariance as in (\ref{eq:covF}), a contradiction. We deduce that 
$$
\limsup_{n\to\infty} \E\big[ \rho(\dot F[G_n]) \big]
=
\infty 
.
$$
Consequently, taking the limit in (\ref{eq:Lrbound2}) yields 
$
\RRR_{\theta_0}(\rho, F, \VV_0) 
=
\infty
,
$
which concludes the proof.
\end{Proof}
\section{Application to data assimilation for non-linear parabolic PDEs}
\label{sec:Appli}

In this section, we derive the main information geometric properties of generic non-linear and non-Gaussian models, for which we will check the assumptions in Section \ref{sec:AppliReac} and Section \ref{sec:AppliNS} below in a data assimilation context when the regression function describes the evolution over time of a non-linear parabolic dynamical system, with an emphasis on reaction-diffusion systems and Navier-Stokes systems.
To avoid some technicalities and streamline the exposition, we consider the space-time cylinder $\XX=[0,T]\times \Om$, where $\Om$ is some $d$-dimensional domain. For simplicity we consider $\Om=\T^d$, the $d$-dimensional torus, to avoid boundary issues which permits a simple interpretation of duals of Sobolev spaces as `negative-order' Sobolev spaces; see Example \ref{ex:HandSBis} and the comments thereafter.\smallskip

Let $(\VV,\ps{\cdot}{\cdot}_\VV)$ be a separable Hilbert space and $\theta_0\in\Theta$, where $\Theta\subset \VV$ is a Borel measurable subset. Let $\GGG:\Theta\to L^2_\lambda \equiv L^2([0,T]\times \Om, \R^p, d\lambda)$ be a measurable map,  where $\lambda$ is a probability measure on $[0,T]\times \Om$ with density, abusively denoted by $\lambda(t,x)$, such that 
\begin{equation}\label{eq:BoundsLambda}
\lambda_{\min}
\le 
\lambda(t,x)
\le 
\lambda_{\max}
,\spa 
\forall\ t\in [0,T]\times \Om
,
\end{equation}
for some $0<\lambda_{\min}\le \lambda_{\max}<\infty$. This assumption means that the design in the regression model below will eventually explore all parts of $[0,T]\times \Om$, and further allows to relate the $L^2_\lambda$-norm to its usual $L^2$-antecedent: the two-sided estimates in Proposition \ref{CoroMinimaxHilbert} below can thus be simply expressed in terms of the $L^2([0,T]\times \Om)$-norm. Consider data $(t_i,\omega_i, Y_i)_{i=1}^N$ arising from 
$$
Y_i 
=
\GGG(\theta)(t_i, \omega_i)
+
\ve_i 
,\spa 
i=1,\ldots, N
$$
where $X_i\equiv (t_i, \omega_i)$ are i.i.d. samples drawn from $\lambda$, and $(\ve_i)$ are i.i.d. samples drawn from some density $q_\ve\in L^1(\R^p)$ such that $\sqrt{q_\ve}\in H^1(\R^p)$ independently from $(X_i)$, and let 
$$
\II_\ve 
=
4 \int_{\R^p} (\nabla \sqrt{q_\ve})(y) (\nabla \sqrt{q_\ve})(y)^T\, dy 
.
$$
We can now combine Example \ref{ex:HandSBis}, Example \ref{ex:SeriesG}, Example \ref{ex:SuppHilbert}, and Example \ref{ex:ContBanach} to obtain the following general result on the information geometry of such models, for which we will check the assumptions in Section \ref{sec:AppliReac} and Section \ref{sec:AppliNS} below in the case of reaction-diffusion equations and Navier-Stokes equations.


\begin{Prop}\label{CoroMinimaxHilbert}
Consider the setting of this Section. For some linear subspace $\VV_0\subset \VV$, assume that $\GGG$ admits a continuous linearization $\I_{\theta_0}:(\VV_0, \ps{\cdot}{\cdot}_\VV)\to L^2_\lambda$ satisfying Condition \ref{CondIntroLinModel},
and let $\I_{\theta_0}^* : L^2_\lambda \to (\overline \VV_0, \ps{\cdot}{\cdot}_\VV)$ be its adjoint, where $\overline\VV_0$ denotes the $\|\cdot\|_\VV$-closure of $\VV_0$. Assume that 
$$
\|h\|^2_{\DD^{-\kappa}}
\lesssim 
\int_0^T \int_{\Om} |\I_{\theta_0}[h](t,x)|^2\, dx\, dt 
\lesssim 
\|h\|^2_{\DD^{-\kappa}}
,\spa 
\forall\ h\in\VV_0 
,
$$
where, for some positive sequence $\tau_j \asymp j^{1/\alpha}$, $\alpha>0$, and an orthonormal basis $(e_j)\subset \SS$ of $\overline \VV_0$,
$$
\DD^s
\equiv 
\Big\{ \sum_{j\ge 1} u_j e_j : \|u\|^2_{\DD^s}\equiv \sum_{j\ge 1} \tau_j^s u_j^2 < \infty \Big\} 
,\spa 
s\in\R 
.
$$
Then, the following holds.
\begin{itemize}
    \item[(i)] The model $(P_{\theta}^N : \theta\in\Theta)_{N\ge 1}$ is LAN around $\theta_0$ with respect to the tangent space $\VV_0$, and the Fisher information operator of the model is given by $h\mapsto (\I_{\theta_0}^* \II_\ve \I_{\theta_0})[h]$, for $h\in\VV_0$. 
    
    \item[(ii)] The matrix $\II_\ve$ is positive definite, and the information operator $\I_{\theta_0}^* \II_\ve \I_{\theta_0}$ extends to a homeomorphism between $\DD^{-\kappa}$ and $\DD^\kappa$, regardless of $q_\ve$ and $\II_\ve$.
    
    \item[(iii)] If $\psi\in \DD^\kappa$, then the local asymptotic minimax risk for estimating $F_\psi\equiv \ps{\psi}{\theta_0}_\VV$ with tangent space $\VV_0$ in Euclidean $|\cdot|^2$-loss is lower bounded as
    $$
    \RRR_{\theta_0}(|\cdot|^2, F_\psi, \VV_0)
    \ge 
    \ps{\psi}{(\I_{\theta_0}^* \II_\ve \I_{\theta_0})^{-1}[\psi]}_\VV
    ,
    $$
    and the risk above is infinite if $\psi\in \overline \VV_0\sm \DD^\kappa$. 
    
    \item[(iv)] For centred normal random variables $(\gamma_j)$ with covariance 
    $\E[\gamma_i \gamma_j]=\ps{e_i}{(\I_{\theta_0}^*\II_\ve \I_{\theta_0})^{-1}[e_j]}_\VV$, the series $G=\sum_j \gamma_j e_j$ converges almost surely in $\DD^{-\beta}$ for any $\beta>\kappa+\alpha$, and has covariance 
    $$
    \E\big[\ps{G}{f}_\VV \ps{G}{g}_\VV\big]
    =
    \ps{f}{(\I_{\theta_0}^* \II_\ve \I_{\theta_0})^{-1}[g]}_\VV 
    ,\spa 
    \forall\ f,g\in \DD^\beta
    ,\ \beta > \kappa+\alpha
    .
    $$

    \item[(v)] Let $(\B,\|\cdot\|_\B)$ be a Banach space and $F : \Theta \to \B$ a measurable map satisfying Condition~\ref{CondRegularParameters} with linearization $\dot F: \VV_0\to \B$. Further assume that, for some $\beta>\kappa+\alpha$, we have $\|\dot F[h]\|_\B\lesssim \|h\|_{\DD^{-\beta}}$ for all $h\in\VV_0$. Then $\dot F[G]$ defines a centred $\B$-valued Gaussian random variable with covariance as in (\ref{eq:covF}), and for all $s>0$, we have 
    $$
	\RRR_{\theta_0}(\|\cdot\|^s_{\B}, F, \VV_0)
	\geq 
	\E\Big[\|\dot F[G]\|^s_{\B}\Big]
    .
	$$
\end{itemize}
\end{Prop}

Let us stress that the adjoint $\I_{\theta_0}^*$ in Proposition \ref{CoroMinimaxHilbert} is computed with respect to the $L^2_\lambda$-norm, not the $L^2([0,T]\times\Om, \R^p, dx)$-one. However, in view of (\ref{eq:BoundsLambda}), it is enough to establish two-sided estimates on $\|\I_{\theta_0}[h]\|_{L^2([0,T]\times \Om)}$.

\subsection{Reaction-diffusion equations}\label{sec:AppliReac}

For $\Om=\T^d$, with $d\in\{1,2,3\}$, we consider $\GGG(\theta)$ to be the solution $u_\theta :[0,T]\times \Om\to\R$ of the reaction-diffusion equation
\begin{align}\label{eq:RD}
\frac{\partial u_\theta}{\partial t}(t,x)
-
\Delta u_\theta(t,x) &= f(u_\theta(t,x))
~~~ \textrm{on } (0,T]\times\Om, \\ \nonumber
u_\theta(0,x)&=\theta(x)\spa 
~~~~ \textrm {on } \Om
,
\end{align}
where $\theta\in H^1(\Om)$ is the initial condition and the reaction term $f$ is taken in $C^\infty_c(\R)$. A classical theory for the existence, uniqueness and regularity of such solutions can be found, e.g., in \citet*{Robinson2001}. The main properties of this system of equations in a statistical context are also accounted for in \citet*{Nickl2024}. Theorem 6 there entails that, for any $\theta\in H^1(\Om)$, there exists a unique strong solution $u_\theta\in L^2([0,T]\times \Om)$ to the equation above; see Section 3.1 in the same reference for details on the notions of strong and weak solutions in this setting. \smallskip 

We consider the separable Hilbert space $(\VV, \ps{\cdot}{\cdot}_\VV)=(L^2(\Om), \ps{\cdot}{\cdot}_{L^2(\Om)})$, the Borel measurable subset $\Theta=H^1(\Om)\subset L^2(\Om)$, and $L^2_\lambda=L^2([0,T]\times \Om, \R, d\lambda)$ (\ie $p=1$). In view of~(\ref{eq:BoundsLambda}), $\GGG:\Theta\mapsto L^2_\lambda$ is thus a well-defined map. For $\theta_0\in C^\infty(\Om)$ and tangent space $\VV_0 = H^2(\Om)$, Theorem 7 in \citet*{Nickl2024} with $\bar\gamma\equiv 2 > \max\{1, d/2\}$ entails that $\GGG:\Theta\to L^2_\lambda$ satisfies Condition~\ref{CondIntroLinModel} with linearization $\I_{\theta_0} : H^2(\Om)\to L^2_\lambda$, where $\I_{\theta_0}[h]$ is given, for all $h\in H^2(\Om)$, by the solution $U=U_h$ to the linearized PDE 
\begin{align*}
\frac{\partial U}{\partial t}(t,x)
-
\Delta U(t,x) - f'(u_{\theta_0}(t,x)) U(t,x) &= 0\spa\spa\textrm{on } (0,T]\times\Om, \\ 
U(0,x)&=h(x)\spa 
~~ \textrm {on } \Om
.
\end{align*}
The linear operator $\I_{\theta_0}$ is continuous from $(H^2(\Om), \ps{\cdot}{\cdot}_{L^2(\Om)})$ to $L^2_\lambda$ by virtue of Proposition 5 in \citet*{Nickl2024} with $a=-1$ and potential $V=f'(u_{\theta_0})$. In fact, this result and Lemma 4 there yield the two-sided estimate 
\begin{equation}\label{eq:RDEstimate}
\|h\|_{H^1(\Om)^*}
\lesssim 
\|\I_{\theta_0}[h]\|_{L^2([0,T]\times \Om)} 
\lesssim 
\|h\|_{H^1(\Om)^*} 
,\spa 
\forall\ h\in H^2(\Om) 
,
\end{equation}
where 
$$
\|h\|_{H^1(\Om)^*}
=
\sup_{\substack{f\in H^1(\Om)\\ \|f\|_{H^1(\Om)}\le 1}} \Big| \int_{\Om} h(x)f(x)\, dx \Big|
.
$$
Each Sobolev space $H^k(\Om)$, $k\in\N$, has equivalent norm
$$
\|u\|^2_{H^k(\Om)}
\asymp 
\sum_{j\ge 0} (1+\lambda_j)^k \ps{u}{e_j}^2_{L^2(\Om)}
,
$$
where $(\lambda_j)$ are the eigenvalues of the periodic Laplacian $-\Delta$ over $\Om=\T^d$, ordered so that $0=\lambda_0\le \lambda_{j-1}\le \lambda_{j}\simeq j^{2/d}$, with corresponding eigenfunctions $(e_j)$ forming an orthonormal basis of $L^2(\Om)$; see also Section 1.3 in \citet*{Nickl2024}. The definition of the previous series norm can be extended to any $k\in\R$, and one sees directly that the $\|\cdot\|_{H^1(\Om)^*}$-norm is equivalent to the series norm above with $k=-1$. Consequently, the assumptions of Proposition \ref{CoroMinimaxHilbert} above hold with $\kappa=1$,  $\tau_j=1+\lambda_j$ and $\alpha=d/2$.

\begin{Prop}\label{PropRD}
    For a fixed $f\in C^\infty_c(\Om)$, let $\GGG(\theta):[0,T]\times \Om\to \R$ be the unique strong solution to the reaction-diffusion equation (\ref{eq:RD}) with initial condition $\theta\in H^1(\Om)$. If $\theta_0\in C^\infty(\Om)$, then Proposition~\ref{CoroMinimaxHilbert} holds for parameter space $\Theta=H^1(\Om)$, tangent space $\VV_0 = H^2(\Om)$ and corresponding $\overline \VV_0 = L^2(\Om)$, with $\kappa=1$, $\alpha=d/2$, $\DD^\kappa = H^1(\Om)$ and $\DD^{-\kappa}=H^1(\Om)^*$.
\end{Prop}

\begin{Remark}
    One can obtain similar results for suitable restrictions of $\VV_0$. For instance, take the subspace $\VV_0=\dot H^2(\Om)$ consisting of those $h\in H^2(\Om)$ such that $\int_\Om h = 0$. Its $\|\cdot\|_{L^2(\Om)}$-closure $\overline \VV_0$ is a strict subspace of $L^2(\Om)$, and the resulting Hilbert space $\HH$ (see Section~\ref{sec:Extend} and Example~\ref{ex:Stability}), providing $\DD^{-\kappa}$ in Proposition \ref{PropRD} above, is a strict subspace of $H^1(\Om)^*$ with dual $\DD^\kappa$ equal to $\dot H^1(\Om)$ consisting of those $h\in H^1(\Om)$ such that $\int_\Om h=0$.
\end{Remark}


To establish the asymptotic efficiency lower bound in Proposition~\ref{CoroMinimaxHilbert}(v) to functionals $F:\Theta\to (\B, \|\cdot\|_\B)$ of the parameter $\theta$ with linearization $\dot F : \VV_0 \to \B$, the main assumption to check is the continuity condition $\|\dot F[h]\|_\B\lesssim \|h\|_{\DD^{-\beta}}$ for some $\beta>\kappa+\alpha$. Because $\alpha>0$, this requirement is stronger than the continuity of $\dot F$ between $(\VV_0, \ps{\cdot}{\cdot}_\HH)$ and $(\B, \|\cdot\|_\B)$ from Theorem \ref{TheorGeneralMinimax}(i). One of the most natural and interesting functionals to consider is arguably the regression function itself, \ie $F(\theta)=u_\theta$. Under some appropriate differentiability, the continuity requirement above thus rewrites $\|\I_{\theta_0}[h]\|_\B\lesssim \|h\|_{H^\beta(\Om)^*}$ for some $\beta>1+d/2$. This cannot hold, for instance, for $\|\cdot\|_\B$ given by the $L^2(\Om\times [0,T])$-norm, in view of (\ref{eq:RDEstimate}). In fact this can only hold for $\|\cdot\|_\B$-norms corresponding to quite weak topologies: in view of Proposition 5 and an adaptation of Lemma 4 in \citet*{Nickl2024}, one expects that, for a fixed $\beta>1+d/2$, $\|\cdot\|_\B$ should be taken as the $L^2([0,T], H^{\beta-1}(\Om)^*)$-norm, \ie 
$$
\|\I_{\theta_0}[h]\|^2_\B 
\equiv 
\int_0^T \|\I_{\theta_0}[h](t,\cdot)\|^2_{H^{\beta-1}(\Om)^*}\, dt
.
$$
As already discussed in \citet*{NicTit2024}, \citet*{Nickl2024}, and \citet*{KonNic2025}, these limitations on the choice of $\B$ are unavoidable when $\|\I_{\theta_0}[h]\|_\B$ aggregates information on $\I_{\theta_0}[h]$ at times including the `initial' time $t=0$. This is due to the smoothing nature of the linearized flow above, of which $\I_{\theta_0}[h]$ is a solution, and the exact blow-up rate of the `semi-group' $S_t[h]=\I_{\theta_0}[h](t,\cdot)$ as $t\downarrow 0$. This smoothing phenomenon, however, allows to consider much stronger norms $\B$ than previously provided one restricts to (strictly) positive times $t>0$. This is the content of the next proposition.

\begin{Prop}\label{PropFunctRD}
    Consider the setting of Proposition \ref{PropRD}.
    Fix $0<t_0<t_1<\infty$ and define the Banach space $\B\equiv (C([t_0, t_1]\times\Om), \|\cdot\|_{L^\infty})$.
    Then, Proposition~\ref{CoroMinimaxHilbert}(v) applies to the functional 
    $$
    F:H^2(\Om) \to C([t_0, t_1]\times \Om),\ \theta\mapsto u_\theta
    ,
    $$
    for any $\beta\ge 3$ with corresponding $\DD^{-\beta}=H^\beta(\Om)^*$, and continuous linearization at $\theta_0$ given by 
    $$
    \dot F : (H^2(\Om), \ps{\cdot}{\cdot}_{H^\beta(\Om)^*}) \to (C([t_0, t_1]\times \Om),\|\cdot\|_{L^\infty}),\ h\mapsto \I_{\theta_0}[h]
    .
    $$
\end{Prop}

\begin{Proof}{Proposition \ref{PropFunctRD}}
	Observe that $F$ is well-defined since $u_\theta\in C([t_0, t_1],H^2(\Om))\embed C([t_0, t_1]\times \Om)$ by virtue of Proposition 4 in \citet*{Nickl2024} with $\gamma\equiv 2$ and $\theta\equiv \theta_0\in C^\infty(\Om)$, and the Sobolev embedding $H^2(\Om)\embed C(\Om)$ since $d\le 3$. An adaptation of Theorem 7 there yields 
    $$
    \sup_{t\in [t_0, t_1]} \| u_{\theta_0+sh}(t)-u_{\theta_0}(t) - \I_{\theta_0}[sh](t)\|_{L^\infty(\Om)}
    =
    o(s)
    ,\spa 
    s\to 0 
    .
    $$
    Indeed, letting $w\equiv u_{\theta_0+sh}-u_{\theta_0}$ and using that 
    \begin{align*}
    m(t)
    &\equiv 
    f(u_{\theta_0+sh}(t))-f(u_{\theta_0}(t))-f'\big(u_{\theta_0}(t)\big) \big(u_{\theta_0+sh}(t)-u_{\theta_0}(t)\big)
    \\[2mm]
    &=
    \frac12 \big(u_{\theta_0+sh}(t)-u_{\theta_0}(t)\big)^2 \int_0^1 (1-\tau^2) f''\big(u_{\theta_0}(t)+\tau(u_{\theta_0+sh}(t)-u_{\theta_0}(t))\big)\, d\tau
    \\[2mm]
    &=
    \frac12 w(t)^2 \int_0^1 (1-\tau^2) f''\big(u_{\theta_0}(t)+\tau w(t)\big)\, d\tau
    ,
    \end{align*}
    by Taylor's theorem,
    one replaces (105) in \citet*{Nickl2024} with the estimate 
    $$
    \sup_{t\in [t_0, t_1]} \|r(t)\|_{L^\infty(\Om)}
    \lesssim 
    \sup_{t\in [t_0, t_1]} \|r(t)\|_{H^2(\Om)}
    \lesssim 
    \sup_{t\in [t_0, t_1]}\|m(t)\|_{H^2(\Om)}
    ,
    $$
    where we used the Sobolev embedding $H^2(\Om)\embed L^\infty(\Om)$ and the parabolic regularity estimate from Proposition 5 with $a=2$. Now, we can upper bound $\|m(t)\|_{H^2(\Om)}$ uniformly in $t\in [t_0, t_1]$ as
    \begin{eqnarray*}
    \lefteqn{
    \hspace{-20mm}
    \|m(t)\|_{H^2(\Om)}
    \lesssim
    \|f\|_{C^4(\Om)} \|w(t)^2\|_{H^2(\Om)}
    \lesssim 
    \|w(t)^2\|_{L^2(\Om)} + \|\Delta(w(t)^2)\|_{L^2(\Om)} 
    }
    \\[2mm]
    &&
    \lesssim 
    \|w(t)\|^2_{L^4(\Om)}
    +
    \|\nabla w(t) \|^2_{L^4(\Om)}
    +
    \|w(t)\Delta w(t) \|_{L^2(\Om)}
    \\[2mm]
    &&
    \lesssim 
    \|w(t)\|^2_{H^1(\Om)}
    +
    \|\nabla w(t)\|^2_{H^1(\Om)}
    +
    \|w(t)\|_{L^\infty(\Om)} \| \Delta w(t)\|_{L^2(\Om)}
    \\[2mm]
    &&
    \lesssim 
    \|w(t)\|^2_{H^2(\Om)}
    \lesssim 
    s^2\|h\|^2_{H^2(\Om)}
    ,
    \end{eqnarray*}
    where we used the multiplier inequality for Sobolev norms (see (17) in \citet*{Nickl2024} with $a\equiv 2>d/2$), the pointwise identity 
    $
    \Delta(g^2) 
    =
    2(|\nabla g|^2 + g\Delta g)
    ,
    $
    the Sobolev embedding $H^1(\Om)\embed L^4(\Om)$ (we have $H^1(\Om)\embed L^6(\Om)$ for $d\in\{1,2,3\}$), the Sobolev emebedding $H^2(\Om)\embed L^\infty(\Om)$, and the inequality (100) in \citet*{Nickl2024} with $\zeta\equiv 2>d/2$.
	Consequently, $F:H^2(\Om)\to C([t_0, t_1]\times \Om)$ admits the linearization $\dot F : H^2(\Om)\to C([t_0, t_1]\times\Om), h\mapsto \I_{\theta_0}[h]$ in the sense of Condition \ref{CondRegularParameters}. Finally observe that, from the smoothing nature of the flow, Lemma 2 in \citet*{Nickl2024} with $b\equiv 2$ and $a\equiv-3$ provides
	$$
	\|\dot F[h]\|_{\B}
	=
	\sup_{t\in [t_0, t_1]} \|\I_{\theta_0}[h]\|_{L^\infty(\Om)}
	\lesssim 
	\sup_{t\in [t_0, t_1]} \|\I_{\theta_0}[h](t)\|_{H^2(\Om)}
	\lesssim 
	\|h\|_{H^3(\Om)^*}
	.
	$$
    We deduce that Proposition \ref{CoroMinimaxHilbert}(v) holds with $\VV_0 = H^2(\Om)$ and $\beta\ge 3> \kappa+\alpha=1+d/2$.
\end{Proof}


\subsection{Navier-Stokes equations}\label{sec:AppliNS}

For $\Om=\T^2$ we consider $\GGG(\theta)$ to be the solution $u_\theta \equiv u :[0,T]\times \Om\to\R^2$ of the (periodic) incompressible 2D Navier-Stokes equations
\begin{align}\label{eq:NS}\nonumber
    \frac{\partial u}{\partial t} - \nu \Delta u + (u\cdot \nabla) u  &= f - \nabla p
    ~~~ \textrm{on } (0,T]\times\Om, \\
    \nabla\cdot u &= 0\spa
    ~~~~~ \textrm {on } [0,T]\times \Om\\[1mm]
    u(0,x)&=\theta(x)\spa
    ~ \textrm {on } \Om, \nonumber
 \end{align}
 where $\theta$ is an initial condition satisfying $\int_\Om \theta = 0$ and $\nabla\cdot \theta=0$, $\nu>0$ is a known viscosity, $f:\Om\to\R^2$ is a given exterior forcing independent of time, and $p:[0,T]\times \Om\to \R$ is a scalar pressure term. A classical theory for the existence, uniqueness and regularity of such solutions can be found, e.g., in \citet*{Robinson2001}. The main properties of this system of equations in a statistical context are also accounted for in \citet*{NicTit2024} and \citet*{KonNic2025}. Using the notation of the latter, the relevant Sobolev spaces for this equation are 
 $$
 \divH^k 
 \equiv 
 H^k(\Om)^2 \cap \Big\{ h : \Om\to\R^2,\ \int_\Om h = 0,\ \nabla\cdot h = 0\Big\} 
 ,\spa 
 k\ge 0 
 ,
 $$
 with the usual convention $H^0(\Om)^2=L^2(\Om)^2$.
 These spaces admit the series norms 
 $$
 \|u\|^2_{\divH^k}
 \equiv 
 \sum_{j\ge 1} \lambda_j^k\ \lvert\ps{u}{e_j}_{L^2(\Om)}\rvert^2
 ,
 $$
 where $0\le \lambda_j\simeq j$ are the eigenvalues of the periodic Laplacian on $\T^2$ and $(e_j)$ is a suitable (complex-valued) orthonormal basis of $\divH\equiv \divH^0$; see Section A.1 in \citet*{KonNic2025}. Proposition A.5 there entails that, for any $\theta\in \divH^1$ and $f\in\divH$, there exists a unique strong solution $u_\theta\in L^2([0,T], \divH)$ to (\ref{eq:NS}) with initial condition $\theta$. Throughout, we fix $f\in \divH^1$. \smallskip 

We consider the separable Hilbert space $(\VV, \ps{\cdot}{\cdot}_\VV)=(L^2(\Om)^2, \ps{\cdot}{\cdot}_{L^2(\Om)})$, the Borel measurable subset $\Theta=\divH^2\subset L^2(\Om)$, and $L^2_\lambda=L^2([0,T]\times \Om, \R^2, d\lambda)$ (\ie $p=2$). In view of~(\ref{eq:BoundsLambda}), $\GGG:\Theta\mapsto L^2_\lambda$ is thus a well-defined map. For $\theta_0\in \divH^2$ and tangent space $\VV_0 = \divH^2$, Proposition~A.8 and Proposition A.9 with $a=-1$ in \citet*{KonNic2025} entail that $\GGG:\Theta\to L^2_\lambda$ is Gâteaux-differentiable at $\theta_0$ in all directions $h\in \VV_0$ with continuous linearization $\I_{\theta_0} : (\divH^2, \ps{\cdot}{\cdot}_{L^2(\Om)})\to L^2_\lambda$.
Proposition A.9 and Proposition A.10 there further yield 
the two-sided estimate 
$$
\|h\|_{(\divH^1)^*}
\lesssim 
\|\I_{\theta_0}[h]\|_{L^2([0,T]\times \Om)} 
\lesssim 
\|h\|_{(\divH^1)^*} 
,\spa 
\forall\ h\in \divH^2 
.
$$
where 
$$
\|h\|_{(\divH^1)^*}
=
\sup_{\substack{g\in \divH^1\\ \|g\|_{\divH^1}\le 1}} \Big| \int_{\Om} h(x)g(x)\, dx \Big|
.
$$
Since the $\{\lambda_j:j\ge 1\}$ satisfy a spectral gap $\lambda_1>0$ (see Section A.1.1. in \citet*{KonNic2025}), the definition of the previous series norm can be extended to any $k\in\R$, and one sees directly that the $\|\cdot\|_{\divH^1}$-norm is equivalent to the series norm above with $k=-1$. Consequently, the assumptions of Proposition \ref{CoroMinimaxHilbert} above hold with $\kappa=1$,  $\tau_j=\lambda_j$ and $\alpha=1$.

\begin{Prop}\label{PropNS}
    Fix $\nu>0$ and $f\in \divH^1$, and let $\GGG(\theta):[0,T]\times \Om\to \R^2$ be the solution to the Navier-Stokes equations (\ref{eq:NS}) with initial condition $\theta\in \divH^2$. If $\theta_0\in \divH^2$, then Proposition~\ref{CoroMinimaxHilbert} holds for tangent space $\VV_0 = \divH^2$ and corresponding $\overline \VV_0 = \divH$, with $\kappa=1$, $\alpha=1$, $\DD^\kappa = \divH^1$ and $\DD^{-\kappa}=(\divH^1)^*$.
\end{Prop}

Let us show how we can apply Proposition~\ref{CoroMinimaxHilbert}(v) to particular functionals $F$. Let $(\B,\|\cdot\|_\B)$ a Banach space and $F:\Theta\mapsto \B$ a map of the form 
$$
F(\theta) 
=
\Phi(u_\theta) 
,
$$
where, for $\tilde \B$ a Banach space, $\Phi:\tilde \B\to \B$ is a Hadamard-differentiable map, \ie for all $v\in\tilde \B$ there exists a bounded linear operator $\dot \Phi_v : \tilde\B\to \B$ such that 
$$
\frac{\Phi(v+t_n h_n)-\Phi(v)}{t_n}
\to 
\dot \Phi_v[h]
\spa 
\textrm{in}\quad \B 
,
$$
for all $h\in\tilde\B$, and sequences $(t_n)\subset \R$ and $(h_n)\subset \tilde\B$ such that $t_n\to 0$ and $h_n\to h$ in $\tilde\B$.

\begin{Prop}\label{PropFunct}
    Fix $b\ge 2$, $\theta_0\in \divH^b$ and $f\in \divH^{b-1}$. Let $\Phi:\tilde\B\to \B$ be a Hadamard-differentiable map between Banach spaces $\tilde\B$ and $\B$. Fix $0<t_0<t_1<\infty$ and further assume that $C([t_0, t_1], \divH^b)\embed \tilde\B$. Then, Proposition~\ref{CoroMinimaxHilbert}(v) applies to the functional 
    $$
    F:\divH^b \to \B,\ \theta\mapsto \Phi(u_\theta)
    ,
    $$
    with $\Theta=\VV_0=\divH^b$ and $\beta=3$.
\end{Prop}

\begin{Proof}{Proposition \ref{PropFunct}}
	Observe that $F$ is well-defined since $u_\theta\in C([t_0, t_1],\divH^b)\embed \tilde\B$ by virtue of Proposition A.6 in \citet*{KonNic2025} with $a\equiv b$. Proposition A.12 there with $a\equiv b$ entails that the map $\divH^b\mapsto C([t_0, t_1],\divH^b),\ \theta\mapsto u_\theta$ is Hadamard-differentiable with Hadamard derivative at any $\theta\in \divH^b$ given by 
	$$
	\I_\theta: \divH^b \to C([t_0, t_1],\divH^b)
	,
	$$
	which is a bounded operator by Proposition A.9 with $a\equiv b$. In particular, the map $\theta\mapsto u_\theta$ is also Hadamard-differentiable as a map between $\divH^b$ and $\tilde \B$. Consequently, the chain rule for Hadamard derivatives entails that the map $F$ is Hadamard-differentiable at $\theta_0$ with Hadamard derivative
	$$
	\dot F : \divH^b \to\B 
	,\ 
	h\mapsto
	\dot\Phi_{u_{\theta_0}}[\I_{\theta_0}[h]]
	.
	$$
	Finally observe that the embedding $C([t_0, t_1],\divH^b)\embed \tilde\B$ and Proposition A.11 in \citet*{KonNic2025} with $a\equiv-3$ provide
	$$
	\|\dot F[h]\|_{\B}
	=
	\|\dot\Phi_{u_{\theta_0}}[\I_{\theta_0}[h]]\|_\B
	\lesssim 
	\|\I_{\theta_0}[h]\|_{\tilde\B}
	\lesssim 
	\sup_{t\in [t_0, t_1]} \|\I_{\theta_0}[h](t)\|_{\dH^b}
	\lesssim 
	\|h\|_{(\divH^3)^*}
	,
	$$
    where $3>\kappa+\alpha=2$.
\end{Proof}

\begin{Example}
    When $\Phi$ is the identity on $\CCC\equiv C([t_0,t_1]\times \Om)^2$, then Proposition \ref{PropFunct} yields the minimax lower bound for estimating the trajectory $\{u_{\theta_0}(t,x) : t\in [t_0, t_1], x\in\Om\}$ by taking $b=2$, $\tilde\B=\B=\CCC$, and $\Phi={\rm Id}_{\CCC}$, by virtue of the Sobolev embedding $\divH^b\embed C(\Om)^2$ since $b>1$. 
\end{Example}

\begin{Example}
    We can also apply Proposition \ref{PropFunct} to other functionals of the Navier-Stokes system. In particular, consider the functional $\theta\mapsto (u_{\theta}\cdot\nabla)u_{\theta}$ that yields the non-linear effect of the dynamics in (\ref{eq:NS}), arising from
$$
\Phi : \tilde\B\to \B,\ 
v\mapsto (v\cdot\nabla)v
$$
with Banach spaces $\tilde\B = C([t_0, t_1],\divH^a)$ and $\B=C([t_0, t_1],\divH^{a-1})$ for some $a\geq 3$, and taking $b\equiv a$ in Proposition \ref{PropFunct}. The functional $\Phi$ is well-defined since 
$$
\sup_{t\in [t_0, t_1]} \|(v(t)\cdot\nabla)v(t)\|_{\dH^{a-1}}
\lesssim 
\sup_{t\in [t_0, t_1]} \|v(t)\|^2_{\dH^a}
.
$$
by virtue of the multiplier inequality for Sobolev norms (see Section A.2 in \citet*{KonNic2025}) together with the fact that $(v\cdot\nabla)v=(\nabla v)v$ (see (22) there). The estimates obtained in the last part of the proof of Proposition A.6 in \citet*{KonNic2025} entail that one can take $\tilde\B=C([t_0, t_1],\divH^{a+1})$, instead, to accommodate the case $a\in \{0,1,2\}$. Since the maps $v\mapsto v$ and $v\mapsto \nabla v_i$, $i\in\{1,2\}$, are Hadamard differentiable between $\tilde\B$ and $\B$, then Leibniz's rule for Hadamard derivatives yields the Hadamard-differentiability of $\Phi$. Consequently, Proposition \ref{PropFunct} holds for this choice of $\Phi$ and yields the local asymptotic minimax lower bound for estimating the non-linearity $(u_\theta\cdot\nabla)u_\theta$ at positive times.
\end{Example}

\appendix

\section{Supplementary material}
\subsection{Support of functions in $H^1(\R^p)$}\label{sec:AppendixSupport}

This section establishes a few technical lemmas needed to establish the results of Section \ref{sec:LAN} for densities $q_\ve$ of the errors in the model (\ref{model}) satisfying  $\sqrt{q_\ve}\in H^1(\R^p)$ without imposing further conditions on their support. \smallskip


Vitali's covering lemma establishes that for any measurale set $E\in\R^p$ there is a countable collection of \emph{disjoint} balls $(B_i)$, not necessarily contained in $E$, such that $E\sm \cup_{i\ge 1} B_i$ has zero Lebesgue measure. The usual construction of such balls does not allow to determine how often or if any of these balls overlap with the complement of $E$. The following lemma establishes that, when $E$ is open, the disjoint balls $(B_i)$ can be chosen as subsets of $E$. This result is well-known (see, e.g., Theorem 1.26 in \citet*{EvansGariepy2015}) but we rather present here a simple dyadic argument due to Eric Wofsey on the math StackExcange\footnote{\url{https://math.stackexchange.com/questions/1881977/filling-a-unit-cube-with-countable-balls}}. 
Let us recall that a dyadic cube in $\R^p$ is an open cube of the form 
$$
\prod_{i=1}^p \Big( \frac{j}{2^\ell}, \frac{j+1}{2^\ell} \Big) 
,\spa 
j\in \Z,\ \ell\in \N
.
$$
An important property is that when two of these cubes intersect, then necessarily one of them is a subset of the other. Because $\ell\in\N$, the sides of these cubes all have length no larger than $1$.

\begin{Lem}\label{LemVitali}
    Let $E\subset\R^p$ be an open subset. Then, there exists a countable collection of disjoints open balls $(B_i)$ such that $B_i\subset E$ for all $i$, and $E\sm \cup_{i\ge 1} B_i$ has zero Lebesgue measure.
\end{Lem}

\begin{Proof}{Proof of Lemma \ref{LemVitali}}
    Let $D$ be the set of points $(x_1,\ldots, x_p)\in E$ having at least one dyadic coordinate, \ie there is $i\in\{1,\ldots, p\}$ such that $x_i = j2^{-\ell}$ for some $j\in \Z$ and $\ell\in\N$. Since $D$ is a subset of a countable union of hyperplanes, it has zero Lebesgue measure. Since $E$ is an open set, for each $x\in E\sm D$ there is a dyadic cube $C(x)\subset E$ such that $x\in C(x)$. In addition, since all dyadic cubes are uniformly bounded, we can choose $C(x)$ to be the maximal such cube; in particular $C(y)=C(x)$ for all $y\in C(x)$. Since $C(x)\cap D=\emptyset$, we have 
    \begin{equation}\label{eq:Union}
    E\sm D
    =
    \bigcup_{x\in E\sm D} C(x) 
    .
    \end{equation}
    Since there are at most countably many distinct cubes $C(x)$, the union in (\ref{eq:Union}) reduces to a countable index subset $J\subset E\sm D$ for which $C(x)\ne C(x')$ for all $x,x'\in J$ with $x\ne x'$. Next, note that the cubes in the collection $\{C(x) : x\in J\}$ are disjoint. Indeed, for $x,x'\in J$ such that $C(x)\ne C(x')$, if $C(x)$ and $C(x')$ had a non-empty intersection, $y\in C(x)\cap C(x')$, say, then the maximality would entail that $C(x)=C(y)=C(x')$, a contradiction. Consequently, it is enough to show for any dyadic cube $C$ there exists an at most countable collection $(B_i)$ of disjoint balls with $B_i\subset C$ and such that $C\sm \cup_{i\ge 1} B_i$ has zero Lebesgue measure. Let us start with an open ball $B\subset C$, and write $m$ for the Lebesgue measure. This ball fills a $m(B)/m(C)$ fraction of the volume of $C$, so there is a fraction $1-m(B)/m(C)$ of the volume of $C$ left to fill. Since $C\sm \overline B$ is open, after subtraction of negligible set, it can be covered by countably many dyadic cubes as above. Since each of these cubes is a rescaling of $C$, fit in each of them a rescaled version of $B$ so that the a fraction $m(B)/m(C)$ of the previous volume left to be filled has now been filled. Proceed by induction to construct, at each step $n$, an at most countable collection of disjoints open balls leaving a $(1-m(B)/m(C))^n$ fraction of the  volume of $C$ left to fill. The limiting collection is still at most countable, made of disjoint open balls included in $C$, and has Lebesgue measure $m(C)$, which establishes the result.
\end{Proof}

The next result is the key to establish Proposition \ref{PropQMD} in Section \ref{sec:LAN} for general densities $q_\ve\in H^1(\R^p)$ without further assumptions on the support of $q_\ve$, and crucially relies on the covering result in Lemma \ref{LemVitali}. 

\begin{Prop}\label{LemSupport}
    Let $f\in H^1(\R^p)$. Then the weak gradient $\nabla f$ of $f$ admits a version that vanishes on the set $\{f=0\}$, \ie $\nabla f = (\nabla f)\1[f>0]$.
\end{Prop}

\begin{Proof}{Lemma \ref{LemSupport}}
Throughout, denote by $m$ the Lebesgue measure on $\R^p$. Fix an open subset $\Om\subset \R^p$ such that $\{f>0\}\subset \Om$ and, by virtue of Lemma \ref{LemSupport}, let $(B_k)$ be an at most countable collection of disjoints balls such that $B_k\subset \Om$ and $m(\Om\sm \cup_{k\ge 1} B_k)=0$. Note that $\cup_{k\ge 1} B_k$ is an open subset, and fix $\psi\in C^\infty(\cup_{k\ge 1} B_k)$. Since the $B_k$'s are disjoint, then $\psi_{|B_k}\in C^\infty_c(B_k)$. By weak-differentiability, Fubini's theorem, and Green's formula, we have 
\begin{eqnarray*}
\lefteqn{
\int_{\R^p} (\nabla f) \psi\, dx 
=
- \int_{\R^p} f \nabla \psi\, dx
=
- \int_{\Om} f  \nabla \psi\, dx
}
\\[2mm]
&&
=
- \int_{\cup_k B_k} f \nabla \psi\, dx
=
- \sum_{k\ge 1} \int_{B_k} f \nabla \psi\, dx 
=
\sum_{k\ge 1} \int_{B_k} (\nabla f) \psi\, dx 
\\[2mm]
&&\hspace{7mm}
=
\int_{\cup_k B_k} (\nabla f) \psi\, dx 
=
\int_{\R^p} (\nabla f) \1_\Om \psi\, dx 
.
\end{eqnarray*}
Notice that the step $- \int_{B_k} f\nabla \psi\, dx = \int_{B_k} (\nabla f) \psi\, dx$ holds when $f\in C^1(\R^p)$ by Green's formula since $\psi$ has compact support in $B_k$, hence also when $f\in H^1(\R^p)$ by approximation in $H^1(\R^p)$ since this yields a simultaneous approximation of $f$ and $\nabla f$ also in $L^2(B_k)$. Since any function of the form $\psi\1_\Om$, where $\psi \in C^\infty_c(\R^p)$, can be approximated in $L^2(\R^p)$ by functions $\psi_j\in C^\infty_c(\cup_k B_k)$, we deduce that 
\begin{equation}\label{eq:IPPpsiOmega}
\int_{\R^p} (\nabla f) \psi\, dx 
=
\int_{\R^p} (\nabla f) \1_\Om \psi\, dx 
\end{equation}
holds for all $\psi\in C^\infty_c(\R^p)$. The outer regularity of the Lebesgue measure entails that there exists a sequence $(\Om_k)$ of open subsets such that $\{f>0\}\subset \Om_k$ and $m(\Om_k\sm \{f>0\})\to 0$. Taking limits in (\ref{eq:IPPpsiOmega}) applied to $\psi\in C^\infty_c(\R^d)$ fixed and $\Om_k$ entails that 
\begin{equation}
\int_{\R^p} f \nabla \psi\, dx
=
-\int_{\R^p} (\nabla f) \psi\, dx 
=
-\int_{\R^p} (\nabla f) \1[f>0] \psi\, dx
,\spa 
\forall\ \psi\in C^\infty_c(\R^p) 
.
\end{equation}
We deduce that $(\nabla f) \1[f>0]\in L^2(\R^p)$ is a version of $\nabla f$, which concludes the proof.
\end{Proof}

\subsection{Information operators on normed vector spaces}\label{sec:InvertibilityBanach}

Let $(V,\|\cdot\|_V)$ be a normed vector space and $(B,\ps{\cdot}{\cdot}_B)$ a Hilbert space, and assume that we have a linear continuous and injective map 
$$
A:(V, \|\cdot\|_V)\to (B, \ps{\cdot}{\cdot}_B)
.
$$
In the setting of statistical model (\ref{modelDensity}), the operator $A$ should be thought of as $\II_\ve^{1/2} \I_{\theta_0}$ (see Section \ref{sec:LAN}).
Since $A$ is continuous on $(V,\|\cdot\|_V)$, it admits a unique extension $A: (\overline V,\|\cdot\|_{\overline V})\to B$ to the $\|\cdot\|_V$-completion $\overline V$ of $V$; see Section \ref{sec:Completion} for an explicit construction and basic properties. Since this extension defines a continuous and linear operator between Banach spaces, $A$ admits an adjoint operator
$$
A^* : B \to \overline V^* 
,
$$
characterized by 
\begin{equation}\label{eqA:AdjointA}
(A^*[b])(v) 
\equiv 
\ps{b}{A[v]}_B
,\spa 
v\in \overline V,\ b\in B 
.
\end{equation}
The goal of this subsection is to characterize the mapping properties of the operator 
$$
A^* A : \overline V \to \overline V^*
.
$$

\subsubsection{Extension to the natural domain}

For this purpose, we introduce the `natural' norm on $V$ arising from $A$: since $A$ is injective, we can endow $V$ with the norm induced by $A$ through
$$
\|v\|_{V_A}
\equiv 
\|A[v]\|_{B}
,\spa 
v\in V 
.
$$
The resulting map $\|\cdot\|_{V_A}:V\to [0,\infty)$ is indeed a norm on $V$ since $A$ is linear and injective. Note that we have 
$$
\|v\|_{V_A} 
= 
\|A[v]\|_B 
\lesssim 
\|v\|_V
,\spa 
v\in V
,
$$
by continuity  of $A:(V,\|\cdot\|_V)\to B$. 
Denoting by $\overline V_A$ the completion of $V$ with respect to the $\|\cdot\|_{V_A}$-norm, we then have the continuous embedding 
\begin{equation}\label{VtoVAEmbed}
(\overline V,\|\cdot\|_{\overline V})
\embed 
(\overline V_A, \|\cdot\|_{\overline V_A})
.
\end{equation}
Notice that $\|v\|_{\overline V} = \|v\|_V$ and $\|v\|_{\overline V_A} = \|v\|_{V_A}$ when $v\in V$ (see Section \ref{sec:Completion}). Now, the space $V$ is dense in $\overline V$ and $\overline V_A$ by construction, and the (densely defined) linear map $A:V\to B$ satisfies 
$$
\|A[v]\|_B
=
\|v\|_{\overline V_A} 
,\spa 
v\in V 
.
$$
Since $B$ is a Banach space and $V$ is dense in $\overline V_A$, then $A$ admits a unique extension 
\begin{equation}\label{eqA:AExtended}
A:(\overline V_A, \|\cdot\|_{\overline V_A}) \to (B,\|\cdot\|_{B})
\end{equation}
as a linear and continuous operator satisfying
$$
\|A[v]\|_B
=
\|v\|_{\overline V_A}
,\spa 
v\in \overline V_A
,
$$
where, for $v\in \overline V_A\sm V$, $A[v]$ is defined as the (unique) limit of an arbitrary sequence $(v_k)\subset V$ such that $\|v_k-v\|_{\overline V_A}\to 0$. In particular, $A : \overline V_A \to B$ is an \emph{injective} linear and continuous operator between Banach spaces, and $A^* A$ also extends to a continuous operator 
\begin{equation}\label{eqA:A*Aextended}
A^* A : \overline V_A \to \overline V^* 
.
\end{equation}
Note that $\overline V_A$ is in fact a Hilbert space: since $A:\overline V_A\to B$ is injective, then 
$$
\ps{u}{v}_{\overline V_A}
\equiv 
\ps{A[u]}{A[v]}_B 
,\spa 
u,v\in \overline V_A 
,
$$
defines an inner-product on $\overline V_A$ such that $\ps{u}{u}_{\overline V_A}=\|u\|^2_{V_A}$. Remark that this does not necessarily make $\overline V$ a Hilbert space.









\subsubsection{Characterization of the range and surjectivity}

For all $b\in B$ and $v\in V$, observe that
\begin{equation}\label{eqA:ContA*toS}
|A^*[b](v)|
=
\lvert\ps{b}{A[v]}_B\rvert 
\le
\|b\|_B \|A[v]\|_{B} 
=
\|b\|_B \|v\|_{V_A}
.
\end{equation}
In particular, $A^*[b]$ belongs to the subset
\begin{equation}\label{eqA:DefS}
S
\equiv 
\Big\{ \ell\in \overline V^* : \|\ell\|_S\equiv \sup_{\substack{v\in V\\ \|v\|_{V_A} \le 1}} |\ell(v)| < \infty \Big\} 
,
\end{equation}
for all $b\in B$, which we informally call the $V$-dual of $\overline V_A$. Consequently, we have the following inclusions of sets 
\begin{equation}\label{eqA:InclusionSets}
A^* A(\overline V_A) 
\subset 
A^*(B)
\subset 
S
.
\end{equation}
The following result establishes that the previous inclusions are in fact equalities of sets.

\begin{Prop}\label{PropEqualityRangesA}
We have 
$
A^*(B) 
= 
A^* A(\overline V_A) 
=
S
.
$
\end{Prop}

\begin{Proof}{Proposition \ref{PropEqualityRangesA}}
In view of (\ref{eqA:InclusionSets}), it is enough to establish that $S\subset A^* A(\overline V_A)$. Then, fix $\ell\in S$, \ie $\ell:\overline V\to \R$ is linear and continuous, and 
$$
|\ell(v)| 
\le 
\|\ell\|_S \|v\|_{V_A}
,\spa 
v\in V
.
$$
Since $V$ is dense in $\overline V_A$, then $\ell$ extends to a linear and continuous map $\ell:(\overline V_A, \|\cdot\|_{\overline V_A})\to \R$, \ie $\ell\in \overline V_A^*$. Since $\overline V_A$ is a Hilbert space, Riesz representation theorem entails that there exists $\bar \ell\in \overline V_A$ such that 
$$
\ell(v) 
=
\ps{\bar \ell}{v}_{\overline V_A}
=
\ps{A[\bar \ell]}{A[v]}_{B}
=
(A^* A[\bar \ell])(v)
,\spa 
v\in \overline V_A 
,
$$
where we used the adjoint property $(A^*[b])(v)=\ps{b}{A[v]}_B$ from (\ref{eqA:AdjointA}) with $b=A[\bar\ell]\in B$ in view of (\ref{eqA:AExtended}). In particular, $\ell$ and $A^* A[\bar\ell]$ coincide on $\overline V$, hence $\ell = A^* A[\bar \ell]$ in $\overline V^*$. We deduce that $\ell \in A^* A[\overline V_A]$. Together with (\ref{eqA:InclusionSets}), this yields the result. 
\end{Proof}
\vspace{3mm}
The following result further establishes that $(S,\|\cdot\|_S)$ is a Banach space.

\begin{Prop}\label{PropBanachA}
    The space $(S,\|\cdot\|_S)$ from (\ref{eqA:DefS}) is a Banach space. 
\end{Prop}

\begin{Proof}{Proposition \ref{PropBanach}}
    The fact that $S$ is a linear space follows immediately from the definition (\ref{eqA:DefS}) and the triangle inequality. The homogeneity and triangle inequality for $\|\cdot\|_S$ are obvious, and if $w\in S$ satisfies $\|w\|_S=0$, then $w(v)=0$ for all $v\in V$ with $\|A[v]\|_B\le 1$; since $w$ is linear, this implies that $w(v)=0$ for all $v\in V$, hence $w(v)=0$ for all $v\in \overline V$, to the effect that $w=0$ in $\overline{V}^*$. We now establish that $S$ is complete by showing that it is in linear isometric correspondence with the Banach space $\overline V_A^*$. Fix $w\in S$. Then $w\in \overline V^*$ and for all $v\in V$, we have
    $$
    |w(v)| 
    \le 
    \|w\|_S \|v\|_{\overline V_A} 
    .
    $$
    Since $V$ is dense in $\overline V_A$, then $w$ extends uniquely to a continuous and linear map $\overline w:\overline V_A\to \R$ such that $\|\overline w\|_{\overline V_A^*} \le \|w\|_S$. Recalling the embedding $\overline V\embed \overline V_A$ from (\ref{VtoVAEmbed}), we further have 
    $$
    \|\overline w\|_{\overline V_A^*}
    =
    \sup_{\substack{v\in \overline V_A\\ \|v\|_{\overline V_A} \le 1}} |\overline w(v)|
    \ge 
    \sup_{\substack{v\in \overline V\\ \|v\|_{\overline V_A} \le 1}} |\overline w(v)|
    =
    \sup_{\substack{v\in \overline V\\ \|v\|_{\overline V_A} \le 1}} |w(v)|
    \ge 
    \|w\|_S 
    .
    $$
    We deduce that $w$ extends uniquely to an element $\overline w\in \overline V_A^*$ with $\|w\|_S=\|\overline w\|_{\overline V_A^*}$. Now fix an arbitrary element $\overline w\in \overline V_A^*$, \ie $\overline w:\overline V_A\to \R$ is a continuous and linear map. We then have, for all $v\in V$
    $$
    |\overline w(v)| 
    \lesssim 
    \|v\|_{\overline V_A}
    \lesssim
    \|v\|_{\overline V}
    ,
    $$
    by the embedding $\overline V\embed \overline V_A^*$ from (\ref{VtoVAEmbed}). In particular, the restriction $w:\overline V\to \R$ of $\overline w$ to $\overline V$ defines an element of $\overline V^*$ such that 
    $$
    \|w\|_S 
    =
    \sup_{\substack{v\in V\\ \|v\|_{\overline V_A} \le 1}} |w(v)|
    =
    \sup_{\substack{v\in V\\ \|v\|_{\overline V_A} \le 1}} |\overline w(v)|
    =
    \sup_{\substack{v\in \overline V_A\\ \|v\|_{\overline V_A} \le 1}} |\overline w(v)|
    =
    \|\overline w\|_{\overline V_A}
    <
    \infty 
    ,
    $$
    where we used the fact that $V$ is dense in $\overline V_A$. We deduce that the map 
    $$
    \phi : 
    (S,\|\cdot\|_S)\to (\overline V_A^*, \|\cdot\|_{\overline V_A^*}),\ w\to \overline w
    $$
    defines a linear and isometric bijection. Since $\overline V_A^*$ is a Banach space, as the topological dual of the normed vector space $(V,\|\cdot\|_{V_A})$---see, e.g., Proposition 6.1.9 in \citet*{BogRealandFunc}---we deduce that $(S,\|\cdot\|_S)$ is a Banach space.
\end{Proof}
\vspace{3mm}

Although $(S,\|\cdot\|_S)$ is a Banach space, and the ranges of $A^*:B\to \overline V^*$ and $A^* A : \overline V_A\to \overline V^*$ are equal to $S$, this does not necessarily mean that $A^*$ and $A^* A$ have closed range as $S$ need not be a closed subset of $(\overline V^*, \|\cdot\|_{\overline V^*})$.

\subsubsection{Injectivity and invertibility}

We have seen that the range of the continuous linear operator $A^* : B\to \overline V^*$ is equal to a Banach space $(S,\|\cdot\|_S)$, given by (\ref{eqA:DefS}), and (\ref{eqA:ContA*toS}) establishes that 
$$
A^* : (B,\|\cdot\|_B)\to (S,\|\cdot\|_S)
$$ 
is a surjective and continuous mapping between Banach spaces. In particular, 
$$
A^* A : (\overline V_A,\|\cdot\|_{\overline V_A}) \to (S,\|\cdot\|_S)
$$
is also a surjective and continuous mapping between Banach spaces.

\begin{Theor}\label{TheorI*IHomeoA}
    The linear operator $A^* A : (\overline V_A,\|\cdot\|_{\overline V_A}) \to (S,\|\cdot\|_S)$ is a homeomorphism.
\end{Theor}

\begin{Proof}{Theorem \ref{TheorI*IHomeoA}}
    We already established that $A^* A : (\overline V_A,\|\cdot\|_{\overline V_A}) \to (S,\|\cdot\|_S)$ defines a surjective linear and continuous operator. Since $\overline V_A$ and $S$ are Banach spaces, it follows from the inverse mapping theorem that $A^* A$ is a homeomorphism if and only if it is a continuous bijection. Since it is surjective, it remains to establish injectivity. Since $A^*(B)\subset S$, we have 
    $$
    |(A^*[b])(v)|
    \le 
    \|A^*[b]\|_S \|v\|_{\overline V_A}
    ,\spa 
    b\in B,\ v\in V 
    .
    $$
    In particular, for all $v\in V$ and $b\equiv A[v]\in B$, we have
    $$
    \|v\|^2_{\overline V_A}
    =
    \ps{A[v]}{A[v]}_B 
    =
    \ps{b}{A[v]}_B
    =
    (A^*[b])(v)
    \leq 
    \|A^* b\|_{S} \|v\|_{\overline V_A}
    .
    $$
    We deduce that 
    $$
    \|v\|_{\overline V_A}
    \le 
    \|A^* A[v]\|_S 
    ,\spa 
    v\in V
    .
    $$
    This inequality can be extended to all of $\overline V_A$ by continuity since $V$ is dense in $\overline V_A$. The previous display yields injectivity of $A^* A : \overline V_A\to S$, hence concludes the proof.
\end{Proof}
\vspace{3mm}

In particular, $S$ can be endowed with a Hilbert space structure through 
$$
\ps{\ell}{\ell'}_{S}
\equiv 
\ell((A^* A)^{-1}[\ell'])
=
\ps{A(A^* A)^{-1}[\ell]}{A(A^* A)^{-1}[\ell']}_B
,\spa 
\ell, \ell'\in S  
.
$$

\subsubsection{Constructive completions}\label{sec:Completion}

One can construct the unique (up to isometry) closures $\overline V$ and $\overline V_A$ of $V$ with respect to the $\|\cdot\|_V$-norm and the $\|\cdot\|_{V_A}$-norm, respectively. A possible construction goes as follows, illustrated on the $\|\cdot\|_V$-norm. Consider the topological dual $V^*$ of the normed space $(V,\|\cdot\|_V)$, equipped with the operator norm 
$$
\|\ell\|_{V^*}
\equiv 
\sup_{\substack{v\in V\\ \|v\|_V\le 1}} |\ell(v)| 
,\spa 
\ell\in V^*
,
$$
and the second dual $V^{**}$ obtained as the dual of the normed space $(V^*, \|\cdot\|_{V^*})$, with corresponding norm
$$
\|j\|_{V^{**}} 
\equiv 
\sup_{\substack{\ell\in V^*\\ \|\ell\|_{V^*} \le 1}} |j(\ell)|
,\spa 
j\in V^{**}
.
$$
Then, consider the linear map
$
\JJ : V \mapsto V^{**}
,\ v\mapsto \JJ_v
,
$
where $\JJ_v(\ell)\equiv \ell(v)$ for all $v\in V$ and $\ell\in V^*$. Then, $V^{**}$ is a Banach space and $\JJ$ is a linear and isometric embedding of $V$ into the Banach space $V^{**}$; see, e.g., Proposition 6.4.3 in \citet*{BogRealandFunc}. Therefore, the closure $(\overline V, \|\cdot\|_{\overline V})$ of $(V,\|\cdot\|_V)$ can be obtained as the closure of $\JJ(V)$ in $V^{**}$. Similarly, for $\overline V_A$ we construct a linear isometry $\JJ_A: (V,\|\cdot\|_{V_A}) \mapsto V_A^{**}$ between $V_A\equiv (V,\|\cdot\|_{V_A})$ and its second dual. In particular, we have
$$
\|v\|_{\overline V_A} 
=
\|\JJ_A(v)\|_{V_A^{**}}
=
\sup_{\substack{\ell\in V_A^*\\ \|\ell\|_{\rm op}\le 1}} |\ell(v)|
=
\|v\|_{V_A}
=
\|A[v]\|_{B}
,
\spa 
v\in V
.
$$
Hence, the norm in the closed space $\overline V_A$, when restricted to the original space $V$, coincides with the norm $\|\cdot\|_{V_A}$ on $V$, and the corresponding result holds for the norm $\|\cdot\|_{\overline V}$ as well. Consequently, we have 
$
\|v\|_{\overline V_A}
=
\|v\|_{V_A} 
\lesssim 
\|v\|_V 
=
\|v\|_{\overline V}
$ for all $v\in V$. This yields the continuous embedding
$
(\overline V, \|\cdot\|_{\overline V})
\embed 
(\overline V_A, \|\cdot\|_{\overline V_A})
.
$

\section*{Acknowledgements}
The author acknowledges funding from an ERC Advanced Grant (UKRI G116786).

\bibliographystyle{agsm}
\bibliography{mybib}

\bigskip\bigskip

\textsc{Department of Pure Mathematics \& Mathematical Statistics}

\textsc{University of Cambridge}, Cambridge, UK

Email: dk738@cam.ac.uk

\end{document}